\documentclass[11pt]{article}
\usepackage[tbtags]{amsmath}
\usepackage{amssymb}
\usepackage{amsthm}
\usepackage[misc]{ifsym}
\usepackage{cases}
\numberwithin{equation}{section}   %数学公式，equation环境中将按照section排序
\normalsize
\title{\bf Stochastic Linear Quadratic Stackelberg Differential Game with Overlapping Information
\thanks{Shi acknowledges the financial support from the National Natural Science Fund of China (Grant No. 11571205). Wang acknowledges the financial support from the National Natural Science Funds of China (Grant No. 11371228, 61422305, 61633015), and the Natural Science Fund for Distinguished Young Scholars of Shandong Province of China (Grant No. JQ201418). Xiong acknowledges the financial support from the Southern University of Science and Technology Start-Up Fund (Grant No. Y01286220). Part of the content of this paper was presented on the 36th Chinese Control Conference, Dalian, P. R. China, July 26-28, 2017.}}
\author{\normalsize Jingtao Shi,\thanks{\it School of Mathematics, Shandong University, Jinan 250100, P. R. China E-mail: shijingtao@sdu.edu.cn}\quad
Guangchen Wang\thanks{\it School of Control Science and Engineering, Shandong University, Jinan 250061, P. R. China E-mail: wguangchen@sdu.edu.cn}\quad
Jie Xiong\thanks{\it Department of Mathematics, Southern University of Science and Technology, Shenzhen 518055, P. R. China E-mail: xiongj@sustc.edu.cn}}
\date{}
\begin{document}
\maketitle \noindent{\bf Abstract:}\quad This paper is concerned with the stochastic linear quadratic Stackelberg differential game with overlapping information, where the diffusion terms contain the control and state variables. Here the term ``overlapping" means that there are common part between the follower's and the leader's information, while they have no inclusion relation. Optimal controls of the follower and the leader are obtained by the stochastic maximum principle, the direct calculation of the derivative of the cost functional and stochastic filtering. A new system of Riccati equations is introduced to represent the state estimate feedback of the Stackelberg equilibrium strategy. A special solvable case is then studied and is applied to the continuous-time principal-agent problem.

\vspace{2mm}


\vspace{2mm}


\section{Introduction}

Throughout this paper, we denote by $\mathbb{R}^n$ the Euclidean space of $n$-dimensional vectors, by $\mathbb{R}^{n\times d}$ the space
of $n\times d$ matrices, by $\mathcal{S}^n$ the space of $n\times n$ symmetric matrices. $\langle\cdot,\cdot\rangle$ and $|\cdot|$ denote the scalar product
and norm in the Euclidean space, respectively. $\top$ appearing in the superscripts denotes the transpose of a matrix. $f_x,f_{xx}$ denote the partial derivative and twice partial derivative with respect to $x$ for a differentiable function $f$, respectively.

\subsection{Motivation}

First, we present the following example which motivates us to study the problem in this paper.

{\it Example 1.1:} (Continuous time principal-agent problem) The principal contracts with the agent to manage a production process, whose output $Y(\cdot)$ evolves as
\begin{equation}\label{output equation}
\left\{
\begin{aligned}
dY(t)&=Be(t)dt+\sigma_1dW_1(t)+\sigma_2dW_2(t)+\sigma_3dW_3(t),\ t\in[0,T],\\
 Y(0)&=Y_0\in\mathbb{R},
\end{aligned}
\right.
\end{equation}
where $e(\cdot)\in A\subset\mathbb{R}$ is the agent's effort choice, $B$ represents the productivity of effort, and there are three additive shocks (due to the three independent Brownian motions $W_1(\cdot),W_2(\cdot),W_3(\cdot)$) to the output. The output of the production adds to the principal's asset $y(\cdot)$, which earns a risk free return $r$, and out of which he pays the agent $s(\cdot)\in S\subset\mathbb{R}$ and withdraws his own consumption $d(\cdot)\in\mathbb{R}$. Thus the principal's asset evolves as
\begin{equation}\label{asset equation of the principal}
\left\{
\begin{aligned}
dy(t)&=\big[ry(t)+Be(t)-s(t)-d(t)\big]dt+\sigma_1dW_1(t)+\sigma_2dW_2(t)+\sigma_3dW_3(t),\ t\in[0,T],\\
 y(0)&=y_0\in\mathbb{R},
\end{aligned}
\right.
\end{equation}
where $y_0$ is the initial asset. In addition, the agent has his own wealth $m(\cdot)$, out of which he consumes $c(\cdot)$, then
\begin{equation}\label{asset equation of the agent}
\left\{
\begin{aligned}
dm(t)&=\big[rm(t)+s(t)-c(t)\big]dt+\bar{\sigma}_1dW_1(t)+\bar{\sigma}_2dW_2(t)+\bar{\sigma}_3dW_3(t),\ t\in[0,T],\\
 m(0)&=m_0\in\mathbb{R},
\end{aligned}
\right.
\end{equation}
The agent earns the same rate of return $r$ on his savings, gets income flows due to his payment $s(\cdot)$, and draws down wealth to consume. In the above $\sigma_i,\bar{\sigma}_i,i=1,2,3$ are all constants. At the terminal time $T$, the principal makes a final payment $s(T)$ and the agent chooses consumption based on this payment and his terminal wealth $m(T)$.

We consider an optimal implementable contract problem in the so-called ``hidden savings" information structure (Williams \cite{Wil15}). In this problem, the principal can observe his asset $y(\cdot)$ and the agent's initial wealth $m_0$, but cannot monitor the agent's effort $e(\cdot)$, consumption $c(\cdot)$ and wealth $m(\cdot)$. The principal must provide incentives for the agent to put forth the desired amount of the effort. For any $s(\cdot),d(\cdot)$, the agent first chooses his effort $e^*(\cdot)$ and consumption $c^*(\cdot)$ such that his preference
\begin{equation}\label{cost of the principal}
 J_1\big(e(\cdot),c(\cdot),s(\cdot),d(\cdot)\big)=\frac{1}{2}\mathbb{E}\left[\int_0^T\big[c^2(t)-e^2(t)+m^2(t)\big]dt+m^2(T)\right]
\end{equation}
is maximized. The above $(e^*(\cdot),c^*(\cdot))$ is called an implementable contract if it meets the recommended actions of the principal's, which is based on the principal's observable wealth $y(\cdot)$. Then, the principal selects his payment $s^*(\cdot)$ and consumption $d^*(\cdot)$, to maximize his preference
\begin{equation}\label{cost of the agent}
 J_2\big(e^*(\cdot),c^*(\cdot),s(\cdot),d(\cdot)\big)=\frac{1}{2}\mathbb{E}\left[\int_0^T\big[d^2(t)-s^2(t)+y^2(t)\big]dt+y^2(T)\right].
\end{equation}
Noting that in \cite{Wil15}, exponential preferences are introduced while here we consider the quadratic case. For $t>0$, let
\begin{equation*}
\mathcal{F}_t\triangleq\sigma\big\{W_1(s),W_2(s),W_3(s),0\leq s\leq t\big\}
\end{equation*}
which contains all the information up to time $t$. Let
\begin{equation*}
\mathcal{G}^1_t\triangleq\sigma\big\{W_1(s),W_3(s),0\leq s\leq t\big\}
\end{equation*}
contains the information available to the agent, and
\begin{equation*}
\mathcal{G}^2_t\triangleq\sigma\big\{W_2(s),W_3(s);0\leq s\leq t\big\}
\end{equation*}
contains the information available to the principal, up to time $t$ respectively. Obviously, the information available to them at time $t$ are asymmetric while possess the overlapping part. In the problem, for any $s(\cdot),d(\cdot)$, first the agent solves the following optimization problem:
\begin{equation}\label{cost functional agent}
       J_1\big(e^*(\cdot),c^*(\cdot),s(\cdot),d(\cdot)\big)=\max\limits_{e,c}J_1\big(e(\cdot),c(\cdot),s(\cdot),d(\cdot)\big),
\end{equation}
where $(e^*(\cdot),c^*(\cdot))$ is a $\mathcal{G}^1_t$-adapted process pair. Then the principal solves the following optimization problem:
\begin{equation}\label{cost functional principal}
       J_2\big(e^*(\cdot),c^*(\cdot),s^*(\cdot),d^*(\cdot)\big)=\max\limits_{s,d}J_2\big(e^*(\cdot),c^*(\cdot),s(\cdot),d(\cdot)\big),
\end{equation}
where $(s^*(\cdot),d^*(\cdot))$ is a $\mathcal{G}^2_t$-adapted process pair. This formulates a {\it stochastic linear quadratic (LQ) Stackelberg differential game with overlapping information}. In this setting, the agent is the follower and the principal is the leader. Any process quadruple $(e^*(\cdot),c^*(\cdot),s^*(\cdot),d^*(\cdot))$ satisfying the above two equalities is called a {\it Stackelberg equilibrium strategy}. In \cite{Wil15}, a solvable continuous time principal-agent model is considered under three information structures (full information, hidden actions and hidden savings) and the corresponding optimal contract problems are solved explicitly. But it can not cover our model. For more information for the principal-agent problem, please refer to the monograph by Cvitani\'{c} and Zhang \cite{CZ13}.

Other examples which motivated us to study the problem in this paper can be found in the insider trading model (\O ksendal \cite{Ok06}), the cooperative advertising and pricing problem (He et al. \cite{HPS09}), the continuous time manufacturer-newsvendor problem (\O ksendal et al. \cite{OSU13}), the LQ Nash differential game with asymmetric information (Chang and Xiao \cite{CX14}), and optimal reinsurance arrangement problem between the insurer and the reinsurer (Chen and Shen \cite{CS18}), etc. We will not give their detail statement for the space limitation.

\subsection{Problem formulation}

Inspired by the examples above, we study the stochastic LQ Stackelberg differential game with overlapping information in this paper.

Let $(\Omega,\mathcal{F},\mathbb{P})$ be a complete probability space, on which a standard three-dimensional Brownian motion $\{W_1(t),W_2(t),W_3(t)\}_{0\leq t\leq T}$ is defined, where $T>0$ is a finite time duration. Let $\{\mathcal{F}_t\}_{0\leq t\leq T}$ be the natural filtration generated by $(W_1(\cdot),W_2(\cdot),W_3(\cdot))$ which satisfies the usual conditions and $\mathcal{F}_T=\mathcal{F}$.

We consider an $\mathbb{R}^n$-valued state process $x^{u_1,u_2}(\cdot)$ which satisfies the linear {\it stochastic differential equation} (SDE)
\begin{equation}\label{state equation}
\left\{
\begin{aligned}
     dx^{u_1,u_2}(t)&=\big[A_0(t)x^{u_1,u_2}(t)+B_0(t)u_1(t)+C_0(t)u_2(t)\big]dt\\
                    &\quad+\sum\limits_{i=1}^3\big[A_i(t)x^{u_1,u_2}(t)+B_i(t)u_1(t)+C_i(t)u_2(t)\big]dW_i(t),\ t\in[0,T],\\
      x^{u_1,u_2}(0)&=x_0.
\end{aligned}
\right.
\end{equation}
Here $u_1(\cdot)$ is the follower's control process and $u_2(\cdot)$ is the leader's control process, which are $\mathbb{R}^{k_1}$ and $\mathbb{R}^{k_2}$-valued, respectively. For $i=0,1,2,3$, $A_i(\cdot)\in\mathbb{R}^{n\times{d_i}}$, $B_i(\cdot)\in\mathbb{R}^{n\times{k_1}}$ and $C_i(\cdot)\in\mathbb{R}^{n\times{k_2}}$ are all matrix-valued processes and $x_0\in\mathbb{R}^n$. We define the admissible control sets of the follower and the leader, as follows.
\begin{equation}\label{admissible control sets}
\begin{aligned}
\mathcal{U}_i:=\Big\{u_i(\cdot)\big|u_i(\cdot):\Omega\times[0,T]\rightarrow\mathbb{R}^{k_i}\mbox{ is }\mathcal{G}^i_t\mbox{-adapted and }\sup\limits_{0\leq t\leq T}\mathbb{E}|u_i(t)|^2<\infty\Big\},\ i=1,2.
\end{aligned}
\end{equation}
Here $\mathcal{G}^i_t\triangleq\sigma\{W_i(s),W_3(s);0\leq s\leq t\},i=1,2$ denotes the information of the follower and the leader, respectively.

We now formulate the problem by the following two steps. In step 1, the follower choose a $u_1^*(\cdot)\in\mathcal{U}_1$, which depends on the control $u_2(\cdot)$ of the leader, to minimize the cost functional
\begin{equation}\label{cost functional-follower}
\begin{aligned}
 J_1(u_1(\cdot),u_2(\cdot))&=\frac{1}{2}\mathbb{E}\bigg[\int_0^T\Big(\big\langle Q_1(t)x^{u_1,u_2}(t),x^{u_1,u_2}(t)\big\rangle+\big\langle N_1(t)u_1(t),u_1(t)\big\rangle\Big)dt\\
             &\qquad\qquad+\big\langle G_1x^{u_1,u_2}(T),x^{u_1,u_2}(T)\big\rangle\bigg].
\end{aligned}
\end{equation}
Here $Q_1(\cdot)\in\mathbb{R}^{n\times n},N_1(\cdot)\in\mathbb{R}^{{k_1}\times{k_1}}$ are nonnegative matrices-valued processes and $G_1$ is a nonnegative $\mathbb{R}^{n\times n}$-valued matrix. In step 2, the leader takes into account the follower's optimal control $u_1^*(\cdot)$ in his cost functional, and selects an optimal control $u_2^*(\cdot)\in\mathcal{U}_2$ which will minimize
\begin{equation}\label{cost functional-leader}
\begin{aligned}
J_2(u_1^*(\cdot),u_2(\cdot))&=\frac{1}{2}\mathbb{E}\bigg[\int_0^T\Big[\big\langle Q_2(t)x^{u_1^*,u_2}(t),x^{u_1^*,u_2}(t)\big\rangle+\big\langle N_2(t)u_2(t),u_2(t)\big\rangle\Big]dt\\
              &\qquad\qquad+\big\langle G_2x^{u_1^*,u_2}(T),x^{u_1^*,u_2}(T)\big\rangle\bigg],
\end{aligned}
\end{equation}
where $x^{u_1^*,u_2}(\cdot)$ denotes the optimal state of the follower which is the solution to (\ref{state equation}) with respect to $u_1^*(\cdot)$. $Q_2(\cdot)\in\mathbb{R}^{n\times n}$ is a nonnegative matrix-valued process, $N_2(\cdot)\in\mathbb{R}^{{k_2}\times{k_2}}$ is a positive matrix-valued process and $G_2$ is a nonnegative $\mathbb{R}^{n\times n}$-valued matrix. In this general LQ model, the information of the leader and the follower have overlapping part, due to the structure of the admissible control sets. The target of this paper is to give the conditions of its Stackelberg equilibrium strategy $(u_1^*(\cdot),u_2^*(\cdot))\in\mathcal{U}_1\times\mathcal{U}_2$.

Note that some Stackelberg differential games with partially observable information can be put into the above LQ model by the Girsanov transformation. See, for example, Shi et al. \cite{SWX16}, Wang et al. \cite{WXX16}.

\subsection{Literature review and the contribution of this paper}

In recent years, Stackelberg (also known as leader-follower) game has been an active topic, in the research of nonzero-sum games. Compared with its Nash counterpart, Stackelberg game has many appealing properties, which are useful both in theory and applications. The Stackelberg solution to the game is obtained when one of the players is forced to wait until the other player announces his decision, before making his own decision. Problems of this nature arise frequently in economics, where decisions must be made by two parties and one of them is subordinated to the other, and hence must wait for the other party's decision before formulating its own. The research of Stackelberg game can be traced back to the pioneering work by Stackelberg \cite{S52} in static competitive economics. Simann and Cruz \cite{SC73} studied the dynamic LQ Stackelberg differential game, and the Stackelberg strategy was expressed in terms of Riccati-like differential equations. Bagchi and Basar \cite{BB81} investigated the stochastic LQ Stackelberg differential game, where the diffusion term of the Ito-type state equation does not contain the state and control variables. Existence and uniqueness of its Stackelberg solution are established, and the leader's optimal strategy is solved as a nonstandard stochastic control problem and is shown to satisfy a particular integral equation. Yong \cite{Yong02} extended the stochastic LQ Stackelberg differential game to a rather general framework, where the coefficients could be random matrices, the control variables could enter the diffusion term of the state equation and the weight matrices for the controls in the cost functionals need not to be positive definite. The problem of the leader is first described as a stochastic control problem of a {\it forward-backward stochastic differential equation} (FBSDE). Moreover, it is shown that the open-loop solution admits a state feedback representation if a new stochastic Riccati equation is solvable. \O ksendal et al. \cite{OSU13} proved a maximum principle for the Stackelberg differential game when the noise is described as an Ito-L\'{e}vy process, and found applications to a continuous time manufacturer-newsvendor model. Bensoussan et al. \cite{BCS15} proposed several solution concepts in terms of the players' information sets, for the stochastic Stackelberg differential game with the control-independent diffusion term, and derived the maximum principle under the adapted closed-loop memoryless information structure. Xu and Zhang \cite{XZ16} studied both discrete- and continuous-time stochastic Stackelberg differential games with time delay. By introducing a new costate, a necessary and sufficient condition for the existence and uniqueness of the Stackelberg equilibrium was presented and was designed in terms of three decoupled and symmetric Riccati equations. Some recent progress about Stackelberg games can be seen in a review paper by Li and Sethi \cite{LS17} and the references therein.

However, the above literatures do not consider the feature of asymmetric information in Stackelberg differential game, which we believe, to our best knowledge, that it is a nature and important feature from the point of view of theory and applications. In fact, there are some literatures about asymmetric information game theory. For example, \cite{SC73} considered a non-zero sum velocity-controlled pursuit-evasion game, where the pursuer's information is always later in time than that of the evader's, which is in some sense of time asymmetry. \O ksendal \cite{Ok06} solved a universal optimal consumption rate problems with insider trading, where the consumer is called an insider when he has more information than what can be obtained by observing the driving process. That is a kind of information asymmetry with respect to the driving process. Cardaliaguet and Rainer \cite{CR09} investigated a two-player zero-sum stochastic differential game in which the players have an asymmetric information on the random payoff. Lempa and Matom\"{a}ki \cite{LM13} studied a Dynkin game with asymmetric information. The players have asymmetric information on the random expiry time, namely only one of the players is able to observe its occurrence. Chang and Xiao \cite{CX14} studied an LQ nonzero sum differential game problem with asymmetric information, where different $\sigma$-algebra generated by different Brownian motions are introduced to represent the asymmetric information of the two players. Nash equilibrium points are obtained for several classes of asymmetric information by stochastic maximum principle and technique of completion of squares. Shi et el. \cite{SWX16} solved a stochastic leader-follower differential game with asymmetric information, where the information available to the follower is based on some sub-$\sigma$-algebra of that available to the leader. Stochastic maximum principles and verification theorems with partial information were obtained. An LQ stochastic leader-follower differential game with noisy observation was solved via measure transformation, stochastic filtering, where not all the diffusion coefficients contain the state and control variables. In a companion paper by Shi et al. \cite{SWX17}, an LQ stochastic Stackelberg differential game with asymmetric information was researched, where the control variables enter both diffusion coefficients of the state equation, via some {\it forward-backward stochastic differential filtering equations} (FBSDFEs). Shi and Wang \cite{SW16} considered another kind of LQ leader-follower stochastic differential game, where the information available to the leader is a sub-$\sigma$-algebra of the filtration generated by the underlying Brownian motion. Wang et. al. \cite{WXX16} focused on an LQ non-zero sum differential game problem derived by the BSDE with asymmetric information. Three classes of observable filtrations are described to classify the information available to the two players. Using the filters of FBSDEs, feedback Nash equilibrium points with observable information generated by Brownian motions were obtained.

In this paper, we consider the stochastic LQ Stackelberg differential game with overlapping information. The LQ problems constitute an extremely important class of optimal control or differential game problems, since they can model many problems in applications, and also reasonably approximate nonlinear control or game problems (\cite{CLZ98}, \cite{LZ13}). The novelty of the formulation and the contribution in this paper is the following.

(i) In our framework, both information filtration available to the leader and the follower could be sub-$\sigma$-algebras of the complete information filtration naturally generated by the random noise source. Specifically, the system noise is described by three independent Brownian motions $W_1(t),W_2(t),W_3(t)$, from which the filtration generated denotes the complete information up to time $t$. The information of the follower comes from the filtration generated by $W_1(t),W_3(t)$, while the information of the leader comes from the filtration generated by $W_2(t),W_3(t)$. This framework is more suitable and interesting to illustrate some game problems in reality.

(ii) The general case that the diffusion terms contain the control and state variables is considered. As is well known in stochastic control and differential game theory, this brings us rather intrinsic mathematical difficulty and technical demanding, especially for the problem of the leader. We overcome the difficulty by the maximum principle approach, the direct calculation of the derivative of the cost functional, and stochastic filtering technique. A new system of high-dimensional Riccati equations is introduced to represent the state estimate feedback of the Stackelberg equilibrium strategy, though its general solvability is very difficult to verify.

(iii) A special solvable case when the diffusion terms are control independent is considered. In this case the system of Riccati equations can be proved to be solvable uniquely, which is used to represent the state feedback form of the Stackelberg equilibrium stratery.

(iv) A continuous-time principal-agent problem is solved by applying the theoretical results. The Stackelberg equilibrium strategy of the principal and the agent are represented explicitly.

We refer to Frankowska et al. \cite{FZZ17} and the references therein for more details on the recent progress for the study of maximum principles for stochastic systems.

The rest of this paper is organized as follows. In Section 2, the problem formulated in Section 1.2 are solved in the two subsections. In subsection 2.1, the follower's problem is considered, while the leader's problem is studied in Subcection 2.2. The Stackelberg equilibrium strategy is derived. A special case with control independent diffusion terms is completely solved in Section 3. In Section 4, the results in the previous sections is applied to a continuous-time principal-agent problem. Some concluding remarks are given in Section 5.

\section{Main results}

In this section, we will deal with the problems of the follower and the leader in two subsections, respectively. First, we introduce the following lemma, which belongs to Xiong \cite{Xiong08} and will play a fundemental role in this paper.

\vspace{1mm}

\noindent{\bf Lemma 2.1}\quad{\it Let $f(\cdot),g(\cdot)$ be $\mathcal{F}_t$-adapted processes, satisfying $\mathbb{E}\int_0^T|f(s)|ds+\mathbb{E}\int_0^T|g(s)|^2ds<\infty$. Then
\begin{equation}
\begin{aligned}
&\mathbb{E}\left[\int_0^tf(s)ds\bigg|\mathcal{G}^i_t\right]=\int_0^t\mathbb{E}\big[f(s)|\mathcal{G}^i_t\big]ds,\\
&\mathbb{E}\left[\int_0^tg(s)dW^i(s)\bigg|\mathcal{G}^i_t\right]=\int_0^t\mathbb{E}\big[g(s)|\mathcal{G}^i_t\big]dW^i(s),\ i=1,2,3,
\end{aligned}
\end{equation}
and for $i,j=1,2,3$,
\begin{equation}
\mathbb{E}\left[\int_0^tg(s)dW^i(s)\bigg|\mathcal{G}^j_t\right]=0,\ i\neq j.
\end{equation}}

\noindent{\it Proof.}\quad Please refer to Lemma 5.4 of \cite{Xiong08}. $\Box$

For any $\mathcal{F}_t$-adapted process $\xi(\cdot)$, we denote by $$\hat{\xi}(t):=\mathbb{E}[\xi(t)|\mathcal{G}^1_t],\ \check{\xi}(t)\triangleq\mathbb{E}[\xi(t)|\mathcal{G}^2_t]$$ and $$\check{\hat{\xi}}(t)\triangleq\mathbb{E}\big[\mathbb{E}[\xi(t)|\mathcal{G}^1_t]\big|\mathcal{G}^2_t\big]\equiv\mathbb{E}\big[\mathbb{E}[\xi(t)|\mathcal{G}^2_t]\big|\mathcal{G}^1_t\big]$$ its optimal filtering estimates.

\subsection{Problem of the follower}

In this subsection, we try to find the necessary condition for the optimal control of the follower. For given leader's control $u_2$, let us assume that there exists a $\mathcal{G}^1_t$-adapted optimal control $u_1^*(\cdot)$ of the follower, and the corresponding optimal state is $x^{u_1^*,u_2}(\cdot)$ as before. We define the follower's Hamiltonian function as
\begin{equation}\label{Hamiltonian function of the follower}
\begin{aligned}
&H_1\big(t,x,u_1,u_2,q,k_1,k_2,k_3\big)\triangleq\langle q,A_0(t)x+B_0(t)u_1+C_0(t)u_2\rangle\\
&\quad+\sum\limits_{i=1}^3\langle k_i,A_i(t)x+B_i(t)u_1+C_i(t)u_2\rangle-\frac{1}{2}Q_1(t)|x|^2-\frac{1}{2}N_1(t)|u_1|^2.
\end{aligned}
\end{equation}
The maximum principle (See, for example, \cite{SWX16}) yields that
\begin{equation}\label{optimal control of the follower}
N_1(t)u_1^*(t)=B_0^\top(t)\hat{q}(t)+\sum\limits_{i=1}^3B_i^\top(t)\hat{k}_i(t),
\end{equation}
where the $\mathcal{F}_t$-adapted process quadruple $(q(\cdot),k_1(\cdot),k_2(\cdot),k_3(\cdot))\in\mathbb{R}^n\times\mathbb{R}^n\times\mathbb{R}^n\times\mathbb{R}^n$ satisfies the adjoint {\it backward SDE} (BSDE)
\begin{equation}\label{adjoint equation of the follower}
\left\{
\begin{aligned}
-dq(t)=&\big[A_0(t)q(t)+\sum\limits_{i=1}^3A_i(t)k_i(t)-Q_1(t)x^{u_1^*,u_2}(t)\big]dt-\sum\limits_{i=1}^3k_i(t)dW_i(t),\ t\in[0,T],\\
  q(T)=&-G_1x^{u_1^*,u_2}(T).
\end{aligned}
\right.
\end{equation}
Taking clue from the terminal condition, we try to find
\begin{equation}\label{supposed form of q(t)}
q(t)=-P_1(t)x^{u_1^*,u_2}(t)-\phi(t),
\end{equation}
for some $\mathbb{R}^{n\times n}$-valued, deterministic, differentiable function $P_1(\cdot)$ with $P_1(T)=G_1$, and $\mathbb{R}^n$-valued, $\mathcal{F}_t$-adapted process $\phi(\cdot)$ which satisfies the BSDE
\begin{equation}\label{supposed equation of phi(t)}
\left\{
\begin{aligned}
d\phi(t)&=\alpha(t)dt+\beta_1(t)dW_1(t)+\beta_3(t)dW_3(t),\ t\in[0,T],\\
 \phi(T)&=0.
\end{aligned}
\right.
\end{equation}
In the above equation, $\alpha(\cdot),\beta_1(\cdot),\beta_3(\cdot)$ are all $\mathbb{R}^n$-valued, $\mathcal{F}_t$-adapted processes. Applying It\^{o}'s formula to (\ref{supposed form of q(t)}), we get
\begin{equation}\label{applying Ito's formula to q(t)}
\begin{aligned}
 -dq(t)&=\big[\dot{P}_1(t)x^{u_1^*,u_2}(t)+P_1(t)A_0(t)x^{u_1^*,u_2}(t)+P_1(t)B_0(t)u_1^*(t)+P^1(t)C_0(t)u_2(t)+\alpha(t)\big]dt\\
      &\quad+\sum\limits_{i=1,3}\Big\{P_1(t)\big[A^i(t)x^{u_1^*,u_2}(t)+B^i(t)u_1^*(t)+C^i(t)u_2(t)\big]+\beta^i(t)\Big\}dW_i(t)\\
      &\quad+P_1(t)\big[A^2(t)x^{u_1^*,u_2}(t)+B^2(t)u_1^*(t)+C^2(t)u_2(t)\big]dW_2(t).
\end{aligned}
\end{equation}
Comparing (\ref{applying Ito's formula to q(t)}) with (\ref{adjoint equation of the follower}), we have
\begin{equation}\label{comparing dWt}
\left\{
\begin{aligned}
k_1(t)&=-P_1(t)\big[A_1(t)x^{u_1^*,u_2}(t)+B_1(t)u_1^*(t)+C_1(t)u_2(t)\big]-\beta_1(t),\\
k_2(t)&=-P_1(t)\big[A_2(t)x^{u_1^*,u_2}(t)+B_2(t)u_1^*(t)+C_2(t)u_2(t)\big],\\
k_3(t)&=-P_1(t)\big[A_3(t)x^{u_1^*,u_2}(t)+B_3(t)u_1^*(t)+C_3(t)u_2(t)\big]-\beta_3(t),\\
\end{aligned}
\right.
\end{equation}
and
\begin{equation}\label{comparing dt}
\begin{aligned}
&A_0(t)q(t)+\sum\limits_{i=1}^3A_i(t)k_i(t)-Q_1(t)x^{u_1^*,u_2}(t)\\
&=\dot{P}_1(t)x^{u_1^*,u_2}(t)+P_1(t)A_0(t)x^{u_1^*,u_2}(t)+P_1(t)B_0(t)u_1^*(t)+P_1(t)C_0(t)u_2(t)+\alpha(t).
\end{aligned}
\end{equation}
Taking $\mathbb{E}[\cdot|\mathcal{G}^1_t]$ on both sides of (\ref{supposed form of q(t)}), (\ref{comparing dWt}) and (\ref{comparing dt}), we get
\begin{equation}\label{optimal filter for q}
\hat{q}(t)=-P_1(t)\hat{x}^{u^{1*},u^2}(t)-\hat{\phi}(t),
\end{equation}
\begin{equation}\label{optimal filter for k123}
\left\{
\begin{aligned}
\hat{k}_1(t)&=-P_1(t)\big[A_1(t)\hat{x}^{u^{1*},u^2}(t)+B_1(t)u_1^*(t)+C_1(t)\hat{u}_2(t)\big]-\hat{\beta}_1(t),\\
\hat{k}_2(t)&=-P_1(t)\big[A_2(t)\hat{x}^{u^{1*},u^2}(t)+B_2(t)u_1^*(t)+C_2(t)\hat{u}_2(t)\big],\\
\hat{k}_3(t)&=-P_1(t)\big[A_3(t)\hat{x}^{u^{1*},u^2}(t)+B_3(t)u_1^*(t)+C_3(t)\hat{u}_2(t)\big]-\hat{\beta}_3(t),\\
\end{aligned}
\right.
\end{equation}
and
\begin{equation}\label{filter}
\begin{aligned}
&A_0(t)\hat{q}(t)+\sum\limits_{i=1}^3A_i(t)\hat{k}_i(t)-Q_1(t)\hat{x}^{u_1^*,u_2}(t)\\
&=\dot{P}_1(t)\hat{x}^{u_1^*,u_2}(t)+P_1(t)A_0(t)\hat{x}^{u_1^*,u_2}(t)+P_1(t)B_0(t)u_1^*(t)+P_1(t)C_0(t)\hat{u}_2(t)+\hat{\alpha}(t).
\end{aligned}
\end{equation}
Applying Lemma 2.1 to (\ref{state equation}) corresponding to $u_1^*(\cdot)$ and (\ref{adjoint equation of the follower}) with $\mathbb{E}[\cdot|\mathcal{G}^1_t]$, we derive the follower's optimal filtering equation
\begin{equation}\label{optimal filter equation}
\left\{
\begin{aligned}
     d\hat{x}^{u_1^*,u_2}(t)&=\big[A_0(t)\hat{x}^{u_1^*,u_2}(t)+B_0(t)u_1^*(t)+C_0(t)\hat{u}_2(t)\big]dt\\
                            &\quad+\sum\limits_{i=1,3}\big[A_i(t)\hat{x}^{u_1^*,u_2}(t)+B_i(t)u_1^*(t)+C_i(t)\hat{u}_2(t)\big]dW_i(t),\\
                -d\hat{q}(t)&=\big[A_0(t)\hat{q}(t)+\sum\limits_{i=1}^3A_i(t)\hat{k}_i(t)-Q_1(t)\hat{x}^{u_1^*,u_2}(t)\big]dt\\
                            &\quad-\hat{k}_1(t)dW_1(t)-\hat{k}_3(t)dW_3(t),\ t\in[0,T],\\
      \hat{x}^{u_1^*,u_2}(0)&=x_0,\quad \hat{q}(T)=-G_1\hat{x}^{u_1^*,u_2}(T).
\end{aligned}
\right.
\end{equation}
Putting (\ref{optimal filter for q}), (\ref{optimal filter for k123}) into (\ref{optimal control of the follower}), we get
\begin{equation}\label{optimal control of the follower-feedback}
\begin{aligned}
u_1^*(t)=
&-\Big[N_1(t)+\sum\limits_{i=1}^3B_i^\top(t)P_1(t)B_i(t)\Big]^{-1}\bigg\{\Big[B_0^\top(t)P_1(t)+\sum\limits_{i=1}^3B_i^\top(t)P_1(t)A_i(t)\Big]\hat{x}^{u_1^*,u_2}(t)\\
&\quad+B_0^\top(t)P_1(t)\hat{\phi}(t)+B_1^\top(t)\hat{\beta}_1(t)+B_3^\top(t)\hat{\beta}_3(t)+\Big(\sum\limits_{i=1}^3B_i^\top(t)P_1(t)C_i(t)\Big)\hat{u}_2(t)\bigg\},
\end{aligned}
\end{equation}
where we have assumed that

\vspace{1mm}

\noindent{\bf (A2.1)}\quad{\it $\overline{N}_1(t)\triangleq N_1(t)+\sum\limits_{i=1}^3B_i^\top(t)P_1(t)B_i(t)>0$, $\forall t\in[0,T]$.}

\vspace{1mm}

Substituting (\ref{optimal filter for q}), (\ref{optimal filter for k123}) and (\ref{optimal control of the follower-feedback}) into (\ref{filter}), we obtain the Riccati's type equation
\begin{equation}\label{Riccati equation}
\left\{
\begin{aligned}
 &\dot{P}_1(t)+P_1(t)A_0(t)+A_0^\top(t)P_1(t)+\sum\limits_{i=1}^3A_i^\top(t)P_1(t)A_i(t)+Q_1(t)\\
 &-\Big[P_1(t)B_0(t)+\sum\limits_{i=1}^3A_i^\top(t)P_1(t)B_i(t)\Big]\overline{N}_1^{-1}(t)\Big[B_0^\top(t)P_1(t)+\sum\limits_{i=1}^3B_i^\top(t)P_1(t)A_i(t)\Big]=0,\\
 &P_1(T)=G_1,
\end{aligned}
\right.
\end{equation}
which admits a unique solution by Theorem 7.10, Chapter 6 of Yong and Zhou \cite{YZ99}. Then
\begin{equation}\label{alpha(t)-explicit form}
\hat{\alpha}(t)=-L_0(t)\hat{\phi}(t)-L_1(t)\hat{\beta}_1(t)-L_3(t)\hat{\beta}_3(t)-L_4(t)\hat{u}_2(t),
\end{equation}
where ($t$ is omitted for simplification)
{\small\begin{equation}
\left\{
\begin{aligned}
 L_0\triangleq&\overline{N}_1^{-1}\Big(P_1B_0+\sum\limits_{i=1}^3A_i^\top P_1B_i\Big)B_0^\top P-A_0,\\
 L_j\triangleq&\overline{N}_1^{-1}\Big(P_1B_0+\sum\limits_{i=1}^3A_i^\top P_1B_i\Big)B_j^\top-A_j,\ j=1,3,\\
 L_4\triangleq&\overline{N}_1^{-1}\Big(P_1B_0+\sum\limits_{i=1}^3A_i^\top P_1B_i\Big)\Big(\sum\limits_{i=1}^3B_i^\top P_1C_i\Big)-P_1C_0-\sum\limits_{i=1}^3A_i^\top P_1C_i.
\end{aligned}
\right.
\end{equation}}
Applying Lemma 2.1 again to BSDE (\ref{supposed equation of phi(t)}), we have
\begin{equation}\label{BSDFE}
\left\{
\begin{aligned}
-d\hat{\phi}(t)&=\big[L_0(t)\hat{\phi}(t)+L_1(t)\hat{\beta}_1(t)+L_3(t)\hat{\beta}_3(t)+L_4(t)\hat{u}_2(t)\big]dt\\
               &\quad-\hat{\beta}_1(t)dW_1(t)-\hat{\beta}_3(t)dW_3(t),\ t\in[0,T],\\
  \hat{\phi}(T)&=0.
\end{aligned}
\right.
\end{equation}
For given $u_2(\cdot)$, (\ref{BSDFE}) admits a unique $\mathcal{G}^1_t$-adapted solution triple $(\hat{\phi}(\cdot),\hat{\beta}_1(\cdot),\hat{\beta}_3(\cdot))$ by the standard BSDE theory (See, for example, El Karoui et al. \cite{EPQ97}). Putting (\ref{optimal control of the follower-feedback}) into the forward equation in (\ref{optimal filter equation}), we get
\begin{equation}\label{optimal filter equation of the follower}
\left\{
\begin{aligned}
  d\hat{x}^{u_1^*,u_2}(t)&=\bigg\{\Big[A_0(t)-B_0(t)\overline{N}_1^{-1}(t)\Big(B_0^\top(t)P_1(t)+\sum\limits_{i=1}^3B_i^\top(t)P_1(t)A_i(t)\Big)\Big]\hat{x}^{u_1^*,u_2}(t)\\
                         &\qquad-B_0(t)\overline{N}_1^{-1}(t)B_0^\top(t)P_1(t)\hat{\phi}(t)-B_0(t)\overline{N}_1^{-1}(t)B_1^\top(t)\hat{\beta}_1(t)-B_0(t)\overline{N}_1^{-1}(t)\\
                         &\qquad\times B_3^\top(t)\hat{\beta}_3(t)+\Big[C_0(t)-B_0(t)\overline{N}_1^{-1}(t)\Big(\sum\limits_{i=1}^3B_i^\top(t)P_1(t)C_i(t)\Big)\Big]\hat{u}_2(t)\bigg\}dt\\
                         &\quad+\sum\limits_{i=1,3}\bigg\{\Big[A_i(t)-B_i(t)\overline{N}_1^{-1}(t)\Big(B_0^\top(t)P_1(t)
                          +\sum\limits_{i=1}^3B_i^\top(t)P_1(t)A_i(t)\Big)\Big]\hat{x}^{u_1^*,u_2}(t)\\
                         &\qquad-B_i(t)\overline{N}_1^{-1}(t)B_0^\top(t)P_1(t)\hat{\phi}(t)-B_i(t)\overline{N}_1^{-1}(t)B_1^\top(t)\hat{\beta}_1(t)\\
                         &\qquad-B_i(t)\overline{N}_1^{-1}(t)B_3^\top(t)\hat{\beta}_3(t)\bigg\}dW_i(t),\ t\in[0,T],\\
   \hat{x}^{u_1^*,u_2}(0)&=x_0,
\end{aligned}
\right.
\end{equation}
which admits a unique $\mathcal{G}^1_t$-adapted solution $\hat{x}^{u^{1*},u^2}(\cdot)$, from (\ref{BSDFE}). In fact, for given $u^2(\cdot)$, we can verify the solvability of (\ref{optimal filter equation}). The optimal control $u_1^*(\cdot)$ is expressed by (\ref{optimal control of the follower-feedback}).

Moreover, it is easy to check that the concavity/convexity conditions in the verification theorem (Please refer to Proposition 2.2 of \cite{SWX16}) hold, then $u_1^*(\cdot)$ given by (\ref{optimal control of the follower-feedback}) is really optimal. We summarize the above argument in the following theorem.

\vspace{1mm}

\noindent{\bf Theorem 2.1}\quad{\it Let ${\bf (A2.1)}$ hold and $P_1(\cdot)$ satisfy (\ref{Riccati equation}). For chosen $u_2(\cdot)$ of the leader, $u_1^*(\cdot)$ defined by (\ref{optimal control of the follower-feedback}) is an optimal control of the follower, where $(\hat{x}^{u_1^*,u_2}(\cdot),\hat{\phi}(\cdot),\hat{\beta}_1(\cdot),\hat{\beta}_3(\cdot))$ is determined by (\ref{BSDFE}) and (\ref{optimal filter equation of the follower}).}

\subsection{Problem of the leader}

In this subsection, since the follower's optimal control $u_1^*(\cdot)$ by (\ref{optimal control of the follower-feedback}) is a linear functional of $\hat{x}^{u_1^*,u_2}(\cdot),\hat{\phi}(\cdot),\hat{\beta}_1(\cdot),\hat{\beta}_3(\cdot)$ and $\hat{u}_2(\cdot)$, the leader's state equation now writes
\begin{equation}\label{state equation-leader}
\left\{
\begin{aligned}
    dx^{u_2}(t)&=\Big[A_0(t)x^{u_2}(t)+L_{01}(t)\hat{x}^{u_2}(t)+L_{02}(t)\hat{\phi}(t)+L_{03}(t)\hat{\beta}_1(t)+L_{04}(t)\hat{\beta}_3(t)\\
               &\qquad+C_0(t)u_2(t)+L_{05}(t)\hat{u}_2(t)\Big]dt+\sum\limits_{i=1}^3\Big[A_i(t)x^{u_2}(t)+L_{i1}(t)\hat{x}^{u_2}(t)\\
               &\quad+L_{i2}(t)\hat{\phi}(t)+L^{i3}(t)\hat{\beta}^1(t)+L^{i4}(t)\hat{\beta}^3(t)+C^i(t)u^2(t)+L_{i5}(t)\hat{u}_2(t)\Big]dW_i(t),\\
-d\hat{\phi}(t)&=\big[L_0(t)\hat{\phi}(t)+L_1(t)\hat{\beta}_1(t)+L_3(t)\hat{\beta}_3(t)+L_4(t)\hat{u}_2(t)\big]dt\\
               &\qquad-\hat{\beta}_1(t)dW_1(t)-\hat{\beta}_3(t)dW_3(t),\ t\in[0,T],\\
     x^{u_2}(0)&=x_0,\quad \hat{\phi}(T)=0,
\end{aligned}
\right.
\end{equation}
where we denote $x^{u_2}\equiv x^{u_1^*,u_2}$, $\hat{x}^{u_2}\equiv\hat{x}^{u_1^*,u_2}$ and for $j=0,1,2,3$,
{\small\begin{equation}
\left\{
\begin{aligned}
 L_{j1}\triangleq&-B_j\overline{N}_1^{-1}\Big(B_0^\top P_1+\sum\limits_{i=1}^3B_i^\top P_1A_i\Big),\ L_{j2}\triangleq-B_j\overline{N}_1^{-1}B_0^\top P_1,\ L_{j3}\triangleq-B_j\overline{N}_1^{-1}B_1^\top,\\
 L_{j4}\triangleq&-B_j\overline{N}_1^{-1}B_3^\top,\ L_{j5}\triangleq-B_j\overline{N}_1^{-1}\Big(\sum\limits_{i=1}^3B_i^\top P_1C_i\Big).
\end{aligned}
\right.
\end{equation}}
The problem of the leader is to select a $\mathcal{G}^2_t$-adapted optimal control $u_2^*(\cdot)$ such that the cost functional
\begin{equation}\label{cost functional-leader-simple}
\begin{aligned}
&J_2(u_2(\cdot))\equiv J_2(u_1^*(\cdot),u_2(\cdot))\\
&=\frac{1}{2}\mathbb{E}\left[\int_0^T\Big(\big\langle Q_2(t)x^{u_2}(t),x^{u_2}(t)\big\rangle+\big\langle N_2(t)u_2(t),u_2(t)\big\rangle\Big)dt+\big\langle G_2x^{u_2}(T),x^{u_2}(T)\big\rangle\right]
\end{aligned}
\end{equation}
is minimized.

Suppose that there exists a $\mathcal{G}^2_t$-adapted optimal control $u_2^*(\cdot)$ of the leader, and his optimal state is $(x^*(\cdot),\hat{\phi}^*(\cdot),\hat{\beta}_1^*(\cdot),\hat{\beta}_3^*(\cdot))\equiv(x^{u_2^*}(\cdot),\hat{\phi}^*(\cdot),\hat{\beta}_1^*(\cdot),\hat{\beta}_3^*(\cdot))$. Next, we will derive the necessary condition for $u_2^*(\cdot)$, by a direct calculation of the derivative of the cost functional. We define the leader's Hamiltonian function
\begin{equation}\label{Hamiltonian function of the leader}
\begin{aligned}
& H_2\big(t,x^{u_2},u_2,\phi,\beta_1,\beta_3;p,y,z_1,z_2,z_3\big)\\
&\hspace{-4mm}\triangleq \big\langle y,A_0(t)x^{u_2}+L_{01}(t)\hat{x}^{u_2}+L_{02}(t)\hat{\phi}+L_{03}(t)\hat{\beta}_1+L_{04}(t)\hat{\beta}_3+C_0(t)u_2+L_{05}(t)\hat{u}_2\big\rangle\\
&+\big\langle p,L_0(t)\hat{\phi}+L_1(t)\hat{\beta}_1+L_3(t)\hat{\beta}_3+L_4(t)\hat{u}_2\big\rangle
 +\frac{1}{2}\big\langle Q_2(t)x^{u_2},x^{u_2}\big\rangle+\frac{1}{2}\big\langle N_2(t)u_2,u_2\big\rangle\\
&+\sum\limits_{i=1}^3\big\langle z_i,A_i(t)x^{u_2}+L_{i1}(t)\hat{x}^{u_2}+L_{i2}(t)\hat{\phi}+L_{i3}(t)\hat{\beta}_1+L_{i4}(t)\hat{\beta}_3+C_0(t)u_2+L_{i5}(t)\hat{u}_2\big\rangle,
\end{aligned}
\end{equation}
where the $\mathcal{F}_t$-adapted process quintuple $(p(\cdot),y(\cdot),z_1(\cdot),z_2(\cdot),z_3(\cdot))\in\mathbb{R}^n\times\mathbb{R}^n\times\mathbb{R}^n\times\mathbb{R}^n\times\mathbb{R}^n$ satisfies the adjoint equation
\begin{equation}\label{adjoint equation of the leader}
\left\{
\begin{aligned}
  dp(t)&=\Big[L_{02}(t)y(t)+L_0(t)p(t)+\sum\limits_{i=1}^3L_{i2}(t)z_i(t)\Big]dt\\
       &\quad+\Big[L_{03}(t)y(t)+L_1(t)p(t)+\sum\limits_{i=1}^3L_{i3}(t)z_i(t)\Big]dW_1(t)\\
       &\quad+\Big[L_{04}(t)y(t)+L_3(t)p(t)+\sum\limits_{i=1}^3L_{i4}(t)z_i(t)\Big]dW_3(t),\\
 -dy(t)&=\Big[A_0(t)y(t)+L_{01}(t)\hat{y}(t)+\sum\limits_{i=1}^3A_i(t)z_i(t)+\sum\limits_{i=1}^3L_{i1}(t)\hat{z}_i(t)+Q_2(t)x^*(t)\Big]dt\\
       &\quad-z_1(t)dW_1(t)-z_2(t)dW_2(t)-z_3(t)dW_3(t),\ t\in[0,T],\\
   p(0)&=0,\quad y(T)=G_2x^*(T).
\end{aligned}
\right.
\end{equation}

Without loss of generality, let $x_0\equiv0$, and define the perturbed optimal control $u_2^*(\cdot)+\epsilon u_2(\cdot)$ for sufficiently small $\epsilon>0$, with any $u_2(\cdot)$. Then it is easy to see from the linearity of (\ref{state equation-leader}), that the solution to it is $x^*(\cdot)+\epsilon x^{u_2}(\cdot)$. First we have
\begin{equation*}
\begin{aligned}
&\hspace{-3mm}\widetilde{J}(\epsilon)\triangleq J_2\big(u_2^*(\cdot)+\epsilon u_2(\cdot)\big)=\frac{1}{2}\mathbb{E}\int_0^T\Big[\big\langle Q_2(t)(x^*(t)+\epsilon x^{u_2}(t)),x^*(t)+\epsilon x^{u_2}(t)\big\rangle\\
&+\big\langle N_2(t)(u_2^*(t)+\epsilon u_2(t)),u_2^*(t)+\epsilon u_2(t)\big\rangle\Big]dt+\frac{1}{2}\mathbb{E}\big\langle G_2(x^*(T)+\epsilon x^{u_2}(T)),x^*(T)+\epsilon x^{u_2}(T)\big\rangle.
\end{aligned}
\end{equation*}
Hence
\begin{equation}\label{equality}
\begin{aligned}
0&=\frac{\partial\widetilde{J}(\epsilon)}{\partial\epsilon}\bigg|_{\epsilon=0}
  =\mathbb{E}\int_0^T\Big[\big\langle Q_2(t)x^*(t),x^{u_2}(t)\big\rangle+\big\langle N_2(t)u_2^*(t),u_2(t)\big\rangle\Big]dt+\mathbb{E}\big\langle G_2x^*(T),x^{u_2}(T)\big\rangle\\
 &=\mathbb{E}\int_0^T\Big[\big\langle Q_2(t)x^*(t),x^{u_2}(t)\big\rangle+\big\langle N_2(t)u_2^*(t),u_2(t)\big\rangle\Big]dt+\mathbb{E}\big\langle y(T),x^{u_2}(T)\big\rangle.
\end{aligned}
\end{equation}
Applying It\^{o}'s formula to $\langle y(t),x^{u_2}(t)\rangle-\langle p(t),\hat{\phi}(t)\rangle$, noting (\ref{state equation-leader}) and (\ref{adjoint equation of the leader}), we derive
\begin{equation*}
\begin{aligned}
 &d\langle y(t),x^{u_2}(t)\rangle-d\langle p(t),\hat{\phi}(t)\rangle\\
=&\bigg\langle y(t),\big[L_{01}(t)\hat{x}^{u_2}(t)+C_0(t)u_2(t)+L_{05}(t)\hat{u}_2(t)\big]dt+\sum\limits_{i=1}^3\Big[A_i(t)x^{u_2}(t)+L_{i1}(t)\hat{x}^{u_2}(t)\\
 &\quad+L_{i2}(t)\hat{\phi}(t)+L_{i3}(t)\hat{\beta}_1(t)+L_{i4}(t)\hat{\beta}_3(t)+C_i(t)u_2(t)+L_{i5}(t)\hat{u}_2(t)\Big]dW_i(t)\bigg\rangle\\
 &-\bigg\langle x^{u_2}(t),\Big[L_{01}(t)\hat{y}(t)+\sum\limits_{i=1}^3L_{i1}(t)\hat{z}_i(t)+Q_2(t)x^*(t)\Big]dt-\sum\limits_{i=1}^3z_i(t)dW_i(t)\bigg\rangle\\
 &+\sum\limits_{i=1}^3\Big\langle z_i(t),\big[L_{i1}(t)\hat{x}^{u_2}(t)+C_i(t)u_2(t)+L_{i5}(t)\hat{u}_2(t)\big]dt\Big\rangle+\big\langle p(t),L_4(t)\hat{u}_2(t)dt\big\rangle\\
\end{aligned}
\end{equation*}
\begin{equation*}
\begin{aligned}
 &+\hat{\beta}_1(t)dW_1(t)+\hat{\beta}_3(t)dW_3(t)-\bigg\langle\hat{\phi}(t),\Big[L_{03}(t)y(t)+L_1(t)p(t)+\sum\limits_{i=1}^3L_{i3}(t)z_i(t)\Big]dW_1(t)\\
 &\quad+\Big[L_{04}(t)y(t)+L_3(t)p(t)+\sum\limits_{i=1}^3L_{i4}(t)z_i(t)\Big]dW_3(t)\bigg\rangle.
\end{aligned}
\end{equation*}
Therefore,
\begin{equation*}
\begin{aligned}
  &\hspace{-3mm}\mathbb{E}\big\langle y(T),x^{u_2}(T)\big\rangle=-\mathbb{E}\int_0^T\Big[\big\langle Q_2(t)x^*(t),x^{u_2}(t)\big\rangle+\big\langle N_2(t)u_2^*(t),u_2(t)\big\rangle\Big]dt\\
  &=\mathbb{E}\int_0^T\Big\langle y(t),\big[L_{01}(t)\hat{x}^{u_2}(t)+C_0(t)u_2(t)+L_{05}(t)\hat{u}_2(t)\big]dt\Big\rangle+\mathbb{E}\int_0^T\big\langle p(t),L_4(t)\hat{u}_2(t)dt\big\rangle\\
  &\quad-\mathbb{E}\int_0^T\bigg\langle x^{u_2}(t),\Big[L_{01}(t)\hat{y}(t)+\sum\limits_{i=1}^3L_{i1}(t)\hat{z}_i(t)+Q_2(t)x^*(t)\Big]dt\bigg\rangle\\
  &\quad+\sum\limits_{i=1}^3\mathbb{E}\int_0^T\big\langle z_i(t),\big(L_{i1}(t)\hat{x}^{u_2}(t)+C_i(t)u_2(t)+L_{i5}(t)\hat{u}_2(t)\big)dt\big\rangle.
\end{aligned}
\end{equation*}
Noting that
$$\mathbb{E}\int_0^T\langle\mathbb{E}[\xi|\mathcal{G}^1_t],\eta\rangle dt=\mathbb{E}\int_0^T\langle\xi,\mathbb{E}[\eta|\mathcal{G}^1_t]\rangle dt,\quad
\mathbb{E}\int_0^T\langle\mathbb{E}[\xi|\mathcal{G}^2_t],\eta\rangle dt=\mathbb{E}\int_0^T\langle\xi,\mathbb{E}[\eta|\mathcal{G}^2_t]\rangle dt$$
for any $\mathcal{F}_t$-adapted random variables $\xi,\eta$, we have
\begin{equation*}
\begin{aligned}
0=&\ \mathbb{E}\int_0^T\big\langle N_2(t)u_2^*(t),u_2(t)\big\rangle dt+\mathbb{E}\int_0^T\Big\langle y(t),\big[C_0(t)u_2(t)+L_{05}(t)\hat{u}_2(t)\big]dt\Big\rangle\\
  &+\sum\limits_{i=1}^3\mathbb{E}\int_0^T\Big\langle z_i(t),\big[C_i(t)u_2(t)+L_{i5}(t)\hat{u}_2(t)\big]dt\Big\rangle
   +\mathbb{E}\int_0^T\big\langle p(t),L_4(t)\hat{u}_2(t)dt\big\rangle\\
 =&\ \mathbb{E}\int_0^T\big\langle N_2(t)u_2^*(t),u_2(t)\big\rangle dt+\mathbb{E}\int_0^T\Big\langle L_4^\top(t)\hat{p}(t)+C_0^\top(t)y(t)+L_{05}^\top(t)\hat{y}(t)\\
  &\qquad+\sum\limits_{i=1}^3C_i^\top (t) z_i(t)+\sum\limits_{i=1}^3L_{i5}^\top(t)\hat{z}_i(t),u_2(t)\Big\rangle dt.
\end{aligned}
\end{equation*}
This implies that
\begin{equation}\label{optimal control of the leader}
\begin{aligned}
\hspace{-4mm}u_2^*(t)=-N_2^{-1}(t)\bigg[L_4^\top(t)\check{\hat{p}}(t)+C_0^\top(t)\check{y}(t)+L_{05}^\top(t)\check{\hat{y}}(t)+\sum\limits_{i=1}^3C_i^\top(t)\check{z}_i(t)
+\sum\limits_{i=1}^3L_{i5}^\top(t)\check{\hat{z}}_i(t)\bigg],
\end{aligned}
\end{equation}
where we have used that
$$
\mathbb{E}\int_0^T\langle\xi,\eta\rangle dt=\mathbb{E}\int_0^T\mathbb{E}\big[\langle\xi,\eta\rangle\big|\mathcal{G}^2_t\big]dt=\mathbb{E}\int_0^T\langle\mathbb{E}[\xi|\mathcal{G}^2_t],\eta\rangle dt
$$
for any $\mathcal{F}_t$-adapted random variable $\xi$ and $\mathcal{G}^2_t$-adapted random variable $\eta$.

In the following, we will derive the filtering equation for $\check{p}(\cdot),\check{y}(\cdot),\check{z}_i(\cdot)$ and $\check{\hat{p}}(\cdot),\check{\hat{y}}(\cdot),\check{\hat{z}}_i(\cdot)$. Applying again Lemma 2.1 to (\ref{adjoint equation of the leader}) and (\ref{state equation-leader}) corresponding to $u_2^*(\cdot)$ with $\mathbb{E}[\cdot|\mathcal{G}^2_t]$, we obtain the leader's optimal filtering equation
\begin{equation}\label{optimal filter equation-leader}
\left\{
\begin{aligned}
          d\check{x}^*(t)&=\big[A_0\check{x}^*+L_{01}\check{\hat{x}}^*+L_{02}\check{\hat{\phi}}^*+L_{03}\check{\hat{\beta}}_1^*+L_{04}\check{\hat{\beta}}_3^*+C_0u_2^*+L_{05}\hat{u}_2^*\big]dt\\
                         &\quad+\sum\limits_{i=2,3}\big[A_i\check{x}^*+L_{i1}\check{\hat{x}}^*+L_{i2}\check{\hat{\phi}}^*
                          +L_{i3}\check{\hat{\beta}}_1^*+L_{i4}\check{\hat{\beta}}_3^*+C_iu_2^*+L_{i5}\hat{u}_2^*\big]dW_i(t),\\
            d\check{p}(t)&=\Big[L_{02}\check{y}+L_0\check{p}+\sum\limits_{i=1}^3L_{i2}\check{z}_i\Big]dt+\Big[L_{04}\check{y}+L_3\check{p}
                          +\sum\limits_{i=1}^3L_{i4}\check{z}_i\Big]dW_3(t),\\
           -d\check{y}(t)&=\Big[\big(A_0+L_{01}\big)\check{y}+\sum\limits_{i=1}^3\big(A_i+L_{i1}\big)\check{z}_i+Q_2\check{x}^*\Big]dt-\check{z}_2dW_2(t)-\check{z}_3dW_3(t),\\
-d\check{\hat{\phi}}^*(t)&=\big[L_0\check{\hat{\phi}}^*+L_1\check{\hat{\beta}}_1^*+L_3\check{\hat{\beta}}_3^*+L_4\hat{u}_2^*\big]dt-\check{\hat{\beta}}_3dW_3(t),\ t\in[0,T],\\
           \check{x}^*(0)&=x_0,\quad \check{\hat{\phi}}^*(T)=0,\quad \check{p}(0)=0,\quad \check{y}(T)=G_2\check{x}^*(T).
\end{aligned}
\right.
\end{equation}
Putting (\ref{optimal control of the leader}) into it, we get
\begin{equation}\label{optimal filter equation-leader-without u2}
\left\{
\begin{aligned}
          d\check{x}^*(t)&=\bigg\{A_0\check{x}^*+L_{01}\check{\hat{x}}^*+L_{02}\check{\hat{\phi}}^*+L_{03}\check{\hat{\beta}}_1^*+L_{04}\check{\hat{\beta}}_3^*-C_0N_2^{-1}\Big[L_4^\top\check{\hat{p}}
                          +C_0^\top\check{y}+L_{05}^\top\check{\hat{y}}\\
                         &\qquad+\sum\limits_{j=1}^3C_j^\top(t)\check{z}_j(t)+\sum\limits_{j=1}^3L_{j5}^\top(t)\check{\hat{z}}_j(t)\Big]-L_{05}N_2^{-1}\Big[L_4^\top\check{\hat{p}}+\big(C_0+L_{05}\big)^\top\check{\hat{y}}\\
                         &\qquad+\sum\limits_{j=1}^3\big(C_j+L_{j5}\big)^\top\check{\hat{z}}_j\Big]\bigg\}dt
                          +\sum\limits_{i=2,3}\bigg\{A_i\check{x}^*+L_{i1}\check{\hat{x}}^*+L_{i2}\check{\hat{\phi}}^*+L_{i3}\check{\hat{\beta}}_1^*\\
                         &\qquad+L_{i4}\check{\hat{\beta}}_3^*-C_iN_2^{-1}\Big[L_4^\top\check{\hat{p}}+C_0^\top\check{y}+L_{05}^\top\check{\hat{y}}+\sum\limits_{j=1}^3C_j^\top\check{z}_j+\sum\limits_{j=1}^3L_{j5}^\top\check{\hat{z}}_j\Big]\\
                         &\qquad-L_{i5}N_2^{-1}\Big[L_4^\top\check{\hat{p}}+\big(C_0+L_{05}\big)^\top\check{\hat{y}}+\sum\limits_{j=1}^3\big(C_j+L_{j5}\big)^\top\check{\hat{z}}_j\Big]\bigg\}dW_i(t),\\
            d\check{p}(t)&=\Big[L_{02}\check{y}+L_0\check{p}+\sum\limits_{i=1}^3L_{i2}\check{z}_i\Big]dt
                          +\Big[L_{04}\check{y}+L_3\check{p}+\sum\limits_{i=1}^3L_{i4}\check{z}_i\Big]dW_3(t),\\
           -d\check{y}(t)&=\Big[\big(A_0+L_{01}\big)\check{y}+\sum\limits_{j=1}^3\big(A_j+L_{j1}\big)\check{z}_j+Q_2\check{x}^*\Big]dt-\check{z}_2dW_2(t)-\check{z}_3dW_3(t),\\
-d\check{\hat{\phi}}^*(t)&=\bigg\{L_0\check{\hat{\phi}}^*+L_1\check{\hat{\beta}}_1^*+L_3\check{\hat{\beta}}_3^*-L_4N_2^{-1}\Big[L_4^\top\check{\hat{p}}+\big(C_0+L_{05}\big)^\top\check{\hat{y}}\\
                         &\qquad+\sum\limits_{j=1}^3\big(C_j+L_{j5}\big)^\top\check{\hat{z}}_j\Big]\bigg\}dt-\check{\hat{\beta}}_3dW_3(t),\ t\in[0,T],\\
           \check{x}^*(0)&=x_0,\quad \check{\hat{\phi}}^*(T)=0,\quad \check{p}(0)=0,\quad \check{y}(T)=G_2\check{x}^*(T),
\end{aligned}
\right.
\end{equation}
where the nine tuple $(\check{\hat{x}}^*(\cdot),\check{\hat{\phi}}^*(\cdot),\check{\hat{\beta}}_1^*(\cdot),\check{\hat{\beta}}_3^*(\cdot),\check{\hat{p}}^*(\cdot),\check{\hat{y}}^*(\cdot),
\check{\hat{z}}_1^*(\cdot),\check{\hat{z}}_2^*(\cdot),\check{\hat{z}}_3^*(\cdot))$ satisfies
{\small\begin{equation}\label{optimal filter equation-leader and follower-without u2}
\left\{
\begin{aligned}
    d\check{\hat{x}}^*(t)&=\bigg\{\big(A_0+L_{01}\big)\check{\hat{x}}^*+L_{02}\check{\hat{\phi}}^*+L_{03}\check{\hat{\beta}}_1^*+L_{04}\check{\hat{\beta}}_3^*
                          -\big(C_0+L_{05}\big)N_2^{-1}\Big[L_4^\top\check{\hat{p}}+\big(C_0+L_{05}\big)^\top\check{\hat{y}}\\
                         &\qquad+\sum\limits_{i=1}^3\big(C_i+L_{i5}\big)^\top\check{\hat{z}}_i\Big]\bigg\}dt+\bigg\{\big(A_3+L_{31}\big)\check{\hat{x}}^*+L_{32}\check{\hat{\phi}}^*
                          +L_{33}\check{\hat{\beta}}^{1*}+L_{34}\check{\hat{\beta}}_3^*\\
                         &\qquad-\big(C_3+L_{35}\big)N_2^{-1}\Big[L_4^\top\check{\hat{p}}+\big(C_0+L_{05}\big)^\top\check{\hat{y}}
                          +\sum\limits_{i=1}^3\big(C_3+L_{35}\big)^\top\check{\hat{z}}_3\Big]\bigg\}dW_3(t),\\
      d\check{\hat{p}}(t)&=\Big[L_{02}\check{\hat{y}}+L_0\check{\hat{p}}+\sum\limits_{i=1}^3L_{i2}\check{\hat{z}}_i\Big]dt+\Big[L_{04}\check{\hat{y}}
                          +L_3\check{\hat{p}}+\sum\limits_{i=1}^3L_{i4}\check{\hat{z}}_i\Big]dW_3(t),\\
-d\check{\hat{\phi}}^*(t)&=\bigg\{L_0\check{\hat{\phi}}^*+L_1\check{\hat{\beta}}_1^*+L_3\check{\hat{\beta}}_3^*-L_4N_2^{-1}\Big[L_4^\top\check{\hat{p}}+\big(C_0+L_{05}\big)^\top\check{\hat{y}}
                          +\sum\limits_{i=1}^3\big(C_i+L_{i5}\big)^\top\check{\hat{z}}_i\Big]\bigg\}dt\\
                         &\quad-\check{\hat{\beta}}_3dW_3(t),\\
     -d\check{\hat{y}}(t)&=\bigg[\big(A_0+L_{01}\big)\check{\hat{y}}+\sum\limits_{i=1}^3\big(A_i+L_{i1}\big)\check{\hat{z}}_i+Q_2\check{\hat{x}}^*\bigg]dt
                          -\check{\hat{z}}_3dW_3(t),\ t\in[0,T],\\
     \check{\hat{x}}^*(0)&=x_0,\quad \check{\hat{p}}(0)=0,\quad \check{\hat{\phi}}^*(T)=0,\quad\check{\hat{y}}(T)=G_2\check{\hat{x}}^*(T).
\end{aligned}
\right.
\end{equation}}

Up to now, we have obtained the optimal control $u_2^*(\cdot)$ of the leader by (\ref{optimal control of the leader}). However, this representation relies on the solvability of filtering equations (\ref{optimal filter equation-leader-without u2}) and (\ref{optimal filter equation-leader and follower-without u2}). In the following, we will derive the state estimate feedback representation of (\ref{optimal control of the leader}), via some Riccati type equations. And the solvability of the above filtering equations will be solved as a corollary.

For this target, first we rewrite the optimal state of the leader as
{\small\begin{equation}\label{optimal state equation-leader}
\left\{
\begin{aligned}
          dx^*(t)&=\bigg\{A_0x^*+L_{01}\hat{x}^*+L_{02}\hat{\phi}^*+L_{03}\hat{\beta}_1^*+L_{04}\hat{\beta}_3^*-C^0N_2^{-1}\Big[L_4^\top\check{\hat{p}}+C_0^\top\check{y}+L_{05}^\top\check{\hat{y}}\\
                 &\qquad+\sum\limits_{j=1}^3C_j^\top\check{z}_j+\sum\limits_{j=1}^3L_{j5}^\top\check{\hat{z}}_j\Big]-L_{05}N_2^{-1}\Big[L_4^\top\check{\hat{p}}+\big(C_0+L_{05}\big)^\top\check{\hat{y}}\\
                 &\qquad+\sum\limits_{i=1}^3\bigg\{A_ix^*+L_{i1}\hat{x}^*+L_{i2}\hat{\phi}^*+L_{i3}\hat{\beta}_1^*+L_{i4}\hat{\beta}_3^*\\
                 &\qquad-C_iN_2^{-1}\Big[L_4^\top\check{\hat{p}}+C_0^\top\check{y}+L_{05}^\top\check{\hat{y}}+\sum\limits_{j=1}^3C_j^\top\check{z}_j+\sum\limits_{j=1}^3L_{j5}^\top\check{\hat{z}}_j\Big]\\
                 &\qquad-L_{i5}N_2^{-1}\Big[L_4^\top\check{\hat{p}}+\big(C_0+L_{05}\big)^\top\check{\hat{y}}+\sum\limits_{j=1}^3\big(C_j+L_{j5}\big)^\top\check{\hat{z}}_j\Big]\bigg\}dW_i(t),\\
-d\hat{\phi}^*(t)&=\bigg\{L_0\hat{\phi}^*+L_1\hat{\beta}_1^*+L_3\hat{\beta}_3^*-L_4N_2^{-1}\Big[L_4^\top\check{\hat{p}}+\big(C_0+L_{05}\big)^\top\check{\hat{y}}
                  +\sum\limits_{j=1}^3\big(C_j+L_{j5}\big)^\top\check{\hat{z}}_j\Big]\bigg\}dt\\
                 &\qquad-\hat{\beta}_1^*dW_1(t)-\hat{\beta}_3^*dW_3(t),\ t\in[0,T],\\
           x^*(0)&=x_0,\quad \hat{\phi}^*(T)=0.
\end{aligned}
\right.
\end{equation}}

Now, inspired by \cite{Yong02}, let
\begin{equation}\label{new state}
\begin{aligned}
&X\triangleq\left(\begin{array}{c}x^*\\p\end{array}\right),\ Y\triangleq\left(\begin{array}{c}y\\\hat{\phi}^*\end{array}\right),\ Z_1\triangleq\left(\begin{array}{c}z_1\\\hat{\beta}_1^*\end{array}\right),\
Z_2\triangleq\left(\begin{array}{c}z_2\\0\end{array}\right),\ Z_3\triangleq\left(\begin{array}{c}z_3\\\hat{\beta}_3^*\end{array}\right),
\end{aligned}
\end{equation}
then (\ref{optimal state equation-leader}) and (\ref{adjoint equation of the leader}) can be rewritten as
\begin{equation}\label{optimality system-leader-2 dim}
\left\{
\begin{aligned}
 dX(t)=&\big(\mathcal{A}_0X+\widehat{\mathcal{A}}_0\hat{X}+\overline{\mathcal{A}}_0\check{\hat{X}}+\mathcal{B}_0Y+\mathcal{C}_0\check{Y}+\widetilde{\mathcal{C}}_0\check{\hat{Y}}
        +\mathcal{B}_1^\top Z_1+\widetilde{\mathcal{B}}_1^\top\check{Z}_1\\
       &\ +\overline{\mathcal{C}}_0\check{\hat{Z}}_1+\mathcal{B}_2^\top Z_2+\widetilde{\mathcal{B}}_2^\top\check{Z}_2+\overline{\mathcal{D}}_0\check{\hat{Z}}_2
        +\mathcal{B}_3^\top Z_3+\widetilde{\mathcal{B}}_3^\top\check{Z}_3+\overline{\mathcal{E}}_0\check{\hat{Z}}_3\big)dt\\
       &\ +\sum\limits_{i=1}^3\big(\mathcal{A}_iX+\widehat{\mathcal{A}}_i\hat{X}+\overline{\mathcal{A}}_i\check{\hat{X}}+\mathcal{B}_iY+\widetilde{\mathcal{B}}_i\check{Y}
        +\overline{\mathcal{B}}_i\check{\hat{Y}}+\mathcal{B}_iZ_1+\widetilde{\mathcal{C}}_i\check{Z}_1\\
       &\ +\overline{\mathcal{C}}_i\check{\hat{Z}}_1+\mathcal{D}_iZ_2+\widetilde{\mathcal{D}}_i\check{Z}_2+\overline{\mathcal{D}}_i\check{\hat{Z}}_2
        +\mathcal{E}_iZ_3+\widetilde{\mathcal{E}}_i\check{Z}_3+\overline{\mathcal{E}}_i\check{\hat{Z}}_3\big)dW_i(t),\\
-dY(t)=&\big(\mathcal{Q}_2X+\mathcal{H}_1\check{\hat{X}}+\mathcal{A}_0Y+\mathcal{H}_2\hat{Y}+\overline{\mathcal{A}}_0^\top\check{\hat{Y}}
        +\mathcal{A}_1Z_1+\mathcal{H}_3\hat{Z}_1+\overline{\mathcal{A}}_1^\top\check{\hat{Z}}_1\\
       &\ +\mathcal{A}_2Z_2+\hat{\mathcal{A}}_2\hat{Z}_2+\overline{\mathcal{A}}_2^\top\check{\hat{Z}}_2
        +\mathcal{A}_3Z_3+\hat{\mathcal{A}}_3\hat{Z}_3+\overline{\mathcal{A}}_3^\top\check{\hat{Z}}_3\big)dt-Z_1dW_1(t)\\
       &-Z_2dW_2(t)-Z_3dW_3(t),\ t\in[0,T],\\
  X(0)=&\ X_0,\quad  Y(T)=\mathcal{G}_2X(T),
\end{aligned}
\right.
\end{equation}
where
{\footnotesize\begin{equation*}
\left\{
\begin{aligned}
&X_0\triangleq\left(\begin{array}{c}x_0\\0\end{array}\right),\ \mathcal{Q}_2\triangleq\left(\begin{array}{cc}Q_2&0\\0&0\end{array}\right),\ \mathcal{G}_2\triangleq\left(\begin{array}{cc}G_2&0\\0&0\end{array}\right),\
 \mathcal{H}_1\triangleq\left(\begin{array}{cc}0&0\\0&-L_4N_2^{-1}L_4^\top\end{array}\right),\\
&\mathcal{H}_2\triangleq\left(\begin{array}{cc}L_{01}&0\\0&L_0\end{array}\right),\ \mathcal{H}_3\triangleq\left(\begin{array}{cc}L_{11}&0\\0&L_1\end{array}\right),\
 \mathcal{A}_0\triangleq\left(\begin{array}{cc}A_0&0\\0&L_0\end{array}\right),\ \widehat{\mathcal{A}}_0\triangleq\left(\begin{array}{cc}L_{01}&0\\0&0\end{array}\right),\\
&\overline{\mathcal{A}}_0\triangleq\left(\begin{array}{cc}0&-(C_0+L_{05})N_2^{-1}L_4^\top\\0&0\end{array}\right),\ \mathcal{B}_0\triangleq\left(\begin{array}{cc}0&L_{02}\\L_{02}&0\end{array}\right),\ \mathcal{C}_0\triangleq\left(\begin{array}{cc}-C_0N_2^{-1}C_0^\top&0\\0&0\end{array}\right),\\
&\widetilde{\mathcal{C}}_0\triangleq\left(\begin{array}{cc}-(C_0+L_{05})N_2^{-1}(C_0+L_{05})^\top&0\\0&0\end{array}\right),\
 \overline{\mathcal{C}}_0\triangleq\left(\begin{array}{cc}-C_0N_2^{-1}L_{15}^\top-L_{05}N_2^{-1}(C_1+L_{15})^\top&0\\0&0\end{array}\right),\\
&\overline{\mathcal{D}}_0\triangleq\left(\begin{array}{cc}-C_0N_2^{-1}L_{25}^\top-L_{05}N_2^{-1}(C_2+L_{25})^\top&0\\0&0\end{array}\right),\
 \overline{\mathcal{E}}_0\triangleq\left(\begin{array}{cc}-C_0N_2^{-1}L_{35}^\top-L_{05}N_2^{-1}(C_3+L_{35})^\top&0\\0&0\end{array}\right),\\
&\mathcal{A}_1\triangleq\left(\begin{array}{cc}A_1&0\\0&L_1\end{array}\right),\
 \widehat{\mathcal{A}}_1\triangleq\left(\begin{array}{cc}L_{11}&0\\0&0\end{array}\right),\ \overline{\mathcal{A}}_1\triangleq\left(\begin{array}{cc}0&-(C_1+L_{15})N_2^{-1}L_4^\top\\0&0\end{array}\right),\
 \mathcal{B}_1\triangleq\left(\begin{array}{cc}0&L_{12}\\L_{03}&0\end{array}\right),\\
&\widetilde{\mathcal{B}}_1\triangleq\left(\begin{array}{cc}-C_1N_2^{-1}C_0^\top&0\\0&0\end{array}\right),\
 \overline{\mathcal{B}}_1\triangleq\left(\begin{array}{cc}-(C_1+L_{15})N_2^{-1}(C_0+L_{05})^\top&0\\0&0\end{array}\right),\ \mathcal{C}_1\triangleq\left(\begin{array}{cc}0&L_{13}\\L_{13}&0\end{array}\right),\\
&\widetilde{\mathcal{C}}_1\triangleq\left(\begin{array}{cc}-C_1N_2^{-1}C_1^\top&0\\0&0\end{array}\right),\
 \overline{\mathcal{C}}_1\triangleq\left(\begin{array}{cc}-C_1N_2^{-1}L_{15}-L_{15}N_2^{-1}(C_1+L_{15})^\top&0\\0&0\end{array}\right),\ \mathcal{D}_1\triangleq\left(\begin{array}{cc}0&0\\L_{23}&0\end{array}\right),\\
&\widetilde{\mathcal{D}}_1\triangleq\left(\begin{array}{cc}-C_1N_2^{-1}C_2^\top&0\\0&0\end{array}\right),\
 \overline{\mathcal{D}}_1\triangleq\left(\begin{array}{cc}-C_1N_2^{-1}L_{25}-L_{15}N_2^{-1}(C_2+L_{25})^\top&0\\0&0\end{array}\right),\\
&\mathcal{E}_1\triangleq\left(\begin{array}{cc}0&L_{14}\\L_{33}&0\end{array}\right),\
 \widetilde{\mathcal{E}}_1\triangleq\left(\begin{array}{cc}-C_1N_2^{-1}C_3^\top&0\\0&0\end{array}\right),\
 \overline{\mathcal{E}}_1\triangleq\left(\begin{array}{cc}-C_1N_2^{-1}L_{35}-L_{15}N_2^{-1}(C_3+L_{35})^\top&0\\0&0\end{array}\right),\\
&\mathcal{A}_2\triangleq\left(\begin{array}{cc}A_2&0\\0&0\end{array}\right),\
 \widehat{\mathcal{A}}_2\triangleq\left(\begin{array}{cc}L_{21}&0\\0&0\end{array}\right),\ \overline{\mathcal{A}}_2\triangleq\left(\begin{array}{cc}0&-(C_2+L_{25})N_2^{-1}L_4^\top\\0&0\end{array}\right),\
 \mathcal{B}_2\triangleq\left(\begin{array}{cc}0&L_{22}\\0&0\end{array}\right),\\
\end{aligned}
\right.
\end{equation*}}
and
{\footnotesize\begin{equation*}
\left\{
\begin{aligned}
&\widetilde{\mathcal{B}}_2\triangleq\left(\begin{array}{cc}-C_2N_2^{-1}C_0^\top&0\\0&0\end{array}\right),\
 \overline{\mathcal{B}}_2\triangleq\left(\begin{array}{cc}-(C_2+L_{25})N_2^{-1}(C_0+L_{05})^\top&0\\0&0\end{array}\right),\ \mathcal{C}_2\triangleq\left(\begin{array}{cc}0&L_{23}\\0&0\end{array}\right),\\
&\widetilde{\mathcal{C}}_2\triangleq\left(\begin{array}{cc}-C_2N_2^{-1}C_1^\top&0\\0&0\end{array}\right),\
 \overline{\mathcal{C}}_2\triangleq\left(\begin{array}{cc}-C_2N_2^{-1}L_{15}-L_{25}N_2^{-1}(C_1+L_{15})^\top&0\\0&0\end{array}\right),\ \mathcal{D}_2\equiv0,\\
&\widetilde{\mathcal{D}}_2\triangleq\left(\begin{array}{cc}-C_2N_2^{-1}C_2^\top&0\\0&0\end{array}\right),\
 \overline{\mathcal{D}}_2\triangleq\left(\begin{array}{cc}-C_2N_2^{-1}L_{25}-L_{25}N_2^{-1}(C_2+L_{25})^\top&0\\0&0\end{array}\right),\\
&\mathcal{E}_2\triangleq\left(\begin{array}{cc}0&L_{24}\\0&0\end{array}\right),\
 \widetilde{\mathcal{E}}_2\triangleq\left(\begin{array}{cc}-C_2N_2^{-1}C_3^\top&0\\0&0\end{array}\right),\
 \overline{\mathcal{E}}_2\triangleq\left(\begin{array}{cc}-C_2N_2^{-1}L_{35}-L_{25}N_2^{-1}(C_3+L_{35})^\top&0\\0&0\end{array}\right),\\
&\mathcal{A}_3\triangleq\left(\begin{array}{cc}A_3&0\\0&L_3\end{array}\right),\
 \widehat{\mathcal{A}}_3\triangleq\left(\begin{array}{cc}L_{31}&0\\0&0\end{array}\right),\ \overline{\mathcal{A}}_3\triangleq\left(\begin{array}{cc}0&-(C_3+L_{35})N_2^{-1}L_4^\top\\0&0\end{array}\right),\
 \mathcal{B}_3\triangleq\left(\begin{array}{cc}0&L_{32}\\L_{04}&0\end{array}\right),\\
&\widetilde{\mathcal{B}}_3\triangleq\left(\begin{array}{cc}-C_3N_2^{-1}C_0^\top&0\\0&0\end{array}\right),\
 \overline{\mathcal{B}}_3\triangleq\left(\begin{array}{cc}-(C_3+L_{35})N_2^{-1}(C_0+L_{05})^\top&0\\0&0\end{array}\right),\ \mathcal{C}_3\triangleq\left(\begin{array}{cc}0&L_{33}\\L_{14}&0\end{array}\right),\\
&\widetilde{\mathcal{C}}_3\triangleq\left(\begin{array}{cc}-C_3N_2^{-1}C_1^\top&0\\0&0\end{array}\right),\
 \overline{\mathcal{C}}_3\triangleq\left(\begin{array}{cc}-C_3N_2^{-1}L_{15}-L_{35}N_2^{-1}(C_1+L_{15})^\top&0\\0&0\end{array}\right),\ \mathcal{D}_3\triangleq\left(\begin{array}{cc}0&0\\L_{24}&0\end{array}\right),\\
&\widetilde{\mathcal{D}}_3\triangleq\left(\begin{array}{cc}-C_3N_2^{-1}C_2^\top&0\\0&0\end{array}\right),\
 \overline{\mathcal{D}}_3\triangleq\left(\begin{array}{cc}-C_3N_2^{-1}L_{25}-L_{35}N_2^{-1}(C_2+L_{25})^\top&0\\0&0\end{array}\right),\\
&\mathcal{E}_3\triangleq\left(\begin{array}{cc}0&L_{34}\\L_{34}&0\end{array}\right),\
 \widetilde{\mathcal{E}}_3\triangleq\left(\begin{array}{cc}-C_3N_2^{-1}C_3^\top&0\\0&0\end{array}\right),\
 \overline{\mathcal{E}}_3\triangleq\left(\begin{array}{cc}-C_3N_2^{-1}L_{35}-L_{35}N_2^{-1}(C_3+L_{35})^\top&0\\0&0\end{array}\right).
\end{aligned}
\right.
\end{equation*}}
And (\ref{optimal control of the leader}) can be written as
\begin{equation}\label{optimal control of the leader-2 dim}
\begin{aligned}
u_2^*(t)=-N_2^{-1}\bigg(\mathcal{L}_4^\top\check{\hat{X}}(t)+\mathcal{C}_{05}^\top\check{Y}(t)+\mathcal{L}_{05}^\top\check{\hat{Y}}(t)
+\sum\limits_{i=1}^3\mathcal{C}_{i5}^\top\check{Z}_i(t)+\sum\limits_{i=1}^3\mathcal{L}_{i5}^\top\check{\hat{Z}}_i(t)\bigg),
\end{aligned}
\end{equation}
where
\begin{equation*}
\mathcal{L}_4\triangleq\left(\begin{array}{c}0\\L_4\end{array}\right),\ \mathcal{C}_{i5}\triangleq\left(\begin{array}{c}C_i\\0\end{array}\right),\
 \mathcal{L}_{i5}\triangleq\left(\begin{array}{c}L_{i5}\\0\end{array}\right),\ i=0,1,2,3.
\end{equation*}

We wish to decouple FBSDE (\ref{optimality system-leader-2 dim}). For this target, let
\begin{equation}\label{relation of X and Y}
Y(t)=\mathcal{P}_1(t)X(t)+\mathcal{P}_2(t)\hat{X}(t)+\mathcal{P}_3(t)\check{X}(t)+\mathcal{P}_4(t)\check{\hat{X}}(t),
\end{equation}
where $\mathcal{P}_1(\cdot),\mathcal{P}_2(\cdot),\mathcal{P}_3(\cdot),\mathcal{P}_4(\cdot)$ are all differentiable, deterministic $\mathbb{R}^{2n}\times\mathbb{R}^{2n}$ matrix-valued functions with $\mathcal{P}_1(T)=\mathcal{G}_2$, $\mathcal{P}_2(T)=0$, $\mathcal{P}_3(T)=0$, $\mathcal{P}_4(T)=0$. From (\ref{optimality system-leader-2 dim}), the equations for $\hat{X}(\cdot)$, $\check{X}(\cdot)$, $\check{\hat{X}}(\cdot)$ are
\begin{equation}\label{hat X}
\left\{
\begin{aligned}
d\hat{X}(t)&=\Big[(\mathcal{A}_0+\widehat{\mathcal{A}}_0)\hat{X}+\overline{\mathcal{A}}_0\check{\hat{X}}+\mathcal{B}_0\hat{Y}+(\mathcal{C}_0+\widetilde{\mathcal{C}}_0)\check{\hat{Y}}
            +\mathcal{B}_1^\top\hat{Z}_1+(\widetilde{\mathcal{B}}_1^\top+\overline{\mathcal{C}}_0)\check{\hat{Z}}_1\\
           &\qquad+\mathcal{B}_2^\top\hat{Z}_2+(\widetilde{\mathcal{B}}_2^\top+\overline{\mathcal{D}}_0)\check{\hat{Z}}_2+\mathcal{B}_3^\top \hat{Z}_3+(\widetilde{\mathcal{B}}_3^\top+\overline{\mathcal{E}}_0)\check{\hat{Z}}_3\Big]dt+\sum\limits_{i=1,3}\Big[(\mathcal{A}_i+\widehat{\mathcal{A}}_i)\hat{X}\\
           &\qquad+(\widetilde{\mathcal{A}}_i+\overline{\mathcal{A}}_i)\check{\hat{X}}
            +\mathcal{B}_i\hat{Y}+(\widetilde{\mathcal{B}}_i+\overline{\mathcal{B}}_i)\check{\hat{Y}}+\mathcal{C}_i\hat{Z}_1+(\widetilde{\mathcal{C}}_i+\overline{\mathcal{C}}_i)\check{\hat{Z}}_1\\
           &\qquad+\mathcal{D}_i\hat{Z}_2+(\widetilde{\mathcal{D}}_i+\overline{\mathcal{D}}_i)\check{\hat{Z}}_2
            +\mathcal{E}_i\hat{Z}_3+(\widetilde{\mathcal{E}}_i+\overline{\mathcal{E}}_i)\check{\hat{Z}}_3\Big]dW_i(t),\ t\in[0,T],\\
\hat{X}(0)=&\ X_0,
\end{aligned}
\right.
\end{equation}
\begin{equation}\label{check-X}
\left\{
\begin{aligned}
d\check{X}(t)&=\Big[\mathcal{A}_0\check{X}+(\widehat{\mathcal{A}}_0+\overline{\mathcal{A}}_0)\check{\hat{X}}+(\mathcal{B}_0+\mathcal{C}_0)\check{Y}
              +\widetilde{\mathcal{C}}_0\check{\hat{Y}}+(\mathcal{B}_1+\widetilde{\mathcal{B}}_1)^\top\check{Z}_1+\overline{\mathcal{C}}_0\check{\hat{Z}}_1\\
             &\qquad+(\mathcal{B}_2+\widetilde{\mathcal{B}}_2)^\top\check{Z}_2+\overline{\mathcal{D}}_0\check{\hat{Z}}_2
              +(\mathcal{B}_3+\widetilde{\mathcal{B}}_3)^\top\check{Z}_3+\overline{\mathcal{E}}_0\check{\hat{Z}}_3\Big]dt\\
             &\quad+\sum\limits_{i=2,3}\Big[(\mathcal{A}_i+\widetilde{\mathcal{A}}_i)\check{X}+(\widehat{\mathcal{A}}_i+\overline{\mathcal{A}}_i)\check{\hat{X}}
              +(\mathcal{B}_i+\widetilde{\mathcal{B}}_i)\check{Y}+\overline{\mathcal{B}}_i\check{\hat{Y}}
              +(\mathcal{C}_i+\widetilde{\mathcal{C}}_i)\check{Z}_1\\
             &\qquad+\overline{\mathcal{C}}_i\check{\hat{Z}}_1+(\mathcal{D}_i+\widetilde{\mathcal{D}}_i\big)\check{Z}_2+\overline{\mathcal{D}}_i\check{\hat{Z}}_2
              +(\mathcal{E}_i+\widetilde{\mathcal{E}}_i)\check{Z}_3+\overline{\mathcal{E}}_i\check{\hat{Z}}_3\Big]dW_i(t),\ t\in[0,T],\\
\check{X}(0)=&\ X_0,
\end{aligned}
\right.
\end{equation}
and
\begin{equation}\label{hat and check-X}
\left\{
\begin{aligned}
d\check{\hat{X}}(t)&=\Big[(\mathcal{A}_0+\widehat{\mathcal{A}}_0+\overline{\mathcal{A}}_0)\check{\hat{X}}+(\mathcal{B}_0+\mathcal{C}_0+\widetilde{\mathcal{C}}_0)\check{\hat{Y}}
                    +(\mathcal{B}_1^\top+\widetilde{\mathcal{B}}_1^\top+\overline{\mathcal{C}}_0)\check{\hat{Z}}_1\\
                   &\qquad+(\mathcal{B}_2^\top+\widetilde{\mathcal{B}}_2^\top+\overline{\mathcal{D}}_0)\check{\hat{Z}}_2
                    +(\mathcal{B}_3^\top+\widetilde{\mathcal{B}}_3+\overline{\mathcal{E}}_0)\check{\hat{Z}}_3\Big]dt\\
                   &\quad+\Big[(\mathcal{A}_3+\widetilde{\mathcal{A}}_3+\widehat{\mathcal{A}}_3+\overline{\mathcal{A}}_3)\check{\hat{X}}
                    +(\mathcal{B}_3+\widetilde{\mathcal{B}}_3+\overline{\mathcal{B}}_3)\check{\hat{Y}}
                    +(\mathcal{C}_3+\widetilde{\mathcal{C}}_3+\overline{\mathcal{C}}_3)\check{\hat{Z}}_1\\
                   &\qquad+(\mathcal{D}_3+\widetilde{\mathcal{D}}_3+\overline{\mathcal{D}}_3)\check{\hat{Z}}_2
                    +(\mathcal{E}_3+\widetilde{\mathcal{E}}_3+\overline{\mathcal{E}}_3)\check{\hat{Z}}_3\Big]dW_3(t),\ t\in[0,T],\\
\check{\hat{X}}(0)=&\ X_0.
\end{aligned}
\right.
\end{equation}
Applying It\^{o}'s formula to (\ref{relation of X and Y}), we obtain
\begin{equation*}
\begin{aligned}
      dY(t)&=\Big\{\big(\dot{\mathcal{P}}_1+\mathcal{P}_1\mathcal{A}_0+\mathcal{P}_1\mathcal{B}_0\mathcal{P}_1\big)X
            +\big[\dot{\mathcal{P}}_2+\mathcal{P}_1\widehat{\mathcal{A}}_0+\mathcal{P}_2\big(\mathcal{A}_0+\widehat{\mathcal{A}}_0\big)+\mathcal{P}_1\mathcal{B}_0\mathcal{P}_2\\
           &\qquad+\mathcal{P}_2\mathcal{B}_0\mathcal{P}_1+\mathcal{P}_2\mathcal{B}_0\mathcal{P}_2\big]\hat{X}+\big[\dot{\mathcal{P}}_3+\mathcal{P}_3\mathcal{A}_0
            +\mathcal{P}_1\mathcal{B}_0\mathcal{P}_3+\mathcal{P}_1\mathcal{C}_0(\mathcal{P}_1+\mathcal{P}_2)\\
           &\qquad+\mathcal{P}_3\big(\mathcal{B}_0+\mathcal{C}_0\big)\big(\mathcal{P}_1+\mathcal{P}_3\big)\big]\check{X}
            +\big[\dot{\mathcal{P}}_4+\mathcal{P}_4\big(\mathcal{A}_0+\widehat{\mathcal{A}}_0+\overline{\mathcal{A}}_0\big)+\mathcal{P}_1\overline{\mathcal{A}}_0+\mathcal{P}_2\overline{\mathcal{A}}_0\\
           &\qquad+\mathcal{P}_4\big(\mathcal{B}_0+\mathcal{C}_0+\widetilde{\mathcal{C}}_0\big)\big(\mathcal{P}_1+\mathcal{P}_2+\mathcal{P}_3+\mathcal{P}_4\big)
            +\mathcal{P}_1\mathcal{B}_0\mathcal{P}_4+\mathcal{P}_1\mathcal{C}_0\big(\mathcal{P}_3+\mathcal{P}_4\big)\\
           &\qquad+\mathcal{P}_1\widetilde{\mathcal{C}}_0\big(\mathcal{P}_1+\mathcal{P}_2+\mathcal{P}_3+\mathcal{P}_4\big)
            +\mathcal{P}_2\big(\mathcal{C}_0+\widetilde{\mathcal{C}}_0\big)\big(\mathcal{P}_1+\mathcal{P}_2+\mathcal{P}_3+\mathcal{P}_4\big)
            +\mathcal{P}_2\mathcal{B}_0\big(\mathcal{P}_3+\mathcal{P}_4\big)\\
           &\qquad+\mathcal{P}_3\big(\widehat{\mathcal{A}}_0+\overline{\mathcal{A}}_0\big)+\mathcal{P}_3\big(\mathcal{B}_0+\mathcal{C}_0\big)\big(\mathcal{P}_2+\mathcal{P}_4\big)
            +\mathcal{P}_3\widetilde{\mathcal{C}}_0\big(\mathcal{P}_1+\mathcal{P}_2+\mathcal{P}_3+\mathcal{P}_4\big)\big]\check{\hat{X}}\\
           &\qquad+\mathcal{P}_1\mathcal{B}_1^\top Z_1+\mathcal{P}_2\mathcal{B}_1^\top\hat{Z}_1
            +\big(\mathcal{P}_1\widetilde{\mathcal{B}}_1^\top+\mathcal{P}_3\mathcal{B}_1^\top+\mathcal{P}_3\widetilde{\mathcal{B}}_1^\top\big)\check{Z}_1
            +\mathcal{P}_1\mathcal{B}_2^\top Z_2+\mathcal{P}_2\mathcal{B}_2^\top\hat{Z}_2\\
           &\qquad+\big(\mathcal{P}_1\widetilde{\mathcal{B}}_2^\top+\mathcal{P}_3\mathcal{B}_2^\top+\mathcal{P}_3\widetilde{\mathcal{B}}_2^\top\big)\check{Z}_2
            +\mathcal{P}_1\mathcal{B}_3^\top Z_3+\mathcal{P}_2\mathcal{B}_3^\top\hat{Z}_3
            +\big(\mathcal{P}_1\widetilde{\mathcal{B}}_3^\top+\mathcal{P}_3\mathcal{B}_3^\top+\mathcal{P}_3\widetilde{\mathcal{B}}_3^\top\big)\check{Z}_3\\
           &\qquad+\big[\mathcal{P}_1\overline{\mathcal{C}}_0+\mathcal{P}_2\big(\widetilde{\mathcal{B}}_1^\top+\overline{\mathcal{C}}_0\big)
            +\mathcal{P}_3\overline{\mathcal{C}}_0+\mathcal{P}_4\big(\mathcal{B}_1^\top+\widetilde{\mathcal{B}}_1^\top+\overline{\mathcal{C}}_0\big)\big]\check{\hat{Z}}_1\\
           &\qquad+\big[\mathcal{P}_1\overline{\mathcal{D}}_0+\mathcal{P}_2\big(\widetilde{\mathcal{B}}_2^\top+\overline{\mathcal{D}}_0\big)
            +\mathcal{P}_3\overline{\mathcal{D}}_0+\mathcal{P}_4\big(\mathcal{B}_2^\top+\widetilde{\mathcal{B}}_2^\top+\overline{\mathcal{D}}_0\big)\big]\check{\hat{Z}}_2\\
           &\qquad+\big[\mathcal{P}_1\overline{\mathcal{E}}_0+\mathcal{P}_2\big(\widetilde{\mathcal{B}}_3^\top+\overline{\mathcal{E}}_0\big)
            +\mathcal{P}_3\overline{\mathcal{E}}_0+\mathcal{P}_4\big(\mathcal{B}_3^\top+\widetilde{\mathcal{B}}_3^\top+\overline{\mathcal{E}}_0\big)\big]\check{\hat{Z}}_3\Big\}dt\\
           &\quad+\Big\{\big(\mathcal{P}_1\mathcal{A}_1+\mathcal{P}_1\mathcal{B}_1\mathcal{P}_1\big)X
            +\big[\mathcal{P}_1\widehat{\mathcal{A}}_1+\mathcal{P}_2\big(\mathcal{A}_1+\widehat{\mathcal{A}}_1\big)
            +\mathcal{P}_1\mathcal{B}_1\mathcal{P}_2+\mathcal{P}_2\mathcal{B}_1\big(\mathcal{P}_1+\mathcal{P}_2\big)\big]\hat{X}\\
           &\qquad+\big[\mathcal{P}_1\mathcal{B}_1\mathcal{P}_3+\mathcal{P}_1\widetilde{\mathcal{B}}_1\big(\mathcal{P}_1+\mathcal{P}_3\big)\big]\check{X}
            +\big[\mathcal{P}_1\overline{\mathcal{A}}_1+\mathcal{P}_1\mathcal{B}_1\mathcal{P}_4
            +\mathcal{P}_1\widetilde{\mathcal{B}}_1\big(\mathcal{P}_2+\mathcal{P}_4\big)\\
           &\qquad+\mathcal{P}_1\overline{\mathcal{B}}_1\big(\mathcal{P}_1+\mathcal{P}_2+\mathcal{P}_3+\mathcal{P}_4\big)
            +\mathcal{P}_2\big(\widetilde{\mathcal{A}}_1+\overline{\mathcal{A}}_1\big)+\mathcal{P}_2\mathcal{B}_1\big(\mathcal{P}_3+\mathcal{P}_4\big)\\
           &\qquad+\mathcal{P}_2\big(\widetilde{\mathcal{B}}_1+\overline{\mathcal{B}}_1\big)\big(\mathcal{P}_1+\mathcal{P}_2+\mathcal{P}_3+\mathcal{P}_4\big)\big]\check{\hat{X}}
            +\mathcal{P}_1\mathcal{B}_1Z_1+\mathcal{P}_2\mathcal{C}_1\hat{Z}_1+\mathcal{P}_1\widetilde{\mathcal{C}}_1\check{Z}_1+\mathcal{P}_1\mathcal{D}_1Z_2\\
           &\qquad+\mathcal{P}_2\mathcal{D}_1\hat{Z}_2+\mathcal{P}_1\widetilde{\mathcal{D}}_1\check{Z}_2
            +\mathcal{P}_1\mathcal{E}_1Z_3+\mathcal{P}_2\mathcal{E}_1\hat{Z}_3+\mathcal{P}_1\widetilde{\mathcal{E}}_1\check{Z}_3
            +\big(\mathcal{P}_1\overline{\mathcal{C}}_1+\mathcal{P}_2\widetilde{\mathcal{C}}_1+\mathcal{P}_2\overline{\mathcal{C}}_1\big)\check{\hat{Z}}_1\\
\end{aligned}
\end{equation*}
\begin{equation}\label{Applying Ito's formula to Y}
\begin{aligned}
           &\qquad+\big(\mathcal{P}_1\overline{\mathcal{D}}_1+\mathcal{P}_2\widetilde{\mathcal{D}}_1+\mathcal{P}_2\overline{\mathcal{D}}_1\big)\check{\hat{Z}}_2
            +\big(\mathcal{P}_1\overline{\mathcal{E}}_1+\mathcal{P}_2\widetilde{\mathcal{E}}_1+\mathcal{P}_2\overline{\mathcal{E}}_1\big)\check{\hat{Z}}_3\Big\}dW_1(t)\\
           &\quad+\Big\{\big(\mathcal{P}_1\mathcal{A}_2+\mathcal{P}_1\mathcal{B}_2\mathcal{P}_1\big)X
            +\big(\mathcal{P}_1\widehat{\mathcal{A}}_2+\mathcal{P}_1\mathcal{B}_2\mathcal{P}_2\big)\hat{X}
            +\big[\mathcal{P}_1\mathcal{B}_2\mathcal{P}_3+\mathcal{P}_1\widetilde{\mathcal{B}}_2\big(\mathcal{P}_1+\mathcal{P}_3\big)\\
           &\qquad+\mathcal{P}_3\big(\mathcal{A}_2+\widetilde{\mathcal{A}}_2\big)+\mathcal{P}_3\big(\mathcal{B}_2+\widetilde{\mathcal{B}}_2\big)
            \big(\mathcal{P}_1+\mathcal{P}_3\big)\big]\check{X}+\big[\mathcal{P}_1\overline{\mathcal{A}}_2+\mathcal{P}_1\mathcal{B}_2\mathcal{P}_4
            +\mathcal{P}_1\widetilde{\mathcal{B}}_2\big(\mathcal{P}_2+\mathcal{P}_4\big)\\
           &\qquad+\mathcal{P}_1\overline{\mathcal{B}}_2\big(\mathcal{P}_1+\mathcal{P}_2+\mathcal{P}_3+\mathcal{P}_4\big)
            +\mathcal{P}_3\big(\widehat{\mathcal{A}}_2+\overline{\mathcal{A}}_2\big)
            +\mathcal{P}_3\big(\mathcal{B}_2+\widetilde{\mathcal{B}}_2\big)\big(\mathcal{P}_2+\mathcal{P}_4\big)\\
           &\qquad+\mathcal{P}_3\widetilde{\mathcal{B}}_2\big(\mathcal{P}_1+\mathcal{P}_2+\mathcal{P}_3+\mathcal{P}_4\big)\big]\check{\hat{X}}
            +\mathcal{P}_1\mathcal{B}_2Z_1+\big(\mathcal{P}_1\widetilde{\mathcal{C}}_2+\mathcal{P}_3\mathcal{C}_2
            +\mathcal{P}_3\widetilde{\mathcal{C}}_2\big)\check{Z}_1+\mathcal{P}_1\mathcal{D}_2Z_2\\
           &\qquad+\big(\mathcal{P}_1\widetilde{\mathcal{D}}_2+\mathcal{P}_3\mathcal{D}_2+\mathcal{P}_3\widetilde{\mathcal{D}}_2\big)\check{Z}_2
            +\mathcal{P}_1\mathcal{E}_2Z_3+\big(\mathcal{P}_1\widetilde{\mathcal{E}}_2+\mathcal{P}_3\mathcal{E}_2
            +\mathcal{P}_3\widetilde{\mathcal{E}}_2\big)\check{Z}_3\\
           &\qquad+\big(\mathcal{P}_1\overline{\mathcal{C}}_2+\mathcal{P}_3\overline{\mathcal{C}}_2\big)\check{\hat{Z}}_1+\big(\mathcal{P}_1\overline{\mathcal{D}}_2+\mathcal{P}_3\overline{\mathcal{D}}_2\big)\check{\hat{Z}}_2
            +\big(\mathcal{P}_1\overline{\mathcal{E}}_2+\mathcal{P}_3\overline{\mathcal{E}}_2\big)\check{\hat{Z}}_3\Big\}dW_2(t)\\
           &\quad+\Big\{\big(\mathcal{P}_1\mathcal{A}_3+\mathcal{P}_1\mathcal{B}_3\mathcal{P}_1\big)X
            +\big[\mathcal{P}_1\widehat{\mathcal{A}}_3+\mathcal{P}_1\mathcal{B}_3\mathcal{P}_2
            +\mathcal{P}_2(\mathcal{A}_3+\widehat{\mathcal{A}}_3)+\mathcal{P}_2\mathcal{B}_3\big(\mathcal{P}_1+\mathcal{P}_2\big)\big]\hat{X}\\
           &\qquad+\big[\mathcal{P}_1\mathcal{B}_3\mathcal{P}_3+\mathcal{P}_1\widetilde{\mathcal{B}}_3\big(\mathcal{P}_1+\mathcal{P}_3\big)
            +\mathcal{P}_3\big(\mathcal{A}_3+\widetilde{\mathcal{A}}_3\big)+\mathcal{P}_3\big(\mathcal{B}_3+\widetilde{\mathcal{B}}_3\big)
            \big(\mathcal{P}_1+\mathcal{P}_3\big)\big]\check{X}\\
           &\qquad+\Big[\mathcal{P}_1\overline{\mathcal{A}}_3+\mathcal{P}_1\mathcal{B}_3\mathcal{P}_4+\mathcal{P}_1\widetilde{\mathcal{B}}_3\big(\mathcal{P}_2+\mathcal{P}_4\big)
            +\mathcal{P}_1\overline{\mathcal{B}}_3\big(\mathcal{P}_1+\mathcal{P}_2+\mathcal{P}_3+\mathcal{P}_4\big)+\mathcal{P}_2\big(\widetilde{\mathcal{A}}_3+\overline{\mathcal{A}}_3\big)\\
           &\qquad+\mathcal{P}_2\mathcal{B}_3\big(\mathcal{P}_3+\mathcal{P}_4\big)+\mathcal{P}_2\big(\widetilde{\mathcal{B}}_3+\overline{\mathcal{B}}_3\big)
            \big(\mathcal{P}_1+\mathcal{P}_2+\mathcal{P}_3+\mathcal{P}_4\big)+\mathcal{P}_3\big(\widehat{\mathcal{A}}_3+\overline{\mathcal{A}}_3\big)\\
           &\qquad+\mathcal{P}_3\big(\mathcal{B}_3+\widetilde{\mathcal{B}}_3\big)\big(\mathcal{P}_2+\mathcal{P}_4\big)
            +\mathcal{P}_3\overline{\mathcal{B}}_3\big(\mathcal{P}_1+\mathcal{P}_2+\mathcal{P}_3+\mathcal{P}_4\big)
            +\mathcal{P}_4\big(\mathcal{A}_3+\widetilde{\mathcal{A}}_3+\widehat{\mathcal{A}}_3+\overline{\mathcal{A}}_3\big)\\
           &\qquad+\mathcal{P}_4\big(\mathcal{B}_3+\widetilde{\mathcal{B}}_3+\overline{\mathcal{B}}_3\big)
            \big(\mathcal{P}_1+\mathcal{P}_2+\mathcal{P}_3+\mathcal{P}_4\big)\Big]\check{\hat{X}}
            +\mathcal{P}_1\mathcal{B}_3Z_1+\mathcal{P}_2\mathcal{C}_3\hat{Z}_1+\big(\mathcal{P}_1\widetilde{\mathcal{C}}_3+\mathcal{P}_3\mathcal{C}_3\\
           &\qquad+\mathcal{P}_3\widetilde{\mathcal{C}}_3\big)\check{Z}_1+\mathcal{P}_1\mathcal{D}_3Z_2+\mathcal{P}_2\mathcal{D}_3\hat{Z}_2
            +\big(\mathcal{P}_1\widetilde{\mathcal{D}}_3+\mathcal{P}_3\mathcal{D}_3+\mathcal{P}_3\widetilde{\mathcal{D}}_3\big)\check{Z}_2
            +\mathcal{P}_1\mathcal{E}_3Z_3+\mathcal{P}_2\mathcal{E}_3\hat{Z}_3\\
           &\qquad+\big(\mathcal{P}_1\widetilde{\mathcal{E}}_3+\mathcal{P}_3\mathcal{E}_3+\mathcal{P}_3\widetilde{\mathcal{E}}_3\big)\check{Z}_3
            +\big[\mathcal{P}_1\overline{\mathcal{C}}_3+\mathcal{P}_2\widetilde{\mathcal{C}}_3+\mathcal{P}_2\overline{\mathcal{C}}_3+\mathcal{P}_3\overline{\mathcal{C}}_3
            +\mathcal{P}_4\big(\mathcal{C}_3+\widetilde{\mathcal{C}}_3+\overline{\mathcal{C}}_3\big)\big]\check{\hat{Z}}_1\\
           &\qquad+\big[\mathcal{P}_1\overline{\mathcal{D}}_3+\mathcal{P}_2\widetilde{\mathcal{D}}_3+\mathcal{P}_2\overline{\mathcal{D}}_3+\mathcal{P}_3\overline{\mathcal{D}}_3
            +\mathcal{P}_4\big(\mathcal{D}_3+\widetilde{\mathcal{D}}_3+\overline{\mathcal{D}}_3\big)\big]\check{\hat{Z}}_2\\
           &\qquad+\big[\mathcal{P}_1\overline{\mathcal{E}}_3+\mathcal{P}_2\widetilde{\mathcal{E}}_3+\mathcal{P}_2\overline{\mathcal{E}}_3+\mathcal{P}_3\overline{\mathcal{E}}_3
            +\mathcal{P}_4\big(\mathcal{E}_3+\widetilde{\mathcal{E}}_3+\overline{\mathcal{E}}_3\big)\big]\check{\hat{Z}}_3\Big\}dW_3(t)\\
           &=-\Big\{(\mathcal{Q}_2+\mathcal{A}_0\mathcal{P}_1)X+(\mathcal{A}_0\mathcal{P}_2+\mathcal{H}_1\mathcal{P}_1+\mathcal{H}_1\mathcal{P}_2)\hat{X}
            +\mathcal{A}_0\mathcal{P}_3\check{X}+\big[\mathcal{H}_1+\mathcal{A}_0\mathcal{P}_4+\mathcal{H}_2\mathcal{P}_3\\
           &\qquad+\mathcal{H}_2\mathcal{P}_4+\overline{\mathcal{A}}_0^\top\big(\mathcal{P}_1+\mathcal{P}_2+\mathcal{P}_3+\mathcal{P}_4\big)\big]\check{\hat{X}}
            +\mathcal{A}_1Z_1+\overline{\mathcal{A}}_1^\top\check{\hat{Z}}_1+\mathcal{A}_2Z_2+\widehat{\mathcal{A}}_2\hat{Z}_2+\overline{\mathcal{A}}_2^\top\check{\hat{Z}}_2\\
           &\qquad+\mathcal{A}_3Z_3+\widehat{\mathcal{A}}_3\hat{Z}_3+\overline{\mathcal{A}}_3^\top\check{\hat{Z}}_3\Big\}dt+Z_1dW_1(t)+Z_2dW_2(t)+Z_3dW_3(t).
\end{aligned}
\end{equation}
Comparing the diffusion terms $dW_1(t),dW_2(t),dW_3(t)$ on both sides of (\ref{Applying Ito's formula to Y}) respectively, we have
\begin{equation*}
\begin{aligned}
     Z_1(t)&=\big(\mathcal{P}_1\mathcal{A}_1+\mathcal{P}_1\mathcal{B}_1\mathcal{P}_1\big)X
            +\big[\mathcal{P}_1\widehat{\mathcal{A}}_1+\mathcal{P}_2\big(\mathcal{A}_1+\widehat{\mathcal{A}}_1\big)
            +\mathcal{P}_1\mathcal{B}_1\mathcal{P}_2+\mathcal{P}_2\mathcal{B}_1\big(\mathcal{P}_1+\mathcal{P}_2\big)\big]\hat{X}\\
           &\quad+\big[\mathcal{P}_1\mathcal{B}_1\mathcal{P}_3+\mathcal{P}_1\widetilde{\mathcal{B}}_1\big(\mathcal{P}_1+\mathcal{P}_3\big)\big]\check{X}
            +\big[\mathcal{P}_1\overline{\mathcal{A}}_1+\mathcal{P}_1\mathcal{B}_1\mathcal{P}_4
            +\mathcal{P}_1\widetilde{\mathcal{B}}_1\big(\mathcal{P}_2+\mathcal{P}_4\big)\\
           &\quad+\mathcal{P}_1\overline{\mathcal{B}}_1\big(\mathcal{P}_1+\mathcal{P}_2+\mathcal{P}_3+\mathcal{P}_4\big)
            +\mathcal{P}_2\big(\widetilde{\mathcal{A}}_1+\overline{\mathcal{A}}_1\big)+\mathcal{P}_2\mathcal{B}_1\big(\mathcal{P}_3+\mathcal{P}_4\big)\\
           &\quad+\mathcal{P}_2\big(\widetilde{\mathcal{B}}_1+\overline{\mathcal{B}}_1\big)\big(\mathcal{P}_1+\mathcal{P}_2+\mathcal{P}_3+\mathcal{P}_4\big)\big]\check{\hat{X}}
            +\mathcal{P}_1\mathcal{B}_1Z_1+\mathcal{P}_2\mathcal{C}_1\hat{Z}_1+\mathcal{P}_1\widetilde{\mathcal{C}}_1\check{Z}_1+\mathcal{P}_1\mathcal{D}_1Z_2\\
           &\quad+\mathcal{P}_2\mathcal{D}_1\hat{Z}_2+\mathcal{P}_1\widetilde{\mathcal{D}}_1\check{Z}_2
            +\mathcal{P}_1\mathcal{E}_1Z_3+\mathcal{P}_2\mathcal{E}_1\hat{Z}_3+\mathcal{P}_1\widetilde{\mathcal{E}}_1\check{Z}_3
            +\big(\mathcal{P}_1\overline{\mathcal{C}}_1+\mathcal{P}_2\widetilde{\mathcal{C}}_1+\mathcal{P}_2\overline{\mathcal{C}}_1\big)\check{\hat{Z}}_1\\
           &\quad+\big(\mathcal{P}_1\overline{\mathcal{D}}_1+\mathcal{P}_2\widetilde{\mathcal{D}}_1+\mathcal{P}_2\overline{\mathcal{D}}_1\big)\check{\hat{Z}}_2
            +\big(\mathcal{P}_1\overline{\mathcal{E}}_1+\mathcal{P}_2\widetilde{\mathcal{E}}_1+\mathcal{P}_2\overline{\mathcal{E}}_1\big)\check{\hat{Z}}_3,\\
\end{aligned}
\end{equation*}
\begin{equation*}
\begin{aligned}
     Z_2(t)&=\big(\mathcal{P}_1\mathcal{A}_2+\mathcal{P}_1\mathcal{B}_2\mathcal{P}_1\big)X
            +\big(\mathcal{P}_1\widehat{\mathcal{A}}_2+\mathcal{P}_1\mathcal{B}_2\mathcal{P}_2\big)\hat{X}
            +\big[\mathcal{P}_1\mathcal{B}_2\mathcal{P}_3+\mathcal{P}_1\widetilde{\mathcal{B}}_2\big(\mathcal{P}_1+\mathcal{P}_3\big)\\
           &\quad+\mathcal{P}_3\big(\mathcal{A}_2+\widetilde{\mathcal{A}}_2\big)+\mathcal{P}_3\big(\mathcal{B}_2+\widetilde{\mathcal{B}}_2\big)
            \big(\mathcal{P}_1+\mathcal{P}_3\big)\big]\check{X}+\big[\mathcal{P}_1\overline{\mathcal{A}}_2+\mathcal{P}_1\mathcal{B}_2\mathcal{P}_4
            +\mathcal{P}_1\widetilde{\mathcal{B}}_2\big(\mathcal{P}_2+\mathcal{P}_4\big)\\
           &\quad+\mathcal{P}_1\overline{\mathcal{B}}_2\big(\mathcal{P}_1+\mathcal{P}_2+\mathcal{P}_3+\mathcal{P}_4\big)
            +\mathcal{P}_3\big(\widehat{\mathcal{A}}_2+\overline{\mathcal{A}}_2\big)
            +\mathcal{P}_3\big(\mathcal{B}_2+\widetilde{\mathcal{B}}_2\big)\big(\mathcal{P}_2+\mathcal{P}_4\big)\\
           &\quad+\mathcal{P}_3\overline{\mathcal{B}}_2\big(\mathcal{P}_1+\mathcal{P}_2+\mathcal{P}_3+\mathcal{P}_4\big)\big]\check{\hat{X}}
            +\mathcal{P}_1\mathcal{B}_2Z_1+\big(\mathcal{P}_1\widetilde{\mathcal{C}}_2+\mathcal{P}_3\mathcal{C}_2
            +\mathcal{P}_3\widetilde{\mathcal{C}}_2\big)\check{Z}_1+\mathcal{P}_1\mathcal{D}_2Z_2\\
           &\quad+\big(\mathcal{P}_1\widetilde{\mathcal{D}}_2+\mathcal{P}_3\mathcal{D}_2+\mathcal{P}_3\widetilde{\mathcal{D}}_2\big)\check{Z}_2
            +\mathcal{P}_1\mathcal{E}_2Z_3+\big(\mathcal{P}_1\widetilde{\mathcal{E}}_2+\mathcal{P}_3\mathcal{E}_2
            +\mathcal{P}_3\widetilde{\mathcal{E}}_2\big)\check{Z}_3\\
           &\quad+\big(\mathcal{P}_1\overline{\mathcal{C}}_2+\mathcal{P}_3\overline{\mathcal{C}}_2\big)\check{\hat{Z}}_1+\big(\mathcal{P}_1\overline{\mathcal{D}}_2+\mathcal{P}_3\overline{\mathcal{D}}_2\big)\check{\hat{Z}}_2
            +\big(\mathcal{P}_1\overline{\mathcal{E}}_2+\mathcal{P}_3\overline{\mathcal{E}}_2\big)\check{\hat{Z}}_3,\\
\end{aligned}
\end{equation*}
\begin{equation}\label{comparing dW1,dW2,dW3-leader}
\begin{aligned}     Z_3(t)&=\big(\mathcal{P}_1\mathcal{A}_3+\mathcal{P}_1\mathcal{B}_3\mathcal{P}_1\big)X
            +\big[\mathcal{P}_1\widehat{\mathcal{A}}_3+\mathcal{P}_1\mathcal{B}_3\mathcal{P}_2
            +\mathcal{P}_2(\mathcal{A}_3+\widehat{\mathcal{A}}_3)+\mathcal{P}_2\mathcal{B}_3\big(\mathcal{P}_1+\mathcal{P}_2\big)\big]\hat{X}\\
           &\quad+\big[\mathcal{P}_1\mathcal{B}_3\mathcal{P}_3+\mathcal{P}_1\widetilde{\mathcal{B}}_3\big(\mathcal{P}_1+\mathcal{P}_3\big)
            +\mathcal{P}_3\big(\mathcal{A}_3+\widetilde{\mathcal{A}}_3\big)+\mathcal{P}_3\big(\mathcal{B}_3+\widetilde{\mathcal{B}}_3\big)
            \big(\mathcal{P}_1+\mathcal{P}_3\big)\big]\check{X}\\
           &\quad+\Big[\mathcal{P}_1\overline{\mathcal{A}}_3+\mathcal{P}_1\mathcal{B}_3\mathcal{P}_4
            +\mathcal{P}_1\widetilde{\mathcal{B}}_3\big(\mathcal{P}_2+\mathcal{P}_4\big)
            +\mathcal{P}_1\overline{\mathcal{B}}_3\big(\mathcal{P}_1+\mathcal{P}_2+\mathcal{P}_3+\mathcal{P}_4\big)
            +\mathcal{P}_2\big(\widetilde{\mathcal{A}}_3+\overline{\mathcal{A}}_3\big)\\
           &\quad+\mathcal{P}_2\mathcal{B}_3\big(\mathcal{P}_3+\mathcal{P}_4\big)+\mathcal{P}_2\big(\widetilde{\mathcal{B}}_3+\overline{\mathcal{B}}_3\big)
            \big(\mathcal{P}_1+\mathcal{P}_2+\mathcal{P}_3+\mathcal{P}_4\big)+\mathcal{P}_3\big(\widehat{\mathcal{A}}_3+\overline{\mathcal{A}}_3\big)\\
           &\quad+\mathcal{P}_3\big(\mathcal{B}_3+\widetilde{\mathcal{B}}_3\big)\big(\mathcal{P}_2+\mathcal{P}_4\big)
            +\mathcal{P}_3\overline{\mathcal{B}}_3\big(\mathcal{P}_1+\mathcal{P}_2+\mathcal{P}_3+\mathcal{P}_4\big)
            +\mathcal{P}_4\big(\mathcal{A}_3+\widetilde{\mathcal{A}}_3+\widehat{\mathcal{A}}_3+\overline{\mathcal{A}}_3\big)\\
           &\quad+\mathcal{P}_4\big(\mathcal{B}_3+\widetilde{\mathcal{B}}_3+\overline{\mathcal{B}}_3\big)
            \big(\mathcal{P}_1+\mathcal{P}_2+\mathcal{P}_3+\mathcal{P}_4\big)\Big]\check{\hat{X}}
            +\mathcal{P}_1\mathcal{B}_3Z_1+\mathcal{P}_2\mathcal{C}_3\hat{Z}_1+\big(\mathcal{P}_1\widetilde{\mathcal{C}}_3+\mathcal{P}_3\mathcal{C}_3\\
           &\quad+\mathcal{P}_3\widetilde{\mathcal{C}}_3\big)\check{Z}_1+\mathcal{P}_1\mathcal{D}_3Z_2+\mathcal{P}_2\mathcal{D}_3\hat{Z}_2
            +\big(\mathcal{P}_1\widetilde{\mathcal{D}}_3+\mathcal{P}_3\mathcal{D}_3+\mathcal{P}_3\widetilde{\mathcal{D}}_3\big)\check{Z}_2
            +\mathcal{P}_1\mathcal{E}_3Z_3+\mathcal{P}_2\mathcal{E}_3\hat{Z}_3\\
           &\quad+\big(\mathcal{P}_1\widetilde{\mathcal{E}}_3+\mathcal{P}_3\mathcal{E}_3+\mathcal{P}_3\widetilde{\mathcal{E}}_3\big)\check{Z}_3
            +\big[\mathcal{P}_1\overline{\mathcal{C}}_3+\mathcal{P}_2\widetilde{\mathcal{C}}_3+\mathcal{P}_2\overline{\mathcal{C}}_3+\mathcal{P}_3\overline{\mathcal{C}}_3
            +\mathcal{P}_4\big(\mathcal{C}_3+\widetilde{\mathcal{C}}_3+\overline{\mathcal{C}}_3\big)\big]\check{\hat{Z}}_1\\
           &\quad+\big[\mathcal{P}_1\overline{\mathcal{D}}_3+\mathcal{P}_2\widetilde{\mathcal{D}}_3+\mathcal{P}_2\overline{\mathcal{D}}_3+\mathcal{P}_3\overline{\mathcal{D}}_3
            +\mathcal{P}_4\big(\mathcal{D}_3+\widetilde{\mathcal{D}}_3+\overline{\mathcal{D}}_3\big)\big]\check{\hat{Z}}_2\\
           &\quad+\big[\mathcal{P}_1\overline{\mathcal{E}}_3+\mathcal{P}_2\widetilde{\mathcal{E}}_3+\mathcal{P}_2\overline{\mathcal{E}}_3+\mathcal{P}_3\overline{\mathcal{E}}_3
            +\mathcal{P}_4\big(\mathcal{E}_3+\widetilde{\mathcal{E}}_3+\overline{\mathcal{E}}_3\big)\big]\check{\hat{Z}}_3.
\end{aligned}
\end{equation}

Next, we wish to represent each $Z_i(\cdot)$ and its filtering estimates as functionals of the ``state" $X(\cdot)$ and its filtering estimates, from (\ref{comparing dW1,dW2,dW3-leader}). For this target, we need the following four steps.

\vspace{1mm}

{\bf Step 1.}\quad Taking $\mathbb{E}\big[\mathbb{E}[\cdot|\mathcal{G}^2_t]\big|\mathcal{G}^1_t\big]$ on both sides of (\ref{comparing dW1,dW2,dW3-leader}), we derive
\begin{equation}\label{double filter estimate Zi}
\begin{aligned}
      \check{\hat{Z}}_i(t)=\mathcal{M}_{i0}\check{\hat{X}}(t)+\mathcal{M}_{i1}\check{\hat{Z}}_1(t)+\mathcal{M}_{i2}\check{\hat{Z}}_2(t)+\mathcal{M}_{i3}\check{\hat{Z}}_3(t),\ i=1,2,3,
\end{aligned}
\end{equation}
where
{\footnotesize\begin{equation*}
\left\{
\begin{aligned}
\mathcal{M}_{10}&\triangleq\mathcal{P}_1(\mathcal{A}_1+\widehat{\mathcal{A}}_1+\overline{\mathcal{A}}_1)
                 +\mathcal{P}_2(\mathcal{A}_1+\widehat{\mathcal{A}}_1+\widetilde{\mathcal{A}}_1+\overline{\mathcal{A}}_1)
                 +(\mathcal{P}_1+\mathcal{P}_2)(\mathcal{B}_1+\widetilde{\mathcal{B}}_1+\overline{\mathcal{B}}_1)\big(\mathcal{P}_1+\mathcal{P}_2+\mathcal{P}_3+\mathcal{P}_4\big),\\
\mathcal{M}_{20}&\triangleq\mathcal{P}_1(\mathcal{A}_2+\widehat{\mathcal{A}}_2+\overline{\mathcal{A}}_2)
                 +\mathcal{P}_3(\mathcal{A}_2+\widehat{\mathcal{A}}_2+\widetilde{\mathcal{A}}_2+\overline{\mathcal{A}}_2)
                 +(\mathcal{P}_1+\mathcal{P}_3)(\mathcal{B}_2+\widetilde{\mathcal{B}}_2+\overline{\mathcal{B}}_2)\big(\mathcal{P}_1+\mathcal{P}_2+\mathcal{P}_3+\mathcal{P}_4\big),\\
\mathcal{M}_{30}&\triangleq\mathcal{P}_1(\mathcal{A}_3+\widehat{\mathcal{A}}_3+\overline{\mathcal{A}}_3)
                 +(\mathcal{P}_2+\mathcal{P}_3+\mathcal{P}_4)(\mathcal{A}_3+\widehat{\mathcal{A}}_3+\widetilde{\mathcal{A}}_3+\overline{\mathcal{A}}_3)\\
                &\quad+(\mathcal{P}_1+\mathcal{P}_2+\mathcal{P}_3+\mathcal{P}_4)(\mathcal{B}_3+\widetilde{\mathcal{B}}_3+\overline{\mathcal{B}}_3)
                 \big(\mathcal{P}_1+\mathcal{P}_2+\mathcal{P}_3+\mathcal{P}_4\big),\\
\mathcal{M}_{11}&\triangleq\mathcal{P}_1\mathcal{B}_1+\mathcal{P}_2\mathcal{C}_1+\mathcal{P}_1\widetilde{\mathcal{C}}_1
                 +\mathcal{P}_1\overline{\mathcal{C}}_1+\mathcal{P}_2\widetilde{\mathcal{C}}_1+\mathcal{P}_2\overline{\mathcal{C}}_1,\
\mathcal{M}_{12}\triangleq\mathcal{P}_1\mathcal{D}_1+\mathcal{P}_2\mathcal{D}_1+\mathcal{P}_1\widetilde{\mathcal{D}}_1
                 +\mathcal{P}_1\overline{\mathcal{D}}_1+\mathcal{P}_2\widetilde{\mathcal{D}}_1+\mathcal{P}_2\overline{\mathcal{D}}_1,\\
\mathcal{M}_{13}&\triangleq\mathcal{P}_1\mathcal{E}_1+\mathcal{P}_2\mathcal{E}_1+\mathcal{P}_1\widetilde{\mathcal{E}}_1
                 +\mathcal{P}_1\overline{\mathcal{E}}_1+\mathcal{P}_2\widetilde{\mathcal{E}}_1+\mathcal{P}_2\overline{\mathcal{E}}_1,\
\mathcal{M}_{21}\triangleq\mathcal{P}_1\mathcal{B}_2+\mathcal{P}_1\widetilde{\mathcal{C}}_2+\mathcal{P}_3\mathcal{C}_2
                 +\mathcal{P}_3\widetilde{\mathcal{C}}_2+\mathcal{P}_1\overline{\mathcal{C}}_2+\mathcal{P}_3\overline{\mathcal{C}}_2,\\
\mathcal{M}_{22}&\triangleq\mathcal{P}_1\mathcal{D}_2+\mathcal{P}_1\widetilde{\mathcal{D}}_2+\mathcal{P}_3\mathcal{D}_2
                 +\mathcal{P}_3\widetilde{\mathcal{D}}_2+\mathcal{P}_1\overline{\mathcal{D}}_2+\mathcal{P}_3\overline{\mathcal{D}}_2,\
\mathcal{M}_{23}\triangleq\mathcal{P}_1\mathcal{E}_2+\mathcal{P}_1\widetilde{\mathcal{E}}_2+\mathcal{P}_3\mathcal{E}_2
                 +\mathcal{P}_3\widetilde{\mathcal{E}}_2+\mathcal{P}_1\overline{\mathcal{E}}_2+\mathcal{P}_3\overline{\mathcal{E}}_2,\\
\mathcal{M}_{31}&\triangleq\mathcal{P}_1\mathcal{B}_3+\mathcal{P}_1\widetilde{\mathcal{C}}_3+\mathcal{P}_1\overline{\mathcal{C}}_3
                 +(\mathcal{P}_2+\mathcal{P}_3+\mathcal{P}_4)(\mathcal{C}_3+\widetilde{\mathcal{C}}_3+\overline{\mathcal{C}}_3),\\
\mathcal{M}_{32}&\triangleq\mathcal{P}_1\mathcal{D}_3+\mathcal{P}_1\widetilde{\mathcal{D}}_3+\mathcal{P}_1\overline{\mathcal{D}}_3
                 +(\mathcal{P}_2+\mathcal{P}_3+\mathcal{P}_4)(\mathcal{D}_3+\widetilde{\mathcal{D}}_3+\overline{\mathcal{D}}_3),\\
\mathcal{M}_{33}&\triangleq\mathcal{P}_1\mathcal{E}_3+\mathcal{P}_1\widetilde{\mathcal{E}}_3+\mathcal{P}_1\overline{\mathcal{E}}_3
                 +(\mathcal{P}_2+\mathcal{P}_3+\mathcal{P}_4)(\mathcal{E}_3+\widetilde{\mathcal{E}}_3+\overline{\mathcal{E}}_3).
\end{aligned}
\right.
\end{equation*}}
We rewrite (\ref{double filter estimate Zi}) as
\begin{equation}\label{linear equation system-1}
\begin{aligned}
\left(
\begin{array}{ccc}
I_n-\mathcal{M}_{11}&-\mathcal{M}_{12}&-\mathcal{M}_{13}\\
-\mathcal{M}_{21}&I_n-\mathcal{M}_{22}&-\mathcal{M}_{23}\\
-\mathcal{M}_{31}&-\mathcal{M}_{32}&I_n-\mathcal{M}_{33}\\
\end{array}
\right)
\left(
\begin{array}{ccc}
\check{\hat{Z}}_1(t)\\\check{\hat{Z}}_2(t)\\\check{\hat{Z}}_3(t)
\end{array}
\right)=
\left(
\begin{array}{ccc}
\mathcal{M}_{10}\\\mathcal{M}_{20}\\\mathcal{M}_{30}
\end{array}
\right)\check{\hat{X}}(t).
\end{aligned}
\end{equation}
If we assume that

{\bf (A2.2)}\quad the coefficient matrix of (\ref{linear equation system-1}) is invertible, for any $t\in[0,T]$,

\noindent then by Cramer's rule, we have
\begin{equation}\label{double filter estimate Zi-G1G2}
\check{\hat{Z}}_i(t)=(-1)^{i-1}\big(\mathbb{N}_1\big)^{-1}\Big[\mathcal{M}_{10}\mathbf{M}_{1i}-\mathcal{M}_{20}\mathbf{M}_{2i}+\mathcal{M}_{30}\mathbf{M}_{3i}\Big]\check{\hat{X}}(t)
\triangleq\mathcal{N}^i(t)\check{\hat{X}}(t),\ i=1,2,3,
\end{equation}
where $\mathbb{N}_1$ is the determinant of the coefficient of (\ref{linear equation system-1}), and $\mathbf{M}^{ji}(t)$ is the adjoint matrix of the $(j,i)$ element in (\ref{linear equation system-1}), for $j,i=1,2,3$.

\vspace{1mm}

{\bf Step 2.}\quad Taking $\mathbb{E}[\cdot|\mathcal{G}^2_t]$ on both sides of (\ref{comparing dW1,dW2,dW3-leader}), we get
\begin{equation}\label{filter estimate Zi-G2}
\begin{aligned}
\check{Z}_i(t)&=\widetilde{\mathcal{M}}_{i0}\check{X}(t)+\overline{\mathcal{M}}_{i0}\check{\hat{X}}(t)
              +\widetilde{\mathcal{M}}_{i1}\check{Z}_1(t)+\overline{\mathcal{M}}_{i1}\check{\hat{Z}}_1(t)\\
             &\quad+\widetilde{\mathcal{M}}_{i2}\check{Z}_2(t)+\overline{\mathcal{M}}_{i2}\check{\hat{Z}}_2(t)
              +\widetilde{\mathcal{M}}_{i3}\check{Z}_3(t)+\overline{\mathcal{M}}_{i3}\check{\hat{Z}}_3(t),\ i=1,2,3,
\end{aligned}
\end{equation}
where
{\footnotesize\begin{equation*}
\left\{
\begin{aligned}
\widetilde{\mathcal{M}}_{10}&\triangleq\mathcal{P}_1\mathcal{A}_1+\mathcal{P}_1(\mathcal{B}_1+\widetilde{\mathcal{B}}_1)(\mathcal{P}_1+\mathcal{P}_3),\\
 \overline{\mathcal{M}}_{10}&\triangleq\mathcal{P}_1(\widehat{\mathcal{A}}_1+\overline{\mathcal{A}}_1)
                             +\mathcal{P}_2(\mathcal{A}_1+\widehat{\mathcal{A}}_1+\widetilde{\mathcal{A}}_1+\overline{\mathcal{A}}_1)
                             +\mathcal{P}_1(\mathcal{B}_1+\widetilde{\mathcal{B}}_1)\big(\mathcal{P}_2+\mathcal{P}_4\big)\\
                            &\quad+\big(\mathcal{P}_1\overline{\mathcal{B}}_1+\mathcal{P}_2\mathcal{B}_1+\mathcal{P}_2\widetilde{\mathcal{B}}_1
                             +\mathcal{P}_2\overline{\mathcal{B}}_1\big)\big(\mathcal{P}_1+\mathcal{P}_2+\mathcal{P}_3+\mathcal{P}_4\big),\\
\widetilde{\mathcal{M}}_{20}&\triangleq\mathcal{P}_1\mathcal{A}_2+\mathcal{P}_3(\mathcal{A}_2+\widetilde{\mathcal{A}}_2)
                             +(\mathcal{P}_1+\mathcal{P}_3)(\mathcal{B}_2+\widetilde{\mathcal{B}}_2)(\mathcal{P}_1+\mathcal{P}_3),\\
 \overline{\mathcal{M}}_{20}&\triangleq(\mathcal{P}_1+\mathcal{P}_3)(\widehat{\mathcal{A}}_2+\overline{\mathcal{A}}_2)
                             +(\mathcal{P}_1+\mathcal{P}_3)(\mathcal{B}_2+\widetilde{\mathcal{B}}_2)\big(\mathcal{P}_2+\mathcal{P}_4\big)\\
                            &\quad+(\mathcal{P}_1+\mathcal{P}_3)\overline{\mathcal{B}}_2\big(\mathcal{P}_1+\mathcal{P}_2+\mathcal{P}_3+\mathcal{P}_4\big),\\
\widetilde{\mathcal{M}}_{30}&\triangleq\mathcal{P}_1\mathcal{A}_3+\mathcal{P}_3(\mathcal{A}_3+\widetilde{\mathcal{A}}_3)
                             +(\mathcal{P}_1+\mathcal{P}_3)(\mathcal{B}_3+\widetilde{\mathcal{B}}_3)(\mathcal{P}_1+\mathcal{P}_3),\\
 \overline{\mathcal{M}}_{30}&\triangleq(\mathcal{P}_1+\mathcal{P}_3)(\widehat{\mathcal{A}}_3+\overline{\mathcal{A}}_3)
                             +(\mathcal{P}_2+\mathcal{P}_4)(\mathcal{A}_3+\widehat{\mathcal{A}}_3+\widetilde{\mathcal{A}}_3+\overline{\mathcal{A}}_3)
                             +(\mathcal{P}_1+\mathcal{P}_3)(\mathcal{B}_3+\widetilde{\mathcal{B}}_3)\big(\mathcal{P}_2+\mathcal{P}_4\big)\\
                            &\quad+\big[(\mathcal{P}_2+\mathcal{P}_2)(\mathcal{B}_3+\widetilde{\mathcal{B}}_3+\overline{\mathcal{B}}_3)
                             +\mathcal{P}_1\overline{\mathcal{B}}_3+\mathcal{P}_3\widetilde{\mathcal{B}}_3\big]\big(\mathcal{P}_1+\mathcal{P}_2+\mathcal{P}_3+\mathcal{P}_4\big),\\
\widetilde{\mathcal{M}}_{11}&\triangleq\mathcal{P}_1\mathcal{B}_1+\mathcal{P}_1\widetilde{\mathcal{C}}_1,\quad
 \overline{\mathcal{M}}_{11}\triangleq\mathcal{P}_2\mathcal{C}_1+\mathcal{P}_1\overline{\mathcal{C}}_1
                               +\mathcal{P}_2\widetilde{\mathcal{C}}_1+\mathcal{P}_2\overline{\mathcal{C}}_1,\\
\widetilde{\mathcal{M}}_{12}&\triangleq\mathcal{P}_1\mathcal{D}_1+\mathcal{P}_1\widetilde{\mathcal{D}}_1,\quad
 \overline{\mathcal{M}}_{12}\triangleq\mathcal{P}_2\mathcal{D}_1+\mathcal{P}_1\overline{\mathcal{D}}_1
                               +\mathcal{P}_2\widetilde{\mathcal{D}}_1+\mathcal{P}_2\overline{\mathcal{D}}_1,\\
\widetilde{\mathcal{M}}_{13}&\triangleq\mathcal{P}_1\mathcal{E}_1+\mathcal{P}_1\widetilde{\mathcal{E}}_1,\quad
 \overline{\mathcal{M}}_{13}\triangleq\mathcal{P}_2\mathcal{E}_1+\mathcal{P}_1\overline{\mathcal{E}}_1
                               +\mathcal{P}_2\widetilde{\mathcal{E}}_1+\mathcal{P}_2\overline{\mathcal{E}}_1,\\
\widetilde{\mathcal{M}}_{21}&\triangleq\mathcal{P}_1\mathcal{B}_2+\mathcal{P}_1\widetilde{\mathcal{C}}_2+\mathcal{P}_3\mathcal{C}_2+\mathcal{P}_3\widetilde{\mathcal{C}}_2,\quad
 \overline{\mathcal{M}}_{21}\triangleq\mathcal{P}_1\overline{\mathcal{C}}_2+\mathcal{P}_3\overline{\mathcal{C}}_2,\\
\widetilde{\mathcal{M}}_{22}&\triangleq\mathcal{P}_1\mathcal{D}_2+\mathcal{P}_1\widetilde{\mathcal{D}}_2+\mathcal{P}_3\mathcal{D}_2+\mathcal{P}_3\widetilde{\mathcal{D}}_2,\quad
 \overline{\mathcal{M}}_{22}\triangleq\mathcal{P}_1\overline{\mathcal{D}}_2+\mathcal{P}_3\overline{\mathcal{D}}_2,\\
\widetilde{\mathcal{M}}_{23}&\triangleq\mathcal{P}_1\mathcal{E}_2+\mathcal{P}_1\widetilde{\mathcal{E}}_2+\mathcal{P}_3\mathcal{E}_2+\mathcal{P}_3\widetilde{\mathcal{E}}_2,\quad
 \overline{\mathcal{M}}_{23}\triangleq\mathcal{P}_1\overline{\mathcal{E}}_2+\mathcal{P}_3\overline{\mathcal{E}}_2,\\
\widetilde{\mathcal{M}}_{31}&\triangleq\mathcal{P}_1\mathcal{B}_3+\mathcal{P}_1\widetilde{\mathcal{C}}_3+\mathcal{P}_3\widetilde{\mathcal{C}}_3
                             +\mathcal{P}_3\widetilde{\mathcal{C}}_3,\quad
 \overline{\mathcal{M}}_{31}\triangleq(\mathcal{P}_1+\mathcal{P}_3)\overline{\mathcal{C}}_3
                             +(\mathcal{P}_2+\mathcal{P}_4)(\mathcal{C}_3+\widetilde{\mathcal{C}}_3+\overline{\mathcal{C}}_3),\\
\widetilde{\mathcal{M}}_{32}&\triangleq\mathcal{P}_1\mathcal{D}_3+\mathcal{P}_1\widetilde{\mathcal{D}}_3+\mathcal{P}_3\widetilde{\mathcal{D}}_3
                             +\mathcal{P}_3\widetilde{\mathcal{D}}_3,\quad
 \overline{\mathcal{M}}_{32}\triangleq(\mathcal{P}_1+\mathcal{P}_3)\overline{\mathcal{D}}_3
                             +(\mathcal{P}_2+\mathcal{P}_4)(\mathcal{D}_3+\widetilde{\mathcal{D}}_3+\overline{\mathcal{D}}_3),\\
\widetilde{\mathcal{M}}_{33}&\triangleq\mathcal{P}_1\mathcal{E}_3+\mathcal{P}_1\widetilde{\mathcal{E}}_3+\mathcal{P}_3\widetilde{\mathcal{E}}_3
                             +\mathcal{P}_3\widetilde{\mathcal{E}}_3,\quad
 \overline{\mathcal{M}}_{33}\triangleq(\mathcal{P}_1+\mathcal{P}_3)\overline{\mathcal{E}}_3
                             +(\mathcal{P}_2+\mathcal{P}_4)(\mathcal{E}_3+\widetilde{\mathcal{E}}_3+\overline{\mathcal{E}}_3).
 \end{aligned}
\right.
\end{equation*}}
Putting (\ref{double filter estimate Zi-G1G2}) into (\ref{filter estimate Zi-G2}), we get
\begin{equation}\label{filter estimate Zi-G2G2}
\begin{aligned}
\check{Z}^i(t)&=\widetilde{\mathcal{M}}_{i0}\check{X}(t)+\Big[\overline{\mathcal{M}}_{i0}+\overline{\mathcal{M}}_{i1}\mathcal{N}_1
              +\overline{\mathcal{M}}_{i2}\mathcal{N}_2+\overline{\mathcal{M}}_{i3}\mathcal{N}_3\Big]\check{\hat{X}}(t)\\
             &\quad+\widetilde{\mathcal{M}}_{i1}\check{Z}_1(t)+\widetilde{\mathcal{M}}_{i2}\check{Z}_2(t)+\widetilde{\mathcal{M}}_{i3}\check{Z}_3(t)\\
             &\triangleq\widetilde{\mathcal{M}}^{i0}(t)\check{X}(t)+\overline{\mathcal{N}}^{i0}(t)\check{\hat{X}}(t)
              +\widetilde{\mathcal{M}}_{i1}\check{Z}_1(t)+\widetilde{\mathcal{M}}_{i2}\check{Z}_2(t)+\widetilde{\mathcal{M}}_{i3}\check{Z}_3(t),\ i=1,2,3.
\end{aligned}
\end{equation}
We rewrite (\ref{filter estimate Zi-G2G2}) as
\begin{equation}\label{linear equation system-2}
\begin{aligned}
\left(
\begin{array}{ccc}
I_n-\widetilde{\mathcal{M}}_{11}&-\widetilde{\mathcal{M}}_{12}&-\widetilde{\mathcal{M}}_{13}\\
-\widetilde{\mathcal{M}}_{21}&I_n-\widetilde{\mathcal{M}}_{22}&-\widetilde{\mathcal{M}}_{23}\\
-\widetilde{\mathcal{M}}_{31}&-\widetilde{\mathcal{M}}_{32}&I_n-\widetilde{\mathcal{M}}_{33}\\
\end{array}
\right)
\left(
\begin{array}{ccc}
\check{Z}_1(t)\\\check{Z}_2(t)\\\check{Z}_3(t)
\end{array}
\right)=
\left(
\begin{array}{ccc}
\widetilde{\mathcal{M}}_{10}\check{X}(t)+\overline{\mathcal{N}}_{10}\check{\hat{X}}(t)\\\widetilde{\mathcal{M}}_{20}\check{X}(t)+\overline{\mathcal{N}}_{20}\check{\hat{X}}(t)\\
\widetilde{\mathcal{M}}_{30}\check{X}(t)+\overline{\mathcal{N}}_{30}\check{\hat{X}}(t)
\end{array}
\right).
\end{aligned}
\end{equation}
Similarly, if we assume that

{\bf (A2.3)}\quad the coefficient matrix of (\ref{linear equation system-2}) is invertible, for any $t\in[0,T]$,

\noindent then we have
\begin{equation}\label{filter estimate Zi-G2G2G2}
\begin{aligned}
\check{Z}_i(t)&=(-1)^{i-1}\big(\mathbb{N}_2\big)^{-1}\Big[\big(\widetilde{\mathcal{M}}_{10}\check{X}(t)
              +\overline{\mathcal{N}}_{10}\check{\hat{X}}(t)\big)\widetilde{\mathbf{M}}_{1i}
              -\big(\widetilde{\mathcal{M}}_{20}\check{X}(t)+\overline{\mathcal{N}}_{20}\check{\hat{X}}(t)\big)\widetilde{\mathbf{M}}_{2i}\\
             &\qquad\qquad\qquad+\big(\widetilde{\mathcal{M}}_{30}\check{X}(t)+\overline{\mathcal{N}}_{30}(t)\check{\hat{X}}(t)\big)\widetilde{\mathbf{M}}_{3i}\Big]\\
             &=(-1)^{i-1}\big(\mathbb{N}_2\big)^{-1}\Big[\widetilde{\mathcal{M}}_{10}\widetilde{\mathbf{M}}_{1i}
              -\widetilde{\mathcal{M}}_{20}\widetilde{\mathbf{M}}_{2i}+\widetilde{\mathcal{M}}_{30}\widetilde{\mathbf{M}}_{3i}\Big]\check{X}(t)\\
             &\qquad\qquad+\big(\mathbb{N}_2\big)^{-1}\Big[\overline{\mathcal{N}}_{10}\widetilde{\mathbf{M}}_{1i}
              -\overline{\mathcal{N}}_{20}\widetilde{\mathbf{M}}_{2i}+\overline{\mathcal{N}}_{30}\widetilde{\mathbf{M}}_{3i}\Big]\check{X}(t)\\
             &\triangleq\widetilde{\mathcal{N}}_i\check{X}(t)+\overline{\mathcal{N}}_i\check{\hat{X}}(t),\ i=1,2,3,
\end{aligned}
\end{equation}
where $\mathbb{N}_2$ is the determinant of the coefficient of (\ref{linear equation system-2}), and $\widetilde{\mathbf{M}}_{ji}$ is the adjoint matrix of the $(j,i)$ element in (\ref{linear equation system-2}), for $j,i=1,2,3$.

\vspace{1mm}

{\bf Step 3.}\quad Taking $\mathbb{E}[\cdot|\mathcal{G}^1_t]$ on both sides of (\ref{comparing dW1,dW2,dW3-leader}), we obtain
\begin{equation}\label{filter estimate Zi-G1}
\begin{aligned}
\hat{Z}_i(t)&=\widehat{\mathcal{M}}_{i0}\hat{X}(t)+\overline{\overline{\mathcal{M}}}_{i0}\check{\hat{X}}(t)
            +\widehat{\mathcal{M}}_{i1}\hat{Z}_1(t)+\overline{\overline{\mathcal{M}}}_{i1}\check{\hat{Z}}_1(t)\\
           &\quad+\widehat{\mathcal{M}}_{i2}\hat{Z}_2(t)+\overline{\overline{\mathcal{M}}}_{i2}\check{\hat{Z}}_2(t)
            +\widehat{\mathcal{M}}_{i3}\hat{Z}_3(t)+\overline{\overline{\mathcal{M}}}_{i3}\check{\hat{Z}}_3(t),\ i=1,2,3,
\end{aligned}
\end{equation}
where
{\footnotesize\begin{equation*}
\left\{
\begin{aligned}
            \widehat{\mathcal{M}}_{10}&\triangleq(\mathcal{P}_1+\mathcal{P}_2)(\mathcal{A}_1+\widehat{\mathcal{A}}_1)
                                       +(\mathcal{P}_1+\mathcal{P}_2)\mathcal{B}_1(\mathcal{P}_1+\mathcal{P}_2),\\
\overline{\overline{\mathcal{M}}}_{10}&\triangleq(\mathcal{P}_1+\mathcal{P}_2)\overline{\mathcal{A}}_1+\mathcal{P}_2\widetilde{\mathcal{A}}_1
                                       +(\mathcal{P}_1+\mathcal{P}_2)\mathcal{B}_1(\mathcal{P}_3+\mathcal{P}_4)
                                       +(\mathcal{P}_1+\mathcal{P}_2)(\widetilde{\mathcal{B}}_1+\overline{\mathcal{B}}_1)(\mathcal{P}_1+\mathcal{P}_2+\mathcal{P}_3+\mathcal{P}_4),\\
            \widehat{\mathcal{M}}_{20}&\triangleq\mathcal{P}_1(\mathcal{A}_2+\widehat{\mathcal{A}}_2)+\mathcal{P}_1\mathcal{B}_2(\mathcal{P}_1+\mathcal{P}_2),\\
\overline{\overline{\mathcal{M}}}_{20}&\triangleq\mathcal{P}_1\overline{\mathcal{A}}_2+\mathcal{P}_3(\mathcal{A}_2+\widehat{\mathcal{A}}_2
                                       +\widetilde{\mathcal{A}}_2+\overline{\mathcal{A}}_2)+\mathcal{P}_1\mathcal{B}_2(\mathcal{P}_3+\mathcal{P}_4)\\
                                      &\quad+(\mathcal{P}_1\widetilde{\mathcal{B}}_2+\mathcal{P}_1\overline{\mathcal{B}}_2
                                       +\mathcal{P}_3\mathcal{B}_2+\mathcal{P}_3\widetilde{\mathcal{B}}_2+\mathcal{P}_3\overline{\mathcal{B}}_3)
                                       (\mathcal{P}_1+\mathcal{P}_2+\mathcal{P}_3+\mathcal{P}_4),\\
            \widehat{\mathcal{M}}_{30}&\triangleq(\mathcal{P}_1+\mathcal{P}_2)(\mathcal{A}_3+\widehat{\mathcal{A}}_3)
                                       +(\mathcal{P}_1+\mathcal{P}_2)\mathcal{B}_3(\mathcal{P}_1+\mathcal{P}_2),\\
\overline{\overline{\mathcal{M}}}_{30}&\triangleq\mathcal{P}_1\overline{\mathcal{A}}_3+\mathcal{P}_2(\widetilde{\mathcal{A}}_3+\widehat{\mathcal{A}}_3)
                                       +(\mathcal{P}_3+\mathcal{P}_4)(\mathcal{A}_3+\widehat{\mathcal{A}}_3
                                       +\widetilde{\mathcal{A}}_3+\overline{\mathcal{A}}_3)+(\mathcal{P}_1+\mathcal{P}_2)\mathcal{B}_3(\mathcal{P}_3+\mathcal{P}_4)\\
                                      &\quad+\big[(\mathcal{P}_1+\mathcal{P}_2)(\widetilde{\mathcal{B}}_3+\overline{\mathcal{B}}_3)
                                       +(\mathcal{P}_3+\mathcal{P}_4)(\mathcal{B}_3+\widetilde{\mathcal{B}}_2+\overline{\mathcal{B}}_3)\big]
                                       (\mathcal{P}_1+\mathcal{P}_2+\mathcal{P}_3+\mathcal{P}_4),\\
\end{aligned}
\right.
\end{equation*}}
and
{\footnotesize\begin{equation*}
\hspace{-2cm}\left\{
\begin{aligned}
            \widehat{\mathcal{M}}_{11}&\triangleq\mathcal{P}_1\mathcal{B}_1+\mathcal{P}_2\mathcal{C}_1,\quad
 \overline{\overline{\mathcal{M}}}_{11}\triangleq\mathcal{P}_1\widetilde{\mathcal{C}}_1+\mathcal{P}_1\overline{\mathcal{C}}_1
                                       +\mathcal{P}_2\widetilde{\mathcal{C}}_1+\mathcal{P}_2\overline{\mathcal{C}}_1,\\
            \widehat{\mathcal{M}}_{12}&\triangleq\mathcal{P}_1\mathcal{D}_1+\mathcal{P}_2\mathcal{D}_1,\quad
 \overline{\overline{\mathcal{M}}}_{12}\triangleq\mathcal{P}_1\widetilde{\mathcal{D}}_1+\mathcal{P}_1\overline{\mathcal{D}}_1
                                       +\mathcal{P}_2\widetilde{\mathcal{D}}_1+\mathcal{P}_2\overline{\mathcal{D}}_1,\\
            \widehat{\mathcal{M}}_{13}&\triangleq\mathcal{P}_1\mathcal{E}_1+\mathcal{P}_2\mathcal{E}_1,\quad
 \overline{\overline{\mathcal{M}}}_{13}\triangleq\mathcal{P}_1\widetilde{\mathcal{E}}_1+\mathcal{P}_1\overline{\mathcal{E}}_1
                                       +\mathcal{P}_2\widetilde{\mathcal{E}}_1+\mathcal{P}_2\overline{\mathcal{E}}_1,\\
            \widehat{\mathcal{M}}_{21}&\triangleq\mathcal{P}_1\mathcal{B}_2,\quad
 \overline{\overline{\mathcal{M}}}_{21}\triangleq\mathcal{P}_1\widetilde{\mathcal{C}}_2+\mathcal{P}_1\overline{\mathcal{C}}_2+\mathcal{P}_3\mathcal{C}_2
                                       +\mathcal{P}_3\widetilde{\mathcal{C}}_2+\mathcal{P}_3\overline{\mathcal{C}}_2,\\
            \widehat{\mathcal{M}}_{22}&\triangleq\mathcal{P}_1\mathcal{D}_2,\quad
 \overline{\overline{\mathcal{M}}}_{22}\triangleq\mathcal{P}_1\widetilde{\mathcal{D}}_2+\mathcal{P}_1\overline{\mathcal{D}}_2+\mathcal{P}_3\mathcal{D}_2
                                       +\mathcal{P}_3\widetilde{\mathcal{D}}_2+\mathcal{P}_3\overline{\mathcal{D}}_2,\\
            \widehat{\mathcal{M}}_{23}&\triangleq\mathcal{P}_1\mathcal{E}_2,\quad
 \overline{\overline{\mathcal{M}}}_{23}\triangleq\mathcal{P}_1\widetilde{\mathcal{E}}_2+\mathcal{P}_1\overline{\mathcal{E}}_2+\mathcal{P}_3\mathcal{E}_2
                                       +\mathcal{P}_3\widetilde{\mathcal{E}}_2+\mathcal{P}_3\overline{\mathcal{E}}_2,\\
            \widehat{\mathcal{M}}_{31}&\triangleq\mathcal{P}_1\mathcal{B}_3+\mathcal{P}_2\mathcal{C}_3,\quad
 \overline{\overline{\mathcal{M}}}_{31}\triangleq(\mathcal{P}_1+\mathcal{P}_2)(\widetilde{\mathcal{C}}_3+\overline{\mathcal{C}}_3)
                                       +(\mathcal{P}_3+\mathcal{P}_4)(\mathcal{C}_3+\widetilde{\mathcal{C}}_3+\overline{\mathcal{C}}_3),\\
            \widehat{\mathcal{M}}_{32}&\triangleq\mathcal{P}_1\mathcal{D}_3+\mathcal{P}_2\mathcal{D}_3,\quad
 \overline{\overline{\mathcal{M}}}_{31}\triangleq(\mathcal{P}_1+\mathcal{P}_2)(\widetilde{\mathcal{D}}_3+\overline{\mathcal{D}}_3)
                                       +(\mathcal{P}_3+\mathcal{P}_4)(\mathcal{D}_3+\widetilde{\mathcal{D}}_3+\overline{\mathcal{D}}_3),\\
            \widehat{\mathcal{M}}_{33}&\triangleq\mathcal{P}_1\mathcal{E}_3+\mathcal{P}_2\mathcal{E}_3,\quad
 \overline{\overline{\mathcal{M}}}_{31}\triangleq(\mathcal{P}_1+\mathcal{P}_2)(\widetilde{\mathcal{E}}_3+\overline{\mathcal{E}}_3)
                                       +(\mathcal{P}_3+\mathcal{P}_4)(\mathcal{E}_3+\widetilde{\mathcal{E}}_3+\overline{\mathcal{E}}_3).
\end{aligned}
\right.
\end{equation*}}
Putting (\ref{double filter estimate Zi-G1G2}) into (\ref{filter estimate Zi-G1}), we get
\begin{equation}\label{filter estimate Zi-G1G1}
\begin{aligned}
\hat{Z}_i(t)&=\widehat{\mathcal{M}}_{i0}\hat{X}(t)+\Big[\overline{\overline{\mathcal{M}}}_{i0}+\overline{\overline{\mathcal{M}}}_{i1}\mathcal{N}_1
              +\overline{\overline{\mathcal{M}}}_{i2}\mathcal{N}_2+\overline{\overline{\mathcal{M}}}_{i3}\mathcal{N}_3\Big]\check{\hat{X}}(t)\\
             &\quad+\widehat{\mathcal{M}}_{i1}\hat{Z}_1(t)+\widehat{\mathcal{M}}_{i2}\hat{Z}_2(t)+\widehat{\mathcal{M}}_{i3}\hat{Z}_3(t)\\
             &\triangleq\widehat{\mathcal{M}}_{i0}\hat{X}(t)+\overline{\overline{\mathcal{N}}}_{i0}\check{\hat{X}}(t)
              +\widehat{\mathcal{M}}_{i1}\hat{Z}_1(t)+\widehat{\mathcal{M}}_{i2}\hat{Z}_2(t)+\widehat{\mathcal{M}}_{i3}\hat{Z}_3(t),\ i=1,2,3.
\end{aligned}
\end{equation}
We rewrite (\ref{filter estimate Zi-G1G1}) as
\begin{equation}\label{linear equation system-3}
\begin{aligned}
\left(
\begin{array}{ccc}
I_n-\widehat{\mathcal{M}}_{11}&-\widehat{\mathcal{M}}_{12}&-\widehat{\mathcal{M}}_{13}\\
-\widehat{\mathcal{M}}_{21}&I_n-\widehat{\mathcal{M}}_{22}&-\widehat{\mathcal{M}}_{23}\\
-\widehat{\mathcal{M}}_{31}&-\widehat{\mathcal{M}}_{32}&I_n-\widehat{\mathcal{M}}_{33}\\
\end{array}
\right)
\left(
\begin{array}{ccc}
\hat{Z}_1(t)\\\hat{Z}_2(t)\\\hat{Z}_3(t)
\end{array}
\right)=
\left(
\begin{array}{ccc}
\widehat{\mathcal{M}}_{10}\hat{X}(t)+\overline{\overline{\mathcal{N}}}_{10}\check{\hat{X}}(t)\\
\widehat{\mathcal{M}}_{20}\hat{X}(t)+\overline{\overline{\mathcal{N}}}_{20}\check{\hat{X}}(t)\\
\widehat{\mathcal{M}}_{30}\hat{X}(t)+\overline{\overline{\mathcal{N}}}_{30}\check{\hat{X}}(t)
\end{array}
\right).
\end{aligned}
\end{equation}
Similarly, if we assume that

{\bf (A2.4)}\quad the coefficient matrix of (\ref{linear equation system-3}) is invertible, for any $t\in[0,T]$,

\noindent then we have
\begin{equation}\label{filter estimate Zi-G1G1G1}
\begin{aligned}
 \hat{Z}_i(t)&=(-1)^{i-1}\big(\mathbb{N}_3\big)^{-1}\Big[\big(\widehat{\mathcal{M}}_{10}\hat{X}(t)
              +\overline{\overline{\mathcal{N}}}_{10}\check{\hat{X}}(t)\big)\widehat{\mathbf{M}}_{1i}
              -\big(\widehat{\mathcal{M}}_{20}\hat{X}(t)+\overline{\overline{\mathcal{N}}}_{20}\check{\hat{X}}(t)\big)\widehat{\mathbf{M}}_{2i}\\
             &\qquad\qquad\qquad\qquad+\big(\widehat{\mathcal{M}}_{30}\hat{X}(t)+\overline{\overline{\mathcal{N}}}_{30}\check{\hat{X}}(t)\big)\widehat{\mathbf{M}}_{3i}\Big]\\
             &=(-1)^{i-1}\big(\mathbb{N}_3\big)^{-1}\Big[\widehat{\mathcal{M}}_{10}\widehat{\mathbf{M}}_{1i}
              -\widehat{\mathcal{M}}_{20}\widehat{\mathbf{M}}_{2i}+\widehat{\mathcal{M}}_{30}\widehat{\mathbf{M}}_{3i}\Big]\hat{X}(t)\\
             &\qquad\qquad+\big(\mathbb{N}_3\big)^{-1}\Big[\overline{\overline{\mathcal{N}}}_{10}\widehat{\mathbf{M}}_{1i}
              -\overline{\overline{\mathcal{N}}}_{20}\widehat{\mathbf{M}}_{2i}+\overline{\overline{\mathcal{N}}}_{30}\widehat{\mathbf{M}}_{3i}\Big]\hat{X}(t)\\
             &\triangleq\widehat{\mathcal{N}}_i\hat{X}(t)+\overline{\overline{\mathcal{N}}}_i\check{\hat{X}}(t),\ i=1,2,3,
\end{aligned}
\end{equation}
where $\mathbb{N}_3$ is the determinant of the coefficient of (\ref{linear equation system-3}), and $\widehat{\mathbf{M}}_{ji}$ is the adjoint matrix of the $(j,i)$ element in (\ref{linear equation system-3}), for $j,i=1,2,3$.

\vspace{1mm}

{\bf Step 4.}\quad Putting (\ref{double filter estimate Zi-G1G2}), (\ref{filter estimate Zi-G2G2G2}) and (\ref{filter estimate Zi-G1G1G1}) into (\ref{comparing dW1,dW2,dW3-leader}), we have
\begin{equation}\label{Zi}
\begin{aligned}
Z_i(t)&=\Gamma_{i0}X(t)+\widehat{\Gamma}_{i0}\hat{X}(t)+\widetilde{\Gamma}_{i0}\check{X}(t)+\overline{\Gamma}_{i0}\check{\hat{X}}(t)\\
      &\quad+\mathcal{P}_1\mathcal{B}_iZ_1(t)+\mathcal{P}_1\mathcal{D}_iZ_2(t)+\mathcal{P}_1\mathcal{E}_iZ_3(t),\ i=1,2,3,
\end{aligned}
\end{equation}
where
{\footnotesize\begin{equation*}
\left\{
\begin{aligned}
            \Gamma_{10}&\triangleq\mathcal{P}_1\mathcal{A}_1+\mathcal{P}_1\mathcal{B}_1\mathcal{P}_1,\quad
                   \Gamma_{20}\triangleq\mathcal{P}_1\mathcal{A}_2+\mathcal{P}_1\mathcal{B}_2\mathcal{P}_1,\quad \Gamma_{30}\triangleq\mathcal{P}_1\mathcal{A}_3+\mathcal{P}_1\mathcal{B}_3\mathcal{P}_1,\\
  \widehat{\Gamma}_{10}&\triangleq\mathcal{P}_1\widehat{\mathcal{A}}_1+\mathcal{P}_2\big(\mathcal{A}_1+\widehat{\mathcal{A}}_1\big)
                        +\mathcal{P}_1\mathcal{B}_1\mathcal{P}_2+\mathcal{P}_2\mathcal{B}_1\big(\mathcal{P}_1+\mathcal{P}_2\big)
                        +\mathcal{P}_2\mathcal{C}_1\widehat{\mathcal{N}}_1+\mathcal{P}_2\mathcal{D}_1\widehat{\mathcal{N}}_2+\mathcal{P}_2\mathcal{E}_1\widehat{\mathcal{N}}_3,\\
\widetilde{\Gamma}_{10}&\triangleq\mathcal{P}_1\mathcal{B}_1\mathcal{P}_3+\mathcal{P}_1\widetilde{\mathcal{B}}_1\big(\mathcal{P}_1+\mathcal{P}_3\big)
                        +\mathcal{P}_1\widetilde{\mathcal{C}}_1\widetilde{\mathcal{N}}_1+\mathcal{P}_1\widetilde{\mathcal{D}}_1\widetilde{\mathcal{N}}_2
                        +\mathcal{P}_1\mathcal{E}_1\widetilde{\mathcal{N}}_3,\\
 \overline{\Gamma}_{10}&\triangleq\mathcal{P}_1\overline{\mathcal{A}}_1+\mathcal{P}_1\mathcal{B}_1\mathcal{P}_4
                        +\mathcal{P}_1\widetilde{\mathcal{B}}_1\big(\mathcal{P}_2+\mathcal{P}_4\big)
                        +\mathcal{P}_1\overline{\mathcal{B}}_1\big(\mathcal{P}_1+\mathcal{P}_2+\mathcal{P}_3+\mathcal{P}_4\big)
                        +\mathcal{P}_2\big(\widetilde{\mathcal{A}}_1+\overline{\mathcal{A}}_1\big)\\
                       &\quad+\mathcal{P}_2\mathcal{B}_1\big(\mathcal{P}_3+\mathcal{P}_4\big)
                        +\mathcal{P}_2\big(\widetilde{\mathcal{B}}_1+\overline{\mathcal{B}}_1\big)\big(\mathcal{P}_1+\mathcal{P}_2+\mathcal{P}_3+\mathcal{P}_4\big)
                        +\mathcal{P}_2\mathcal{C}_1\overline{\overline{\mathcal{N}}}_1+\mathcal{P}_1\widetilde{\mathcal{C}}_1\overline{\mathcal{N}}_1\\
                       &\quad
                        +\mathcal{P}_2\mathcal{D}_1\overline{\overline{\mathcal{N}}}_2+\mathcal{P}_1\widetilde{\mathcal{D}}_1\overline{\mathcal{N}}_2
                        +\mathcal{P}_2\mathcal{E}_1\overline{\overline{\mathcal{N}}}_3+\mathcal{P}_1\widetilde{\mathcal{E}}_1\overline{\mathcal{N}}_3
                        +\big(\mathcal{P}_1\overline{\mathcal{C}}_1+\mathcal{P}_2\widetilde{\mathcal{C}}_1+\mathcal{P}_2\overline{\mathcal{C}}_1\big)\mathcal{N}_1,\\
   \widehat{\Gamma}_{20}&\triangleq\mathcal{P}_1\widehat{\mathcal{A}}_2+\mathcal{P}_1\mathcal{B}_2\mathcal{P}_2,\\
\widetilde{\Gamma}_{20}&\triangleq\mathcal{P}_1\mathcal{B}_2\mathcal{P}_3+\mathcal{P}_1\widetilde{\mathcal{B}}_2\big(\mathcal{P}_1+\mathcal{P}_3\big)
                        +\mathcal{P}_3\big(\mathcal{A}_2+\widetilde{\mathcal{A}}_2\big)
                        +\mathcal{P}_3\big(\mathcal{B}_2+\widetilde{\mathcal{B}}_2\big)\big(\mathcal{P}_1+\mathcal{P}_3\big)\\
                       &\quad+\big(\mathcal{P}_1\widetilde{\mathcal{C}}_2+\mathcal{P}_3\mathcal{C}_2+\mathcal{P}_3\widetilde{\mathcal{C}}_2\big)\widetilde{\mathcal{N}}_1
                        +\big(\mathcal{P}_1\widetilde{\mathcal{D}}_2+\mathcal{P}_3\mathcal{D}_2+\mathcal{P}_3\widetilde{\mathcal{D}}_2\big)\widetilde{\mathcal{N}}_2
                        +\big(\mathcal{P}_1\widetilde{\mathcal{E}}_2+\mathcal{P}_3\mathcal{E}_2+\mathcal{P}_3\widetilde{\mathcal{E}}_2\big)\widetilde{\mathcal{N}}_3,\\
 \overline{\Gamma}_{20}&\triangleq\mathcal{P}_1\overline{\mathcal{A}}_2+\mathcal{P}_1\mathcal{B}_2\mathcal{P}_4
                        +\mathcal{P}_1\widetilde{\mathcal{B}}_2\big(\mathcal{P}_2+\mathcal{P}_4\big)+\mathcal{P}_3\big(\widehat{\mathcal{A}}_2+\overline{\mathcal{A}}_2\big)
                        +\mathcal{P}_1\overline{\mathcal{B}}_2\big(\mathcal{P}_1+\mathcal{P}_2+\mathcal{P}_3+\mathcal{P}_4\big)\\
                       &\quad+\mathcal{P}_3\big(\mathcal{B}_2+\widetilde{\mathcal{B}}_2\big)\big(\mathcal{P}_2+\mathcal{P}_4\big)
                        +\mathcal{P}_3\overline{\mathcal{B}}_2\big(\mathcal{P}_1+\mathcal{P}_2+\mathcal{P}_3+\mathcal{P}_4\big)\\
                       &\quad+\big(\mathcal{P}_1\widetilde{\mathcal{C}}_2+\mathcal{P}_3\mathcal{C}_2+\mathcal{P}_3\widetilde{\mathcal{C}}_2\big)\overline{\mathcal{N}}_1
                        +\big(\mathcal{P}_1\widetilde{\mathcal{D}}_2+\mathcal{P}_3\mathcal{D}_2+\mathcal{P}_3\widetilde{\mathcal{D}}_2\big)\overline{\mathcal{N}}_2
                        +\big(\mathcal{P}_1\widetilde{\mathcal{E}}_2+\mathcal{P}_3\mathcal{E}_2+\mathcal{P}_3\widetilde{\mathcal{E}}_2\big)\overline{\mathcal{N}}_3\\
                       &\quad+\big(\mathcal{P}_1\overline{\mathcal{C}}_2+\mathcal{P}_3\overline{\mathcal{C}}_2\big)\mathcal{N}_1
                        +\big(\mathcal{P}_1\overline{\mathcal{D}}_2+\mathcal{P}_3\overline{\mathcal{D}}_2\big)\mathcal{N}_2
                        +\big(\mathcal{P}_1\overline{\mathcal{E}}_2+\mathcal{P}_3\overline{\mathcal{E}}_2\big)\mathcal{N}_3,\\
  \widehat{\Gamma}_{30}&\triangleq\mathcal{P}_1\widehat{\mathcal{A}}_3+\mathcal{P}_1\mathcal{B}_3\mathcal{P}_2
                        +\mathcal{P}_2(\mathcal{A}_3+\widehat{\mathcal{A}}_3)+\mathcal{P}_2\mathcal{B}_3\big(\mathcal{P}_1+\mathcal{P}_2\big)
                        +\mathcal{P}_2\mathcal{C}_3\widehat{\mathcal{N}}_1+\mathcal{P}_2\mathcal{D}_3\widehat{\mathcal{N}}_2+\mathcal{P}_2\mathcal{E}_3\widehat{\mathcal{N}}_3,\\
\widetilde{\Gamma}_{30}&\triangleq\mathcal{P}_1\mathcal{B}_3\mathcal{P}_3+\mathcal{P}_1\widetilde{\mathcal{B}}_3\big(\mathcal{P}_1+\mathcal{P}_3\big)
                        +\mathcal{P}_3\big(\mathcal{A}_3+\widetilde{\mathcal{A}}_3\big)+\mathcal{P}_3\big(\mathcal{B}_3+\widetilde{\mathcal{B}}_3\big)
                        \big(\mathcal{P}_1+\mathcal{P}_3\big)\\
                       &\quad+\big(\mathcal{P}_1\widetilde{\mathcal{C}}_3+\mathcal{P}_3\mathcal{C}_3+\mathcal{P}_3\widetilde{\mathcal{C}}_3\big)\widetilde{\mathcal{N}}_1
                        +\big(\mathcal{P}_1\widetilde{\mathcal{D}}_3+\mathcal{P}_3\mathcal{D}_3+\mathcal{P}_3\widetilde{\mathcal{D}}_3\big)\widetilde{\mathcal{N}}_2
                        +\big(\mathcal{P}_1\widetilde{\mathcal{E}}_3+\mathcal{P}_3\mathcal{E}_3+\mathcal{P}_3\widetilde{\mathcal{E}}_3\big)\widetilde{\mathcal{N}}_3,\\
 \overline{\Gamma}_{30}&\triangleq\mathcal{P}_1\overline{\mathcal{A}}_3+\mathcal{P}_1\mathcal{B}_3\mathcal{P}_4
                        +\mathcal{P}_1\widetilde{\mathcal{B}}_3\big(\mathcal{P}_2+\mathcal{P}_4\big)
                        +\mathcal{P}_1\overline{\mathcal{B}}_3\big(\mathcal{P}_1+\mathcal{P}_2+\mathcal{P}_3+\mathcal{P}_4\big)
                        +\mathcal{P}_2\big(\widetilde{\mathcal{A}}_3+\overline{\mathcal{A}}_3\big)\\
                       &\quad+\mathcal{P}_2\mathcal{B}_3\big(\mathcal{P}_3+\mathcal{P}_4\big)+\mathcal{P}_2\big(\widetilde{\mathcal{B}}_3+\overline{\mathcal{B}}_3\big)
                        \big(\mathcal{P}_1+\mathcal{P}_2+\mathcal{P}_3+\mathcal{P}_4\big)+\mathcal{P}_3\big(\widehat{\mathcal{A}}_3+\overline{\mathcal{A}}_3\big)\\
                       &\quad+\mathcal{P}_3\big(\mathcal{B}_3+\widetilde{\mathcal{B}}_3\big)\big(\mathcal{P}_2+\mathcal{P}_4\big)
                        +\mathcal{P}_3\overline{\mathcal{B}}_3\big(\mathcal{P}_1+\mathcal{P}_2+\mathcal{P}_3+\mathcal{P}_4\big)
                        +\mathcal{P}_4\big(\mathcal{A}_3+\widetilde{\mathcal{A}}_3+\widehat{\mathcal{A}}_3+\overline{\mathcal{A}}_3\big)\\
                       &\quad+\mathcal{P}_4\big(\mathcal{B}_3+\widetilde{\mathcal{B}}_3+\overline{\mathcal{B}}_3\big)
                        \big(\mathcal{P}_1+\mathcal{P}_2+\mathcal{P}_3+\mathcal{P}_4\big)+\mathcal{P}_2\mathcal{C}_3\overline{\overline{\mathcal{N}}}_1
                        +\mathcal{P}_2\mathcal{D}_3\overline{\overline{\mathcal{N}}}_2+\mathcal{P}_2\mathcal{E}_3\overline{\overline{\mathcal{N}}}_3\\
                       &\quad+\big(\mathcal{P}_1\widetilde{\mathcal{C}}_3+\mathcal{P}_3\mathcal{C}_3+\mathcal{P}_3\widetilde{\mathcal{C}}_3\big)\overline{\mathcal{N}}_1
                        +\big(\mathcal{P}_1\widetilde{\mathcal{D}}_3+\mathcal{P}_3\mathcal{D}_3+\mathcal{P}_3\widetilde{\mathcal{D}}_3\big)\overline{\mathcal{N}}_2
                        +\big(\mathcal{P}_1\widetilde{\mathcal{E}}_3+\mathcal{P}_3\mathcal{E}_3+\mathcal{P}_3\widetilde{\mathcal{E}}_3\big)\overline{\mathcal{N}}_3\\
                       &\quad+\big[\mathcal{P}_1\overline{\mathcal{C}}_3+\mathcal{P}_2\widetilde{\mathcal{C}}_3+\mathcal{P}_2\overline{\mathcal{C}}_3
                        +\mathcal{P}_3\overline{\mathcal{C}}_3
                        +\mathcal{P}_4\big(\mathcal{C}_3+\widetilde{\mathcal{C}}_3+\overline{\mathcal{C}}_3\big)\big]\mathcal{N}_1\\
                       &\quad+\big[\mathcal{P}_1\overline{\mathcal{D}}_3+\mathcal{P}_2\widetilde{\mathcal{D}}_3+\mathcal{P}_2\overline{\mathcal{D}}_3
                        +\mathcal{P}_3\overline{\mathcal{D}}_3
                        +\mathcal{P}_4\big(\mathcal{D}_3+\widetilde{\mathcal{D}}_3+\overline{\mathcal{D}}_3\big)\big]\mathcal{N}_2\\
                       &\quad+\big[\mathcal{P}_1\overline{\mathcal{E}}_3+\mathcal{P}_2\widetilde{\mathcal{E}}_3+\mathcal{P}_2\overline{\mathcal{E}}_3
                        +\mathcal{P}_3\overline{\mathcal{E}}_3+\mathcal{P}_4\big(\mathcal{E}_3+\widetilde{\mathcal{E}}_3+\overline{\mathcal{E}}_3\big)\big]\mathcal{N}_3.
\end{aligned}
\right.
\end{equation*}}
We rewrite (\ref{Zi}) as
\begin{equation}\label{linear equation system-4}
\begin{aligned}
&\left(
\begin{array}{ccc}
I_n-\mathcal{P}_1\mathcal{B}_1&-\mathcal{P}_1\mathcal{B}_2&-\mathcal{P}_1\mathcal{B}_3\\
-\mathcal{P}_1\mathcal{D}_1&I_n-\mathcal{P}_1\mathcal{D}_2&-\mathcal{P}_1\mathcal{D}_3\\
-\mathcal{P}_1\mathcal{E}_1&-\mathcal{P}_1\mathcal{E}_2&I_n-\mathcal{P}_1\mathcal{E}_3\\
\end{array}
\right)
\left(
\begin{array}{ccc}
Z_1\\Z_2\\Z_3
\end{array}
\right)\\=
&\left(
\begin{array}{ccc}
\Gamma_{10}X(t)+\widehat{\Gamma}_{10}\hat{X}(t)+\widetilde{\Gamma}_{10}\check{X}(t)+\overline{\Gamma}_{i0}\check{\hat{X}}(t)\\
\Gamma_{20}X(t)+\widehat{\Gamma}_{20}\hat{X}(t)+\widetilde{\Gamma}_{20}\check{X}(t)+\overline{\Gamma}_{20}\check{\hat{X}}(t)\\
\Gamma_{30}X(t)+\widehat{\Gamma}_{30}\hat{X}(t)+\widetilde{\Gamma}_{30}\check{X}(t)+\overline{\Gamma}_{30}\check{\hat{X}}(t)
\end{array}
\right).
\end{aligned}
\end{equation}
Similarly, if we assume that

{\bf (A2.5)}\quad the coefficient matrix of (\ref{linear equation system-4}) is invertible, for any $t\in[0,T]$,

\noindent then we have
\begin{equation}\label{ZiZi}
\begin{aligned}
          &Z_i(t)=(-1)^{i-1}\big(\mathbb{N}_4\big)^{-1}\Big[\big(\Gamma_{10}X(t)+\widehat{\Gamma}_{10}\hat{X}(t)+\widetilde{\Gamma}_{10}\check{X}(t)
           +\overline{\Gamma}_{10}\check{\hat{X}}(t)\big)\overline{\mathbf{M}}_{1i}\\
          &\qquad\qquad-\big(\Gamma_{20}X(t)+\widehat{\Gamma}_{20}\hat{X}(t)+\widetilde{\Gamma}_{20}\check{X}(t)
           +\overline{\Gamma}_{20}\check{\hat{X}}(t)\big)\overline{\mathbf{M}}_{2i}\\
          &\qquad\qquad+\big(\Gamma_{30}X(t)+\widehat{\Gamma}_{30}\hat{X}(t)+\widetilde{\Gamma}_{30}\check{X}(t)
           +\overline{\Gamma}_{30}\check{\hat{X}}(t)\big)\overline{\mathbf{M}}_{3i}(t)\Big]\\
         =&\ (-1)^{i-1}\big(\mathbb{N}_4\big)^{-1}\Big\{\big[\Gamma_{10}\overline{\mathbf{M}}_{1i}-\Gamma_{20}\overline{\mathbf{M}}_{2i}+\Gamma_{30}\overline{\mathbf{M}}_{3i}\big]X(t)
           +\big[\widehat{\Gamma}_{10}\overline{\mathbf{M}}_{1i}-\widehat{\Gamma}_{20}\overline{\mathbf{M}}_{2i}+\widehat{\Gamma}_{30}\overline{\mathbf{M}}_{3i}\big]\hat{X}(t)\\
          &\qquad+\big[\widetilde{\Gamma}_{10}\overline{\mathbf{M}}_{1i}-\widetilde{\Gamma}_{20}\overline{\mathbf{M}}_{2i}+\widetilde{\Gamma}_{30}\overline{\mathbf{M}}_{3i}\big]\check{X}(t)
           +\big[\overline{\Gamma}_{10}\overline{\mathbf{M}}_{1i}-\overline{\Gamma}_{20}\overline{\mathbf{M}}_{2i}+\overline{\Gamma}_{30}\overline{\mathbf{M}}_{3i}\big]\check{\hat{X}}(t)\Big\}\\
\triangleq&\ \Sigma_iX(t)+\widehat{\Sigma}_i\hat{X}(t)+\widetilde{\Sigma}_i\check{X}(t)+\overline{\Sigma}_i\check{\hat{X}}(t),\ i=1,2,3,
\end{aligned}
\end{equation}
where $\mathbb{N}_4$ is the determinant of the coefficient of (\ref{linear equation system-3}), and $\overline{\mathbf{M}}_{ji}$ is the adjoint matrix of the $(j,i)$ element in (\ref{linear equation system-4}), for $j,i=1,2,3$.

After these four steps, we have obtained that
\begin{equation}\label{Zi and its filter estimates}
\begin{aligned}
  Z_i(t)&=\Sigma_iX(t)+\widehat{\Sigma}_i\hat{X}(t)+\widetilde{\Sigma}_i\check{X}(t)+\overline{\Sigma}_i\check{\hat{X}}(t),\ i=1,2,3.
\end{aligned}
\end{equation}

Now, comparing the $dt$ term in (\ref{Applying Ito's formula to Y}) and substituting (\ref{Zi and its filter estimates}) into it, we obtain
{\small\begin{equation}\label{system of Riccati equation}
\left\{
\begin{aligned}
&0=\dot{\mathcal{P}}_1+\mathcal{P}_1\mathcal{A}_0+\mathcal{A}_0\mathcal{P}_1+\mathcal{P}_1\mathcal{B}_0\mathcal{P}_1+\mathcal{P}_1\big(\mathcal{B}_1^\top\Sigma_1+\mathcal{B}_2^\top\Sigma_2+\mathcal{B}_3^\top\Sigma_3\big)
 +\mathcal{A}_1\Sigma_1+\mathcal{A}_2\Sigma_2+\mathcal{A}_3\Sigma_3+\mathcal{Q}_2,\\
&0=\dot{\mathcal{P}}_2+\mathcal{P}_1\widehat{\mathcal{A}}_0+\mathcal{A}_0\mathcal{P}_2+\mathcal{H}_1\mathcal{P}_1+\mathcal{H}_1\mathcal{P}_2+\mathcal{P}_2(\mathcal{A}_0+\widehat{\mathcal{A}}_0)
 +\mathcal{P}_1\mathcal{B}_0\mathcal{P}_2+\mathcal{P}_2\mathcal{B}_0\mathcal{P}_1+\mathcal{P}_2\mathcal{B}_0\mathcal{P}_2\\
&\quad\ +\mathcal{P}_1\big(\mathcal{B}_1^\top\widehat{\Sigma}_1+\mathcal{B}_2^\top\widehat{\Sigma}_2+\mathcal{B}_3\widehat{\Sigma}_3\big)
 +\mathcal{P}_2\big[\mathcal{B}_1^\top(\Sigma_1+\widehat{\Sigma}_1)+\mathcal{B}_2^\top(\Sigma_2+\widehat{\Sigma}_2)+\mathcal{B}_3^\top(\Sigma_3+\widehat{\Sigma}_3)\big]\\
&\quad\ +\mathcal{A}_1\widehat{\Sigma}_1+\mathcal{A}_2\widehat{\Sigma}_2+\mathcal{A}_3\widehat{\Sigma}_3+\widehat{\mathcal{A}}_2(\Sigma_2+\widehat{\Sigma}_2)+\widehat{\mathcal{A}}_3(\Sigma_3+\widehat{\Sigma}_3),\\
&0=\dot{\mathcal{P}}_3+\mathcal{P}_3\mathcal{A}_0+\mathcal{A}_0\mathcal{P}_3+\mathcal{P}_1\mathcal{B}_0\mathcal{P}_3
 +\mathcal{P}_1\mathcal{C}_0(\mathcal{P}_1+\mathcal{P}_2)+\mathcal{P}_3\big(\mathcal{B}_0+\mathcal{C}_0\big)(\mathcal{P}_1+\mathcal{P}_3)\\
&\quad\ +\mathcal{P}_1\big(\mathcal{B}_1^\top\widetilde{\Sigma}_1+\mathcal{B}_2^\top\widetilde{\Sigma}_2+\mathcal{B}_3^\top\widetilde{\Sigma}_3\big)
 +\big(\mathcal{P}_1\widetilde{\mathcal{B}}_1^\top+\mathcal{P}_3\mathcal{B}_1^\top+\mathcal{P}_3\widetilde{\mathcal{B}}_1^\top\big)(\Sigma_1+\widetilde{\Sigma}_1)
 +\mathcal{A}_1\widetilde{\Sigma}_1+\mathcal{A}_2\widetilde{\Sigma}_2\\
&\quad\ +\mathcal{A}_3\widetilde{\Sigma}_3+\big(\mathcal{P}_1\widetilde{\mathcal{B}}_2^\top+\mathcal{P}_3\mathcal{B}_2^\top+\mathcal{P}_3\widetilde{\mathcal{B}}_2^\top\big)(\Sigma_2+\widetilde{\Sigma}_2)
 +\big(\mathcal{P}_1\widetilde{\mathcal{B}}_3^\top+\mathcal{P}_3\mathcal{B}_3^\top+\mathcal{P}_3\widetilde{\mathcal{B}}_3^\top\big)(\Sigma_3+\widetilde{\Sigma}_3),\\
&0=\dot{\mathcal{P}}_4+\mathcal{P}_4(\mathcal{A}_0+\widehat{\mathcal{A}}_0+\overline{\mathcal{A}}_0)
 +\mathcal{A}_0\mathcal{P}_4+\mathcal{H}_2\mathcal{P}_3+\mathcal{H}_2\mathcal{P}_4+\mathcal{P}_3(\widehat{\mathcal{A}}_0+\overline{\mathcal{A}}_0)
 +\mathcal{P}_3(\mathcal{B}_0+\mathcal{C}_0)(\mathcal{P}_2+\mathcal{P}_4)\\
&\quad\ +\mathcal{P}_1\overline{\mathcal{A}}_0+\mathcal{H}_1+\big[\overline{\mathcal{A}}_0^\top+\mathcal{P}_2(\mathcal{C}_0+\widetilde{\mathcal{C}}_0)
 +\mathcal{P}_3\widetilde{\mathcal{C}}_0+\mathcal{P}_4(\mathcal{B}_0+\mathcal{C}_0+\widetilde{\mathcal{C}}_0)\big](\mathcal{P}_1+\mathcal{P}_2+\mathcal{P}_3+\mathcal{P}_4)\\
&\quad\ +\big(\mathcal{A}_1+\mathcal{P}_1\mathcal{B}_1^\top\big)\overline{\Sigma}_1+\big(\mathcal{A}_2+\mathcal{P}_1\mathcal{B}_2^\top\big)\overline{\Sigma}_2
 +\big(\mathcal{A}_3+\mathcal{P}_1\mathcal{B}_3^\top\big)\overline{\Sigma}_3+\mathcal{P}_2\mathcal{B}_1^\top(\widetilde{\Sigma}_1+\overline{\Sigma}_1)\\
&\quad\ +\big(\widehat{\mathcal{A}}_2+\mathcal{P}_2\mathcal{B}_2^\top\big)(\widetilde{\Sigma}_2+\overline{\Sigma}_2)+\big(\widehat{\mathcal{A}}_3+\mathcal{P}_2\mathcal{B}_3^\top\big)(\widetilde{\Sigma}_2+\overline{\Sigma}_2)
 +\big(\mathcal{P}_1\widetilde{\mathcal{B}}_1^\top+\mathcal{P}_3\mathcal{B}_1^\top+\mathcal{P}_3\widetilde{\mathcal{B}}_1^\top\big)(\widehat{\Sigma}_1+\overline{\Sigma}_1)\\
&\quad\ +\big(\mathcal{P}_1\widetilde{\mathcal{B}}_2^\top+\mathcal{P}_3\mathcal{B}_2^\top+\mathcal{P}_3\widetilde{\mathcal{B}}_2^\top\big)(\widehat{\Sigma}_2+\overline{\Sigma}_2)
 +\big(\mathcal{P}_1\widetilde{\mathcal{B}}_3^\top+\mathcal{P}_3\mathcal{B}_3^\top+\mathcal{P}_3\widetilde{\mathcal{B}}_3^\top\big)(\widehat{\Sigma}_3+\overline{\Sigma}_3)\\
&\quad\ +\big[\overline{\mathcal{A}}_1^\top+\mathcal{P}_1\overline{\mathcal{C}}_0+\mathcal{P}_2\big(\widetilde{\mathcal{B}}_1^\top+\overline{\mathcal{C}}_0\big)
 +\mathcal{P}_3\overline{\mathcal{C}}_0+\mathcal{P}_4\big(\mathcal{B}_1^\top+\widetilde{\mathcal{B}}_2^\top+\overline{\mathcal{D}}_0\big)\big](\Sigma_1+\widehat{\Sigma}_1+\widetilde{\Sigma}_1+\overline{\Sigma}_1)\\
&\quad\ +\big[\overline{\mathcal{A}}_2^\top+\mathcal{P}_1\overline{\mathcal{D}}_0+\mathcal{P}_2\big(\widetilde{\mathcal{B}}_2^\top+\overline{\mathcal{D}}_0\big)
 +\mathcal{P}_3\overline{\mathcal{D}}_0+\mathcal{P}_4\big(\mathcal{B}_2^\top+\widetilde{\mathcal{B}}_2^\top+\overline{\mathcal{D}}_0\big)\big](\Sigma_2+\widehat{\Sigma}_2+\widetilde{\Sigma}_2+\overline{\Sigma}_2)\\
&\quad\ +\big[\overline{\mathcal{A}}_3^\top+\mathcal{P}_1\overline{\mathcal{E}}_0+\mathcal{P}_2\big(\widetilde{\mathcal{B}}_3^\top+\overline{\mathcal{E}}_0\big)
 +\mathcal{P}_3\overline{\mathcal{E}}_0+\mathcal{P}_4\big(\mathcal{B}_3^\top+\widetilde{\mathcal{B}}_3^\top+\overline{\mathcal{E}}_0\big)\big](\Sigma_3+\widehat{\Sigma}_3+\widetilde{\Sigma}_3+\overline{\Sigma}_3),\\
&\mathcal{P}_1(T)=\mathcal{G}_2,\ \mathcal{P}_2(T)=0,\ \mathcal{P}_3(T)=0,\ \mathcal{P}_4(T)=0.
\end{aligned}
\right.
\end{equation}}
Notice that $\Sigma_i$ in the above depends on $P_i$, so the solvability of the above complicated and coupled system of Riccati's type equations is very difficult to obtain. We will not discuss this problem at the present paper for some technical reason and leave it open.

Finally, by (\ref{optimal control of the leader-2 dim}), (\ref{relation of X and Y}) and (\ref{Zi and its filter estimates}), we have
\begin{equation}\label{optimal control of the leader-feedback}
\begin{aligned}
u_2^*(t)=&-N_2^{-1}(t)\bigg[\mathcal{L}_4^\top\check{\hat{X}}(t)+\mathcal{C}_{05}^\top\check{Y}(t)+\mathcal{L}_{05}^\top\check{\hat{Y}}(t)
          +\sum\limits_{i=1}^3\mathcal{C}_{i5}^\top\check{Z}_i(t)+\sum\limits_{i=1}^3\mathcal{L}_{i5}^\top\check{\hat{Z}}_i(t)\bigg]\\
        =&-N_2^{-1}(t)\bigg\{\Big[\mathcal{C}_{05}^\top(\mathcal{P}_1+\mathcal{P}_3)
          +\sum\limits_{i=1}^3\mathcal{C}_{i5}^\top(\Sigma_i+\widetilde{\Sigma}_i)\Big]\check{X}(t)+\Big[\mathcal{L}_4^\top
          +\mathcal{C}_{05}^\top(\mathcal{P}_2+\mathcal{P}_4)\\
         &\qquad+\sum\limits_{i=1}^3\mathcal{C}_{i5}^\top(\widehat{\Sigma}_i+\overline{\Sigma}_i)+\mathcal{L}_{05}^\top(\mathcal{P}_1+\mathcal{P}_2+\mathcal{P}_3+\mathcal{P}_4)
          +\sum\limits_{i=1}^3\mathcal{L}_{i5}^\top(\Sigma_i+\widehat{\Sigma}_i+\widetilde{\Sigma}_i+\overline{\Sigma}_i)\Big]\check{\hat{X}}(t)\bigg\}.
\end{aligned}
\end{equation}
And the optimal ``state" $X=\big(x^*,p\big)^\top$ of the leader admits
{\small\begin{equation}\label{close-loop state of the leader}
\left\{
\begin{aligned}
  dX(t)&=\Big\{\big(\mathcal{A}_0+\mathcal{B}_0\mathcal{P}_1+\mathcal{B}_1^\top\Sigma_1+\mathcal{B}_2^\top\Sigma_2+\mathcal{B}_3^\top\Sigma_3\big)X(t)
         +\big(\widehat{\mathcal{A}}_0+\mathcal{B}_0\mathcal{P}_2+\mathcal{B}_1^\top\widehat{\Sigma}_1+\mathcal{B}_2^\top\widehat{\Sigma}_2\\
       &\quad\ +\mathcal{B}_3^\top\widehat{\Sigma}_3\big)\hat{X}(t)+\big[\mathcal{B}_0\mathcal{P}_3+\mathcal{C}_0(\mathcal{P}_1+\mathcal{P}_3)
        +\mathcal{B}_1^\top\widetilde{\Sigma}_1+\mathcal{B}_2^\top\widetilde{\Sigma}_2+\mathcal{B}_3^\top\widetilde{\Sigma}_3+\widetilde{\mathcal{B}}_1^\top(\Sigma_1+\widetilde{\Sigma}_1)\\
       &\quad\ +\widetilde{\mathcal{B}}_2^\top(\Sigma_2+\widetilde{\Sigma}_2)+\widetilde{\mathcal{B}}_3^\top(\Sigma_3+\widetilde{\Sigma}_3)\big]\check{X}(t)
        +\big[\overline{\mathcal{A}}_0+\mathcal{B}_0\mathcal{P}_4+\mathcal{C}_0(\mathcal{P}_2+\mathcal{P}_4)\\
       &\quad\ +\widetilde{\mathcal{C}}_0(\mathcal{P}_1+\mathcal{P}_2+\mathcal{P}_3+\mathcal{P}_4)+\mathcal{B}_1^\top\overline{\Sigma}_1+\mathcal{B}_2^\top\overline{\Sigma}_2+\mathcal{B}_3^\top\overline{\Sigma}_3
        +\widetilde{\mathcal{B}}_1^\top(\widehat{\Sigma}_1+\overline{\Sigma}_1)+\widetilde{\mathcal{B}}_2^\top(\widehat{\Sigma}_2+\overline{\Sigma}_2)\\
       &\quad\ +\widetilde{\mathcal{B}}_3^\top(\widehat{\Sigma}_3+\overline{\Sigma}_3)+\overline{\mathcal{C}}_0(\Sigma_1+\widehat{\Sigma}_1+\widetilde{\Sigma}_1+\overline{\Sigma}_1)
        +\overline{\mathcal{D}}_0(\Sigma_2+\widehat{\Sigma}_2+\widetilde{\Sigma}_2+\overline{\Sigma}_2)\\
       &\quad\ +\overline{\mathcal{E}}_0(\Sigma_3+\widehat{\Sigma}_3+\widetilde{\Sigma}_3+\overline{\Sigma}_3)\big]\check{\hat{X}}(t)\Big\}dt
        +\sum\limits_{i=1}^3\Big\{\big(\mathcal{A}_i+\mathcal{B}_i\mathcal{P}_1+\mathcal{B}_i\Sigma_1+\mathcal{D}_i\Sigma_2+\mathcal{E}_i\Sigma_3\big)X(t)\\
       &\quad+\big(\widehat{\mathcal{A}}_i+\mathcal{B}_i\mathcal{P}_2+\mathcal{B}_i\widehat{\Sigma}_1+\mathcal{D}_i\widehat{\Sigma}_2+\mathcal{E}_i\widehat{\Sigma}_3\big)\hat{X}(t)
        +\big[\mathcal{B}_i\mathcal{P}_3+\widetilde{\mathcal{B}}_i(\mathcal{P}_1+\mathcal{P}_3)+\mathcal{B}_i\widetilde{\Sigma}_1+\mathcal{D}_i\widetilde{\Sigma}_2\\
       &\quad\ +\mathcal{E}_i\widetilde{\Sigma}_3+\widetilde{\mathcal{C}}_i(\Sigma_1+\widetilde{\Sigma}_1)+\widetilde{\mathcal{D}}_i(\Sigma_2+\widetilde{\Sigma}_2)
        +\widetilde{\mathcal{E}}_i(\Sigma_3+\widetilde{\Sigma}_3)\big]\check{X}(t)\\
       &\quad\ +\big[\overline{\mathcal{A}}_i+\mathcal{B}_i\mathcal{P}_4+\widetilde{\mathcal{B}}_i(\mathcal{P}_2+\mathcal{P}_4)
        +\overline{\mathcal{B}}_i(\mathcal{P}_1+\mathcal{P}_2+\mathcal{P}_3+\mathcal{P}_4)+\mathcal{B}_i\overline{\Sigma}_1+\mathcal{D}_i\overline{\Sigma}_2+\mathcal{E}_i\overline{\Sigma}_3\\
       &\quad\ +\widetilde{\mathcal{C}}_i(\widehat{\Sigma}_1+\overline{\Sigma}_1)+\widetilde{\mathcal{D}}_i(\widehat{\Sigma}_2+\overline{\Sigma}_2)
        +\widetilde{\mathcal{E}}_i(\widehat{\Sigma}_3+\overline{\Sigma}_3)+\overline{\mathcal{C}}_i(\Sigma_1+\widehat{\Sigma}_1+\widetilde{\Sigma}_1+\overline{\Sigma}_1)\\
       &\quad\ +\overline{\mathcal{D}}_i(\Sigma_2+\widehat{\Sigma}_2+\widetilde{\Sigma}_2+\overline{\Sigma}_2)
        +\overline{\mathcal{E}}_i(\Sigma_3+\widehat{\Sigma}_3+\widetilde{\Sigma}_3+\overline{\Sigma}_3)\big]\check{\hat{X}}(t)\Big\}dW_i(t),\ t\in[0,T],\\
   X(0)&=X_0,
\end{aligned}
\right.
\end{equation}}
where $\hat{X}$ is determined by
{\small\begin{equation}\label{hat-X-optimal feedback}
\left\{
\begin{aligned}
d\hat{X}(t)&=\Big\{\big[\mathcal{A}_0+\widehat{\mathcal{A}}_0+\mathcal{B}_0(\mathcal{P}_1+\mathcal{P}_2)+\mathcal{B}_1^\top(\Sigma_1+\widehat{\Sigma}_1)+\mathcal{B}_2^\top(\Sigma_2+\widehat{\Sigma}_2)
            +\mathcal{B}_3^\top(\Sigma_3+\widehat{\Sigma}_3)\big]\hat{X}(t)\\
           &\quad\ +\big[\overline{\mathcal{A}}_0+\mathcal{B}_0(\mathcal{P}_3+\mathcal{P}_4)+\mathcal{B}_1^\top(\widetilde{\Sigma}_1+\overline{\Sigma}_1)+\mathcal{B}_2^\top(\widetilde{\Sigma}_2+\overline{\Sigma}_2)
            +\mathcal{B}_3^\top(\widetilde{\Sigma}_3+\overline{\Sigma}_3)\\
           &\quad\ +(\mathcal{C}_0+\widetilde{\mathcal{C}}_0)(\mathcal{P}_1+\mathcal{P}_2+\mathcal{P}_3+\mathcal{P}_4)
            +(\widetilde{\mathcal{B}}_1^\top+\widetilde{\mathcal{B}}_2^\top+\widetilde{\mathcal{B}}_3^\top+\overline{\mathcal{C}}_0)(\Sigma_1+\widehat{\Sigma}_1+\widetilde{\Sigma}_1+\overline{\Sigma}_1)\\
           &\quad\ +\overline{\mathcal{D}}_0(\Sigma_2+\widehat{\Sigma}_2+\widetilde{\Sigma}_2+\overline{\Sigma}_2)
            +\overline{\mathcal{E}}_0(\Sigma_3+\widehat{\Sigma}_3+\widetilde{\Sigma}_3+\overline{\Sigma}_3)\big]\check{\hat{X}}(t)\Big\}dt\\
           &\quad+\sum\limits_{i=1,3}\Big\{\big[\mathcal{A}_i+\widehat{\mathcal{A}}_i+\mathcal{B}_i(\mathcal{P}_1+\mathcal{P}_2)+\mathcal{B}_i(\Sigma_1+\widehat{\Sigma}_1)+\mathcal{D}_i(\Sigma_2+\widehat{\Sigma}_2)
            +\mathcal{E}_i(\Sigma_3+\widehat{\Sigma}_3)\big]\hat{X}(t)\\
           &\quad\ +\big[\overline{\mathcal{A}}_i+\mathcal{B}_i(\mathcal{P}_3+\mathcal{P}_4)
            +\mathcal{B}_i(\widetilde{\Sigma}_1+\overline{\Sigma}_1)+\mathcal{D}_i(\widetilde{\Sigma}_2+\overline{\Sigma}_2)+\mathcal{E}_i(\widetilde{\Sigma}_3+\overline{\Sigma}_3)\\
           &\quad\ +(\widetilde{\mathcal{B}}_i+\overline{\mathcal{B}}_i)(\mathcal{P}_1+\mathcal{P}_2+\mathcal{P}_3+\mathcal{P}_4)
            +(\widetilde{\mathcal{C}}_i+\overline{\mathcal{C}}_i)(\Sigma_1+\widehat{\Sigma}_1+\widetilde{\Sigma}_1+\overline{\Sigma}_1)\\
           &\quad\ +(\widetilde{\mathcal{D}}_i+\overline{\mathcal{D}}_i)(\Sigma_2+\widehat{\Sigma}_2+\widetilde{\Sigma}_2+\overline{\Sigma}_2)
            +(\widetilde{\mathcal{E}}_i+\overline{\mathcal{E}}_i)(\Sigma_3+\widehat{\Sigma}_3+\widetilde{\Sigma}_3+\overline{\Sigma}_3)\big]\check{\hat{X}}(t)\Big\}dW_i(t),\ t\in[0,T],\\
 \hat{X}(0)&=X_0,
\end{aligned}
\right.
\end{equation}}
$\check{X}$ is governed by
{\small\begin{equation}\label{check-X-optimal feedback}
\left\{
\begin{aligned}
d\check{X}(t)&=\Big\{\big[\mathcal{A}_0+(\mathcal{B}_0+\mathcal{C}_0)(\mathcal{P}_1+\mathcal{P}_3)
              +(\mathcal{B}_1+\widetilde{\mathcal{B}}_1)^\top(\Sigma_1+\widetilde{\Sigma}_1)+(\mathcal{B}_2+\widetilde{\mathcal{B}}_2)^\top(\Sigma_2+\widetilde{\Sigma}_2)\\
             &\qquad+(\mathcal{B}_3+\widetilde{\mathcal{B}}_3)^\top(\Sigma_3+\widetilde{\Sigma}_3)\big]\check{X}(t)
              +\big[\widehat{\mathcal{A}}_0+\overline{\mathcal{A}}_0+(\mathcal{B}_0+\mathcal{C}_0)(\mathcal{P}_2+\mathcal{P}_4)\\
             &\qquad+(\mathcal{B}_1+\widetilde{\mathcal{B}}_1)^\top(\widehat{\Sigma}_1+\overline{\Sigma}_1)+(\mathcal{B}_2+\widetilde{\mathcal{B}}_2)^\top(\widehat{\Sigma}_2+\overline{\Sigma}_2)
              +(\mathcal{B}_3+\widetilde{\mathcal{B}}_3)^\top(\widehat{\Sigma}_3+\overline{\Sigma}_3)\\
             &\qquad+\widetilde{\mathcal{C}}_0(\mathcal{P}_1+\mathcal{P}_2+\mathcal{P}_3+\mathcal{P}_4)+\overline{\mathcal{C}}_0(\Sigma_1+\widehat{\Sigma}_1+\widetilde{\Sigma}_1+\overline{\Sigma}_1)
              +\overline{\mathcal{D}}_0(\Sigma_2+\widehat{\Sigma}_2+\widetilde{\Sigma}_2+\overline{\Sigma}_2)\\
             &\qquad+\overline{\mathcal{E}}_0(\Sigma_3+\widehat{\Sigma}_3+\widetilde{\Sigma}_3+\overline{\Sigma}_3)\big]\check{\hat{X}}(t)\Big\}dt
              +\sum\limits_{i=2,3}\Big\{\big[\mathcal{A}_i+(\mathcal{B}_i+\widetilde{\mathcal{B}}_i)(\mathcal{P}_1+\mathcal{P}_3)\\
             &\qquad+(\mathcal{B}_i+\widetilde{\mathcal{C}}_i)(\Sigma_1+\widetilde{\Sigma}_1)+(\mathcal{D}_i+\widetilde{\mathcal{D}}_i)(\Sigma_2+\widetilde{\Sigma}_2)
              +(\mathcal{E}_i+\widetilde{\mathcal{E}}_i)(\Sigma_3+\widetilde{\Sigma}_3)\big]\check{X}(t)\\
             &\qquad+\big[\widehat{\mathcal{A}}_i+\overline{\mathcal{A}}_i+(\mathcal{B}_i+\widetilde{\mathcal{B}}_i)(\mathcal{P}_2+\mathcal{P}_4)
              +\overline{\mathcal{B}}_i(\mathcal{P}_1+\mathcal{P}_2+\mathcal{P}_3+\mathcal{P}_4)\\
             &\qquad+(\mathcal{B}_i+\widetilde{\mathcal{C}}_i)(\widehat{\Sigma}_1+\overline{\Sigma}_1)+(\mathcal{D}_i+\widetilde{\mathcal{D}}_i)(\widehat{\Sigma}_2+\overline{\Sigma}_2)
              +(\mathcal{E}_i+\widetilde{\mathcal{E}}_i)(\widehat{\Sigma}_3+\overline{\Sigma}_3)\\
             &\qquad+\overline{\mathcal{C}}_i(\Sigma_1+\widehat{\Sigma}_1+\widetilde{\Sigma}_1+\overline{\Sigma}_1)+\overline{\mathcal{D}}_i(\Sigma_2+\widehat{\Sigma}_2+\widetilde{\Sigma}_2+\overline{\Sigma}_2)\\
             &\qquad+\overline{\mathcal{E}}_i(\Sigma_3+\widehat{\Sigma}_3+\widetilde{\Sigma}_3+\overline{\Sigma}_3)\big]\check{\hat{X}}(t)\Big\}dW_i(t),\ t\in[0,T],\\
 \check{X}(0)=&\ X_0,
\end{aligned}
\right.
\end{equation}}
and $\check{\hat{X}}$ is given by
{\small\begin{equation}\label{check-hat-X-optimal feedback}
\left\{
\begin{aligned}
d\check{\hat{X}}(t)&=\big[\mathcal{A}_0+\widehat{\mathcal{A}}_0+\overline{\mathcal{A}}_0+(\mathcal{B}_0+\mathcal{C}_0+\widetilde{\mathcal{C}}_0)(\mathcal{P}_1+\mathcal{P}_2+\mathcal{P}_3+\mathcal{P}_4)\\
                   &\quad\ +(\mathcal{B}_1^\top+\widetilde{\mathcal{B}}_1^\top+\overline{\mathcal{C}}_0)(\Sigma_1+\widehat{\Sigma}_1+\widetilde{\Sigma}_1+\overline{\Sigma}_1)
                    +(\mathcal{B}_2^\top+\widetilde{\mathcal{B}}_2^\top+\overline{\mathcal{D}}_0)(\Sigma_2+\widehat{\Sigma}_2+\widetilde{\Sigma}_2+\overline{\Sigma}_2)\\
                   &\quad\ +(\mathcal{B}_3^\top+\widetilde{\mathcal{B}}_3^\top+\overline{\mathcal{E}}_0)(\Sigma_3+\widehat{\Sigma}_3+\widetilde{\Sigma}_3+\overline{\Sigma}_3)\big]\check{\hat{X}}(t)dt\\
                   &\quad+\big[\mathcal{A}_3+\widehat{\mathcal{A}}_3+\overline{\mathcal{A}}_3+(\mathcal{B}_3+\widetilde{\mathcal{B}}_3+\overline{\mathcal{B}}_3)(\mathcal{P}_1+\mathcal{P}_2+\mathcal{P}_3+\mathcal{P}_4)\\
                   &\quad\ +(\mathcal{B}_3+\widetilde{\mathcal{C}}_3+\overline{\mathcal{C}}_3)(\Sigma_1+\widehat{\Sigma}_1+\widetilde{\Sigma}_1+\overline{\Sigma}_1)
                    +(\mathcal{D}_3+\widetilde{\mathcal{D}}_3+\overline{\mathcal{D}}_3)(\Sigma_2+\widehat{\Sigma}_2+\widetilde{\Sigma}_2+\overline{\Sigma}_2)\\
                   &\quad\ +(\mathcal{E}_3+\widetilde{\mathcal{E}}_3+\overline{\mathcal{E}}_3)(\Sigma_3+\widehat{\Sigma}_3+\widetilde{\Sigma}_3+\overline{\Sigma}_3)\big]\check{\hat{X}}(t)dW_3(t),\ t\in[0,T],\\
 \check{\hat{X}}(0)&=X_0.
\end{aligned}
\right.
\end{equation}}

We summarize the above argument in the following theorem.

\vspace{1mm}

\noindent{\bf Theorem 3.2}\quad{\it Let ${\bf (A2.2)\sim(A2.5)}$ hold and $(\mathcal{P}_1(\cdot),\mathcal{P}_2(\cdot),\mathcal{P}_3(\cdot),\mathcal{P}_4(\cdot))$ satisfy (\ref{system of Riccati equation}), $\check{\hat{X}}(\cdot)$ be the $\mathcal{G}^1_t\cap\mathcal{G}^2_t$-adapted solution to (\ref{check-hat-X-optimal feedback}), $\check{X}(\cdot)$ be the $\mathcal{G}^2_t$-adapted solution to (\ref{check-X-optimal feedback}), $\hat{X}(\cdot)$ be the $\mathcal{G}^1_t$-adapted solution to (\ref{hat-X-optimal feedback}), and $X(\cdot)$ be the $\mathcal{F}_t$-adapted solution to (\ref{close-loop state of the leader}). Define $(Y(\cdot),Z_1(\cdot),Z_2(\cdot),Z_3(\cdot))$ by (\ref{relation of X and Y}) and (\ref{comparing dW1,dW2,dW3-leader}), respectively. Then (\ref{optimality system-leader-2 dim}) holds, and $u_2^*(\cdot)$ given by (\ref{optimal control of the leader-feedback}) is a feedback optimal control of the leader.}

\vspace{2mm}

Finally, the optimal control $u_1^*(\cdot)$ of the follower can also be represented in a ``nonanticipating" way. In fact, by (\ref{optimal control of the follower-feedback}), noting (\ref{optimal control of the leader-feedback}), (\ref{new state}), (\ref{relation of X and Y}) and (\ref{Zi and its filter estimates}), we obtain
\begin{equation*}
\begin{aligned}
u_1^*(t)=
 &-\Big[N_1(t)+\sum\limits_{i=1}^3B_i^\top(t)P_1(t)B_i(t)\Big]^{-1}\bigg\{\Big[B_0^\top(t)P_1(t)+\sum\limits_{i=1}^3B_i^\top(t)P_1(t)A_i(t)\Big]\hat{x}^{u_1^*,u_2^*}(t)\\
 &\quad+B_0^\top(t)P_1(t)\hat{\phi}^*(t)+B_1^\top(t)\hat{\beta}_1^*(t)+B_3^\top(t)\hat{\beta}_3^*(t)+\Big[\sum\limits_{i=1}^3B_i^\top(t)P_1(t)C_i(t)\Big]\hat{u}_2^*(t)\bigg\}\\
\end{aligned}
\end{equation*}
\begin{equation}\label{optimal control of the follower-final}
\begin{aligned}
=&-\Big[N_1(t)+\sum\limits_{i=1}^3B_i^\top(t)P_1(t)B_i(t)\Big]^{-1}\bigg\{\left(\begin{array}{cc}B_0^\top(t)P_1(t)+\sum\limits_{i=1}^3B_i^\top(t)P_1(t)A_i(t)&0\end{array}\right)\hat{X}(t)\\
 &\qquad+\left(\begin{array}{cc}0&B_0^\top(t)P_1(t)\end{array}\right)\hat{Y}(t)
  +\left(\begin{array}{cc}0&B_1^\top(t)\end{array}\right)\hat{Z}_1(t)+\left(\begin{array}{cc}0&B_3^\top(t)\end{array}\right)\hat{Z}_3(t)\\
 &\qquad-\Big(\sum\limits_{i=1}^3B_i^\top(t)P_1(t)C_i(t)\Big)N_2^{-1}(t)\Big[\mathcal{C}_{05}^\top(\mathcal{P}_1+\mathcal{P}_3)
  +\sum\limits_{i=1}^3\mathcal{C}_{i5}^\top(\Sigma_i+\widetilde{\Sigma}_i)\\
 &\qquad+\mathcal{L}_4^\top+\mathcal{C}_{05}^\top(\mathcal{P}_2+\mathcal{P}_4)
  +\sum\limits_{i=1}^3\mathcal{C}_{i5}^\top(\widehat{\Sigma}_i+\overline{\Sigma}_i)+\mathcal{L}_{05}^\top(\mathcal{P}_1+\mathcal{P}_2+\mathcal{P}_3+\mathcal{P}_4)\\
 &\qquad+\sum\limits_{i=1}^3\mathcal{L}_{i5}^\top(\Sigma_i+\widehat{\Sigma}_i+\widetilde{\Sigma}_i+\overline{\Sigma}_i)\Big]\check{\hat{X}}(t)\bigg\}\\
=&-\Big[N_1(t)+\sum\limits_{i=1}^3B_i^\top(t)P_1(t)B_i(t)\Big]^{-1}\bigg[\left(\begin{array}{cc}B_0^\top(t)P_1(t)+\sum\limits_{i=1}^3B_i^\top(t)P_1(t)A_i(t)&0\end{array}\right)\\
 &\quad+\left(\begin{array}{cc}0&B_0^\top(t)P_1(t)\end{array}\right)(\mathcal{P}_1+\mathcal{P}_2)
  +\left(\begin{array}{cc}0&B_1^\top(t)\end{array}\right)(\Sigma_1+\widehat{\Sigma}_1)\\
 &\quad+\left(\begin{array}{cc}0&B_3^\top(t)\end{array}\right)(\Sigma_3+\widetilde{\Sigma}_3)\bigg]\hat{X}(t)\\
 &-\Big[N_1(t)+\sum\limits_{i=1}^3B_i^\top(t)P_1(t)B_i(t)\Big]^{-1}\bigg\{\left(\begin{array}{cc}0&B_0^\top(t)P_1(t)\end{array}\right)(\mathcal{P}_3+\mathcal{P}_4)\\
 &\quad+\left(\begin{array}{cc}0&B_1^\top(t)\end{array}\right)(\widetilde{\Sigma}_1+\overline{\Sigma}_1)+\left(\begin{array}{cc}0&B_3^\top(t)\end{array}\right)(\widetilde{\Sigma}_3+\overline{\Sigma}_3)\\
 &\quad-N_2^{-1}(t)\Big(\sum\limits_{i=1}^3B_i^\top(t)P_1(t)C_i(t)\Big)\Big[\mathcal{L}_4^\top+(\mathcal{C}_{05}+\mathcal{L}_{05})^\top(\mathcal{P}_1+\mathcal{P}_2+\mathcal{P}_3+\mathcal{P}_4)\\
 &\quad+\sum\limits_{i=1}^3(\mathcal{C}_{i5}+\mathcal{L}_{i5})^\top(\Sigma_i+\widetilde{\Sigma}_i+\widehat{\Sigma}_i+\overline{\Sigma}_i)\Big]\bigg\}\check{\hat{X}}(t),
\end{aligned}
\end{equation}
which is observable for the follower.

Up to now, the Stackelberg equilibrium strategy $(u_1^*(\cdot),u_2^*(\cdot))$ is obtained, which is represented as the state estimate feedback form in (\ref{optimal control of the follower-final}) and (\ref{optimal control of the leader-feedback}).

\section{A special solvable case: Control independent diffusions}

In this section, we consider the problem for the special $n=1$ case with control independent diffusions and constant parameters. In this case, the problem can be completely solved together with the solution to the system of Riccati equations.

We consider the scalar state process $x^{u_1,u_2}(\cdot)$ which satisfies the linear SDE
\begin{equation}\label{state equation of the follower-control independent diffusion}
\left\{
\begin{aligned}
     dx^{u_1,u_2}(t)&=\big[A_0x^{u_1,u_2}(t)+B_0u_1(t)+C_0u_2(t)\big]dt+A_1x^{u_1,u_2}(t)dW_1(t)\\
                    &\quad+A_2x^{u_1,u_2}(t)dW_2(t)+A_3x^{u_1,u_2}(t)dW_3(t),\ t\in[0,T],\\
      x^{u_1,u_2}(0)&=x_0.
\end{aligned}
\right.
\end{equation}
Here $u_1(\cdot)$ and $u_2(\cdot)$ are both scalar-valued and $A_0,B_0,C_0,A_1,A_2,A_3$ are constants. We define the admissible control sets $\mathcal{U}_1,\mathcal{U}_2$ as in Section 2.

In step 1, for any chosen $u_2(\cdot)$, the follower wishes to select a $u_1^*(\cdot)\in\mathcal{U}_1$ to minimize the cost functional
\begin{equation}\label{cost functional of the follower-control independent diffusion}
\begin{aligned}
 J_1(u_1(\cdot),u_2(\cdot))=\frac{1}{2}\mathbb{E}\bigg[\int_0^T\Big(Q_1\big|x^{u_1,u_2}(t)\big|^2+N_1u_1^2(t)\Big)dt+G_1\big|x^{u_1,u_2}(T)\big|^2\bigg].
\end{aligned}
\end{equation}
Here $Q_1,G_1\geq0,N_1\neq0$ are constants. In step 2, after the follower's optimal control $u_1^*(\cdot)$ is announced, the leader would like to choose a $u_2^*(\cdot)\in\mathcal{U}_2$ to minimize
\begin{equation}\label{cost functional of the leader-control independent diffusion}
\begin{aligned}
J_2(u_1^*(\cdot),u_2(\cdot))=\frac{1}{2}\mathbb{E}\bigg[\int_0^T\Big(Q_2\big|x^{u_1^*,u_2}(t)\big|^2+N_2u_2^2(t)\Big)dt+G_2\big|x^{u_1^*,u_2}(T)\big|^2\bigg],
\end{aligned}
\end{equation}
where $Q_2,G_2\geq0,N_2\neq0$ are constants. We wish to find the Stackelberg equilibrium strategy $(u_1^*(\cdot),u_2^*(\cdot))\in\mathcal{U}_1\times\mathcal{U}_2$.

\subsection{Problem of The Follower}

For given control $u_2(\cdot)$, let $u_1^*(\cdot)$ be a $\mathcal{G}^1_t$-adapted optimal control of the follower, and the corresponding optimal state is $x^{u_1^*,u_2}(\cdot)$. Now the follower's Hamiltonian function (\ref{Hamiltonian function of the follower}) writes
\begin{equation}\label{Hamiltonian function of the follower-control independent diffusion}
\begin{aligned}
&H_1\big(t,x,u_1,u_2,q,k_1,k_2,k_3\big)\triangleq q(A_0x+B_0u_1+C_0u_2)\\
&\quad+A_1k_1x+A_2k_2x+A_3k_3x-\frac{1}{2}Q_1x^2-\frac{1}{2}N_1u_1^2.
\end{aligned}
\end{equation}
And (\ref{optimal control of the follower}) yields that
\begin{equation}\label{optimal control of the follower-control independent diffusion}
0=N_1u_1^*(t)-B_0\hat{q}(t),
\end{equation}
where the $\mathcal{F}_t$-adapted process quadruple $(q(\cdot),k_1(\cdot),k_2(\cdot),k_3(\cdot))$ satisfies the adjoint BSDE
\begin{equation}\label{adjoint equation of the follower-control independent diffusion}
\left\{
\begin{aligned}
-dq(t)=&\big[A_0q(t)+A_1k_1+A_2k_2+A_3k_3-Q_1x^{u_1^*,u_2}(t)\big]dt\\
      &-k_1dW_1(t)-k_2dW_2(t)-k_3dW_3(t),\ t\in[0,T],\\
  q(T)=&-G_1x^{u_1^*,u_2}(T),
\end{aligned}
\right.
\end{equation}
which is a special case of (\ref{adjoint equation of the follower}). Repeat the same approach as in Section 2.1, we obtain the following theorem.

\vspace{1mm}

\noindent{\bf Theorem 3.1}\quad{\it Let $P(\cdot)$ satisfy
\begin{equation}\label{Riccati equation-control independent diffusion}
\left\{
\begin{aligned}
 &\dot{P}(t)+\big(2A_0+A_1^2+A_2^2+A_3^2\big)P(t)-N_1^{-1}B_0^2P^2(t)+Q_1=0,\ t\in[0,T],\\
 &P(T)=G_1.
\end{aligned}
\right.
\end{equation}
For chosen $u_2(\cdot)$ of the leader, $u_1^*(\cdot)$ defined by
\begin{equation}\label{optimal control of the follower-feedback-control independent diffusion}
\begin{aligned}
u_1^*(t)=-N_1^{-1}B_0\big[P(t)\hat{x}^{u_1^*,u_2}(t)+\hat{\phi}(t)\big]
\end{aligned}
\end{equation}
is a feedback optimal control of the follower, where $(\hat{x}^{u_1^*,u_2}(\cdot),\hat{\phi}(\cdot),\hat{\beta}_1(\cdot),\hat{\beta}_3(\cdot))$ is the unique $\mathcal{G}^1_t$-adapted solution to \begin{equation}\label{FBSDFE-control independent diffusion}
\left\{
\begin{aligned}
  d\hat{x}^{u_1^*,u_2}(t)&=\big[\big(A_0-N_1^{-1}B_0^2P(t)\big)\hat{x}^{u_1^*,u_2}(t)-N_1^{-1}B_0^2\hat{\phi}(t)+C_0\hat{u}_2(t)\big]dt\\
                         &\quad+A_1\hat{x}^{u_1^*,u_2}(t)dW_1(t)+A_3\hat{x}^{u_1^*,u_2}(t)dW_3(t),\\
           -d\hat{\phi}(t)&=\big[\big(A_0-N_1^{-1}B_0^2P(t)\big)\hat{\phi}(t)+A_1\hat{\beta}_1(t)+A_3\hat{\beta}_3(t)+P(t)C_0\hat{u}_2(t)\big]dt\\
                         &\quad-\hat{\beta}_1(t)dW_1(t)-\hat{\beta}_3(t)dW_3(t),\ t\in[0,T],\\
      \hat{x}^{u_1,u_2}(0)&=x_0,\quad \hat{\phi}(T)=0.
\end{aligned}
\right.
\end{equation}}

\subsection{Problem of The Leader}

In the following, the leader keeps in mind that the follower takes $u_1^*(\cdot)$ by (\ref{optimal control of the follower-feedback-control independent diffusion}), then his state equation (\ref{state equation-leader}) writes
\begin{equation}\label{state equation of the leader-control independent diffusion}
\left\{
\begin{aligned}
    dx^{u_2}(t)&=\big[A_0x^{u_2}(t)-N_1^{-1}B_0^2P(t)\hat{x}^{u_2}(t)-N_1^{-1}B_0^2\hat{\phi}(t)+C_0u_2(t)\big]dt\\
               &\quad+A_1x^{u_2}(t)dW_1(t)+A_2x^{u_2}(t)dW_2(t)+A_3x^{u_2}(t)dW_3(t),\\
-d\hat{\phi}(t)&=\Big\{\big[A_0-N_1^{-1}B_0^2P(t)\big]\hat{\phi}(t)+A_1\hat{\beta}_1(t)+A_3\hat{\beta}_3(t)+P(t)C_0\hat{u}_2(t)\Big\}dt\\
               &\quad-\hat{\beta}_1(t)dW_1(t)-\hat{\beta}_3(t)dW_3(t),\ t\in[0,T],\\
     x^{u_2}(0)&=x_0,\quad \hat{\phi}(T)=0.
\end{aligned}
\right.
\end{equation}
The problem of the leader is to select a $\mathcal{G}^2_t$-adapted optimal control $u_2^*(\cdot)$ such that the cost functional
\begin{equation}\label{cost functional of the leader-control independent diffusion-simple}
\begin{aligned}
J_2(u_2(\cdot))=\frac{1}{2}\mathbb{E}\left[\int_0^T\big[Q_2|x^{u_2}(t)|^2+N_2u_2^2(t)\big]dt+G_2|x^{u_2}(T)|^2\right]
\end{aligned}
\end{equation}
is minimized.

Suppose that there exists a $\mathcal{G}^2_t$-adapted optimal control $u_2^*(\cdot)$ of the leader, and his optimal state is $(x^*(\cdot),\hat{\phi}^*(\cdot),\hat{\beta}_1^*(\cdot),\hat{\beta}_3^*(\cdot))\equiv(x^{u_2^*}(\cdot),\hat{\phi}^*(\cdot),\hat{\beta}_1^*(\cdot),\hat{\beta}_3^*(\cdot))$. Now the leader's Hamiltonian function
(\ref{Hamiltonian function of the leader}) reduces to
\begin{equation}\label{Hamiltonian function of the leader-control independent diffusion}
\begin{aligned}
&H_2\big(t,x^{u_2},u_2,\phi,\beta_1,\beta_3;p,y,z_1,z_2,z_3\big)\triangleq y\big[A_0x^{u_2}-N_1^{-1}B_0^2P(t)\hat{x}^{u_2}-N_1^{-1}B_0^2\hat{\phi}+C_0u_2\big]\\
&\quad+p\big\{[A_0-N_1^{-1}B_0^2P(t)]\hat{\phi}+A_1\hat{\beta}_1+A_3\hat{\beta}_3+P(t)C_0\hat{u}_2\big\}\\
&\quad+z_1A_1x^{u_2}+z_2A_2x^{u_2}+z_3A_3x^{u_2}+\frac{1}{2}\big[Q_2|x^{u_2}|^2+N_2u_2^2\big],
\end{aligned}
\end{equation}
where the $\mathcal{F}_t$-adapted process quintuple $(p(\cdot),y(\cdot),z_1(\cdot),z_2(\cdot),z_3(\cdot))$ satisfies the adjoint equation
\begin{equation}\label{adjoint equation of the leader-control independent diffusion}
\left\{
\begin{aligned}
  dp(t)&=\Big\{-N_1^{-1}B_0^2y(t)+\big[A_0-N_1^{-1}B_0^2P(t)\big]p(t)\Big\}dt+A_1p(t)dW_1(t)+A_3p(t)dW_3(t),\\
 -dy(t)&=\big[A_0y(t)-N_1^{-1}B_0^2P(t)\hat{y}(t)+A_1z_1(t)+A_2z_2(t)+A_3z_3(t)+Q_2x^*(t)\big]dt\\
       &\quad-z_1(t)dW_1(t)-z_2(t)dW_2(t)-z_3(t)dW_3(t),\ t\in[0,T],\\
   p(0)&=0,\quad y(T)=G_2x^*(T),
\end{aligned}
\right.
\end{equation}
which is a special case of (\ref{adjoint equation of the leader}). Similarly, we have
\begin{equation}\label{optimal control of the leader-control independent diffusion}
\begin{aligned}
u_2^*(t)=-N_2^{-1}C_0\big[\check{y}(t)+P(t)\check{\hat{p}}(t)\big],
\end{aligned}
\end{equation}
where
\begin{equation}\label{optimal filter equation-leader-control independent diffusion}
\left\{
\begin{aligned}
          d\check{x}^*(t)&=\big[A_0\check{x}^*(t)-N_1^{-1}B_0^2P(t)\check{\hat{x}}^*(t)-N_1^{-1}B_0^2\check{\hat{\phi}}^*(t)-N_2^{-1}C_0^2\check{y}(t)\\
                         &\quad-N_2^{-1}C_0^2P(t)\check{\hat{p}}(t)\big]dt+A_2\check{x}^*(t)dW_2(t)+A_3\check{x}^*(t)dW_3(t),\\
            d\check{p}(t)&=\Big\{-N_1^{-1}B_0^2\check{y}(t)+\big[A_0-N_1^{-1}B_0^2P(t)\big]\check{p}(t)\Big\}dt+A_3\check{p}(t)dW_3(t),\\
           -d\check{y}(t)&=\big[A_0\check{y}(t)-N_1^{-1}B_0^2P(t)\check{\hat{y}}(t)+A_1\check{z}_1(t)+A_2\check{z}_2(t)+A_3\check{z}_3(t)+Q_2\check{x}^*(t)\big]dt\\
                         &\quad-\check{z}_2(t)dW_2(t)-\check{z}_3(t)dW_3(t),\\
-d\check{\hat{\phi}}^*(t)&=\Big\{\big[A_0-N_1^{-1}B_0^2P(t)\big]\check{\hat{\phi}}^*(t)+A_1\check{\hat{\beta}}_1^*(t)+A_3\check{\hat{\beta}}_3^*(t)
                          -N_2^{-1}C_0^2P(t)\check{\hat{y}}(t)\\
                         &\quad-N_2^{-1}C_0^2P^2(t)\check{\hat{p}}(t)\Big\}dt-\check{\hat{\beta}}_3^*(t)dW_3(t),\ t\in[0,T],\\
           \check{x}^*(0)&=x_0,\quad \check{p}(0)=0,\quad \check{y}(T)=G_2\check{x}^*(T),\quad \check{\hat{\phi}}^*(T)=0,
\end{aligned}
\right.
\end{equation}
and
\begin{equation}\label{optimal filter equation-leader and follower-control independent diffusion}
\left\{
\begin{aligned}
    d\check{\hat{x}}^*(t)&=\Big\{\big[A_0-N_1^{-1}B_0^2P(t)\big]\check{\hat{x}}^*(t)-N_1^{-1}B_0^2\check{\hat{\phi}}^*(t)-N_2^{-1}C_0^2\check{\hat{y}}(t)\\
                         &\quad-N_2^{-1}C_0^2P(t)\check{\hat{p}}(t)\Big\}dt+A_3\check{\hat{x}}^*(t)dW_3(t),\\
      d\check{\hat{p}}(t)&=\Big\{-N_1^{-1}B_0^2\check{\hat{y}}(t)+\big[A_0-N_1^{-1}B_0^2P(t)\big]\check{\hat{p}}(t)\Big\}dt+A_3\check{\hat{p}}(t)dW_3(t),\\
     -d\check{\hat{y}}(t)&=\Big\{\big[A_0-N_1^{-1}B_0^2P(t)\big]\check{\hat{y}}(t)+A_1\check{\hat{z}}_1(t)+A_2\check{\hat{z}}_2(t)+A_3\check{\hat{z}}_3(t)+Q_2\check{\hat{x}}^*(t)\Big\}dt\\
                         &\quad-\check{\hat{z}}_3(t)dW_3(t),\ t\in[0,T],\\
     \check{\hat{x}}^*(0)&=x_0,\quad \check{\hat{p}}(0)=0,\quad \check{\hat{y}}(T)=G_2\check{\hat{x}}^*(T).
\end{aligned}
\right.
\end{equation}
The solvability of (\ref{adjoint equation of the leader-control independent diffusion}), (\ref{optimal filter equation-leader-control independent diffusion}) and (\ref{optimal filter equation-leader and follower-control independent diffusion}) will be proven in the following context. As in Section 2, we will proceed to represent $u_2^*(\cdot)$ of (\ref{optimal control of the leader-control independent diffusion}) as the state estimate feedback form, via some Riccati type equations.

First, we rewrite the optimal state equation (\ref{state equation of the leader-control independent diffusion}) as
\begin{equation}\label{optimal state equation for the leader-control independent diffusion case-FBSDE}
\left\{
\begin{aligned}
          dx^*(t)&=\big[A_0x^*(t)-N_1^{-1}B_0^2P(t)\hat{x}^*(t)-N_1^{-1}B_0^2\hat{\phi}^*(t)-N_2^{-1}C_0^2\check{y}(t)\\
                 &\quad-N_2^{-1}C_0^2P(t)\check{\hat{p}}(t)\big]dt+A_1x^*(t)dW_1(t)+A_2x^*(t)dW_2(t)+A_3x^*(t)dW_3(t),\\
-d\hat{\phi}^*(t)&=\Big\{\big[A_0-N_1^{-1}B_0^2P(t)\big]\hat{\phi}^*(t)+A_1\hat{\beta}_1^*(t)+A_3\hat{\beta}_3^*(t)-N_2^{-1}C_0^2P(t)\check{\hat{y}}(t)\\
                 &\quad-N_2^{-1}C_0^2P^2(t)\check{\hat{p}}(t)\Big\}dt-\hat{\beta}_1^*(t)dW_1(t)-\hat{\beta}_3(t)^*dW_3(t),\ t\in[0,T],\\
           x^*(0)&=x_0,\quad \hat{\phi}^*(T)=0.
\end{aligned}
\right.
\end{equation}
Define $X,Y,Z_1,Z_2,Z_3$ as (\ref{new state}) and
\begin{equation*}
\left\{
\begin{aligned}
&\mathcal{A}_0\triangleq\left(\begin{array}{cc}A_0&0\\0&A_0-N_1^{-1}B_0^2P(t)\end{array}\right),\
 \mathcal{A}_1\triangleq\left(\begin{array}{cc}A_1&0\\0&A_1\end{array}\right),\
 \mathcal{A}_2\triangleq\left(\begin{array}{cc}A_2&0\\0&0\end{array}\right),\\
&\mathcal{A}_3\triangleq\left(\begin{array}{cc}A_3&0\\0&A_3\end{array}\right),\ \mathcal{B}_0\triangleq\left(\begin{array}{cc}0&-N_1^{-1}B_0^2\\-N_1^{-1}B_0^2&0\end{array}\right),\
 \overline{\mathcal{B}}_0\triangleq\left(\begin{array}{cc}-N_1^{-1}B_0^2P(t)&0\\0&0\end{array}\right),\ \\
&\mathcal{C}_0\triangleq\left(\begin{array}{cc}-N_2^{-1}C_0&0\\0&0\end{array}\right),\widetilde{\mathcal{C}}_0\triangleq\left(\begin{array}{cc}0&-N_2^{-1}C_0^2P(t)\\0&0\end{array}\right),\
 \widehat{\mathcal{C}}_0\triangleq\left(\begin{array}{cc}0&0\\-N_2^{-1}C_0^2P(t)&0\end{array}\right),\\
&\overline{\mathcal{C}}_0\triangleq\left(\begin{array}{cc}0&0\\0&-N_2^{-1}C_0^2P^2(t)\end{array}\right),\ \mathcal{Q}_2\triangleq\left(\begin{array}{cc}Q_2&0\\0&0\end{array}\right),\
 \mathcal{G}_2\triangleq\left(\begin{array}{cc}G_2&0\\0&0\end{array}\right),\
 X_0\triangleq\left(\begin{array}{c}x_0\\0\end{array}\right),\
\end{aligned}
\right.
\end{equation*}
then we have
\begin{equation}\label{optimal control of the leader-control independent diffusion-2 dim}
      u_2^*(t)=-N_2^{-1}\Big[\left(\begin{array}{cc}C_0&0\end{array}\right)\check{Y}(t)+\left(\begin{array}{cc}C_0P(t)&0\end{array}\right)\check{\hat{X}}(t)\Big],
\end{equation}
and
\begin{equation}\label{optimality system-leader-2 dim-without control-control independent diffusion}
\left\{
\begin{aligned}
  dX(t)&=\big[\mathcal{A}_0X(t)+\overline{\mathcal{B}}_0\hat{X}(t)+\widetilde{\mathcal{C}}_0\check{\hat{X}}(t)+\mathcal{B}_0Y(t)+\mathcal{C}_0\check{Y}(t)\big]dt\\
       &\quad+\mathcal{A}_1X(t)dW_1(t)+\mathcal{A}_2X(t)dW_2(t)+\mathcal{A}_3X(t)dW_3(t),\\
 -dY(t)&=\big[\mathcal{Q}_2X(t)+\mathcal{A}_0Y(t)+\overline{\mathcal{B}}_0\hat{Y}(t)+\widehat{\mathcal{C}}_0\check{\hat{Y}}(t)+\overline{\mathcal{C}}_0\check{\hat{X}}(t)+\mathcal{A}_1Z_1(t)\\
       &\quad+\mathcal{A}_2Z_2(t)+\mathcal{A}_3Z_3(t)\big]dt-Z_1dW_1(t)-Z_2dW_2(t)-Z_3dW_3(t),\ t\in[0,T],\\
   X(0)&=X_0,\quad Y(T)=\mathcal{G}_2X(T),
\end{aligned}
\right.
\end{equation}
where the equation for $\hat{X}(\cdot),\check{X}(\cdot),\check{\hat{X}}(\cdot)$ are
\begin{equation}\label{hat X-control independent diffusion}
\left\{
\begin{aligned}
  d\hat{X}(t)&=\big[\big(\mathcal{A}_0+\overline{\mathcal{B}}_0\big)\hat{X}(t)+\widetilde{\mathcal{C}}_0\check{\hat{X}}(t)+\mathcal{B}_0\hat{Y}(t)+\mathcal{C}_0\check{\hat{Y}}(t)\big]dt\\
             &\quad+\mathcal{A}_1\hat{X}(t)dW_1(t)+\mathcal{A}_3\hat{X}(t)dW_3(t),\ t\in[0,T],\\
   \hat{X}(0)&=X_0,
\end{aligned}
\right.
\end{equation}
\begin{equation}\label{check-X-control independent diffusion}
\left\{
\begin{aligned}
  d\check{X}(t)&=\big[\mathcal{A}_0\check{X}(t)+(\overline{\mathcal{B}}_0+\widetilde{\mathcal{C}}_0)\check{\hat{X}}(t)+(\mathcal{B}_0+\mathcal{C}_0)\check{Y}(t)\big]dt\\
               &\quad+\mathcal{A}_2\check{X}(t)dW_2(t)+\mathcal{A}_3\check{X}(t)dW_3(t),\ t\in[0,T],\\
   \check{X}(0)&=X_0,
\end{aligned}
\right.
\end{equation}
and
\begin{equation}\label{hat and check-X-control independent diffusion}
\left\{
\begin{aligned}
  d\check{\hat{X}}(t)&=\big[(\mathcal{A}_0+\overline{\mathcal{B}}_0+\widetilde{\mathcal{C}}_0)\check{\hat{X}}(t)+(\mathcal{B}_0+\mathcal{C}_0)\check{\hat{Y}}(t)\big]dt+\mathcal{A}_3\check{\hat{X}}(t)dW_3(t),\ t\in[0,T],\\
   \check{\hat{X}}(0)&=X_0,
\end{aligned}
\right.
\end{equation}
respectively.

Define $\mathcal{P}_1(\cdot),\mathcal{P}_2(\cdot),\mathcal{P}_3(\cdot),\mathcal{P}_4(\cdot)$ as (\ref{relation of X and Y}), and apply It\^{o}'s formula to it, we obtain
\begin{equation}\label{Applying Ito's formula to Y-control independent diffusion}
\begin{aligned}
      dY(t)=&\Big\{\big(\dot{\mathcal{P}}_1+\mathcal{P}_1\mathcal{A}_0+\mathcal{P}_1\mathcal{B}_0\mathcal{P}_1\big)X(t)
            +\big[\dot{\mathcal{P}}_2+\mathcal{P}_2(\mathcal{A}_0+\overline{\mathcal{B}}_0)+\mathcal{P}_1\mathcal{B}_0\mathcal{P}_2+\mathcal{P}_2\mathcal{B}_0\mathcal{P}_1\\
           &\ +\mathcal{P}_2\mathcal{B}_0\mathcal{P}_2+\mathcal{P}_1\overline{\mathcal{B}}_0\big]\hat{X}(t)
            +\big[\dot{\mathcal{P}}_3+\mathcal{P}_3\mathcal{A}_0+\mathcal{P}_3(\mathcal{B}_0+\mathcal{C}_0)(\mathcal{P}_1+\mathcal{P}_3)+\mathcal{P}_1\mathcal{B}_0\mathcal{P}_3\\
           &\ +\mathcal{P}_1\mathcal{C}_0\mathcal{P}_3+\mathcal{P}_1\mathcal{C}_0\mathcal{P}_1\big]\check{X}(t)
            +\big[\dot{\mathcal{P}}_4+\mathcal{P}_4(\mathcal{B}_0+\mathcal{C}_0)(\mathcal{P}_1+\mathcal{P}_2+\mathcal{P}_3)\\
           &\ +(\mathcal{P}_1+\mathcal{P}_2+\mathcal{P}_3)(\mathcal{B}_0+\mathcal{C}_0)\mathcal{P}_4
            +\mathcal{P}_4\big(\mathcal{A}_0+\overline{\mathcal{B}}_0+\widetilde{\mathcal{C}}_0\big)+\mathcal{P}_4(\mathcal{B}_0+\mathcal{C}_0)\mathcal{P}_4\\
           &\ +(\mathcal{P}_1+\mathcal{P}_2)\widetilde{\mathcal{C}}_0
            +\mathcal{P}_1\mathcal{C}_0\mathcal{P}_2+\mathcal{P}_2\mathcal{C}_0\mathcal{P}_1+\mathcal{P}_2\mathcal{C}_0\mathcal{P}_2+\mathcal{P}_2\mathcal{B}_0\mathcal{P}_3+\mathcal{P}_3\mathcal{B}_0\mathcal{P}_2\\
           &\ +\mathcal{P}_2\mathcal{C}_0\mathcal{P}_3+\mathcal{P}_3\big(\overline{\mathcal{B}}_0+\widetilde{\mathcal{C}}_0\big)
            +\mathcal{P}_3\mathcal{C}_0\mathcal{P}_2\big]\check{\hat{X}}(t)\Big\}dt\\
           &\ +\big[\mathcal{P}_1\mathcal{A}_1X(t)+\mathcal{P}_2\mathcal{A}_1\hat{X}(t)\big]dW_1(t)+\big[\mathcal{P}_1\mathcal{A}_2X(t)+\mathcal{P}_3\mathcal{A}_2\check{X}(t)\big]dW_2(t)\\
           &\ +\big[\mathcal{P}_1\mathcal{A}_3X(t)+\mathcal{P}_2\mathcal{A}_3\hat{X}(t)+\mathcal{P}_3\mathcal{A}_3\check{X}(t)+\mathcal{P}_4\mathcal{A}_3\check{\hat{X}}(t)\big]dW_3(t)\\
          =&-\Big\{\big(\mathcal{Q}_2+\mathcal{A}_0\mathcal{P}_1\big)X(t)+\big(\mathcal{A}_0\mathcal{P}_2+\overline{\mathcal{B}}_0\mathcal{P}_1+\overline{\mathcal{B}}_0\mathcal{P}_2\big)\hat{X}(t)
            +\mathcal{A}_0\mathcal{P}_3\check{X}(t)\\
           &\quad+\big[\mathcal{A}_0\mathcal{P}_4+\overline{\mathcal{B}}_0\mathcal{P}_3+\overline{\mathcal{B}}_0\mathcal{P}_4+\overline{\mathcal{C}}_0
            +\widehat{\mathcal{C}}_0(\mathcal{P}_1+\mathcal{P}_2+\mathcal{P}_3+\mathcal{P}_4)\big]\check{\hat{X}}(t)+\mathcal{A}_1Z_1(t)\\
           &\quad+\mathcal{A}_2Z_2(t)+\mathcal{A}_3Z_3(t)\Big\}dt+Z_1(t)dW_1(t)+Z_2(t)dW_2(t)+Z_3(t)dW_3(t).
\end{aligned}
\end{equation}
Comparing the diffusion terms on both sides of (\ref{Applying Ito's formula to Y-control independent diffusion}), we directly get
\begin{equation}\label{comparing dW1,dW2,dW3-leader-control independent diffusion}
\left\{
\begin{aligned}
            Z_1(t)&=\mathcal{P}_1\mathcal{A}_1X(t)+\mathcal{P}_2\mathcal{A}_1\hat{X}(t),\
            Z_2(t)=\mathcal{P}_1\mathcal{A}_2X(t)+\mathcal{P}_3\mathcal{A}_2\check{X}(t),\\
            Z_3(t)&=\mathcal{P}_1\mathcal{A}_3X(t)+\mathcal{P}_2\mathcal{A}_3\hat{X}(t)+\mathcal{P}_3\mathcal{A}_3\check{X}(t)+\mathcal{P}_4\mathcal{A}_3\check{\hat{X}}(t).
\end{aligned}
\right.
\end{equation}
It is worth to pointing out that, comparing with the four steps in the control-dependent case of Section 2, the current case is rather simple to obtain (\ref{comparing dW1,dW2,dW3-leader-control independent diffusion}). Comparing the drift term on both sides of (\ref{Applying Ito's formula to Y-control independent diffusion}) and substituting (\ref{comparing dW1,dW2,dW3-leader-control independent diffusion}) into it, we obtain
\begin{equation}\label{system of Riccati equation-control independent diffusion}
\left\{
\begin{aligned}
      &0=\dot{\mathcal{P}}_1+\mathcal{P}_1\mathcal{A}_0+\mathcal{A}_0\mathcal{P}_1
       +\mathcal{A}_1\mathcal{P}_1\mathcal{A}_1+\mathcal{A}_2\mathcal{P}_1\mathcal{A}_2
       +\mathcal{A}_3\mathcal{P}_1\mathcal{A}_3+\mathcal{P}_1\mathcal{B}_0\mathcal{P}_1+\mathcal{Q}_2,\\
      &0=\dot{\mathcal{P}}_2+\mathcal{P}_2(\mathcal{A}_0+\overline{\mathcal{B}}_0)+(\mathcal{A}_0+\overline{\mathcal{B}}_0)\mathcal{P}_2+\mathcal{A}_1\mathcal{P}_2\mathcal{A}_1+\mathcal{A}_3\mathcal{P}_2\mathcal{A}_3
       +\mathcal{P}_1\mathcal{B}_0\mathcal{P}_2+\mathcal{P}_2\mathcal{B}_0\mathcal{P}_1\\
      &\qquad+\mathcal{P}_2\mathcal{B}_0\mathcal{P}_2+\mathcal{P}_1\overline{\mathcal{B}}_0+\overline{\mathcal{B}}_0\mathcal{P}_1,\\
      &0=\dot{\mathcal{P}}_3+\mathcal{P}_3\mathcal{A}_0+\mathcal{A}_0\mathcal{P}_3+\mathcal{A}_2\mathcal{P}_3\mathcal{A}_2+\mathcal{A}_3\mathcal{P}_3\mathcal{A}_3
       +\mathcal{P}_3(\mathcal{B}_0+\mathcal{C}_0)\mathcal{P}_1+\mathcal{P}_1(\mathcal{B}_0+\mathcal{C}_0)\mathcal{P}_3\\
      &\qquad+\mathcal{P}_3(\mathcal{B}_0+\mathcal{C}_0)\mathcal{P}_3+\mathcal{P}_1\mathcal{C}_0\mathcal{P}_1,\\
      &0=\dot{\mathcal{P}}_4+\mathcal{P}_4\big(\mathcal{A}_0+\overline{\mathcal{B}}_0+\widetilde{\mathcal{C}}_0\big)
       +\big(\mathcal{A}_0+\overline{\mathcal{B}}_0+\widehat{\mathcal{C}}_0\big)\mathcal{P}_4+\mathcal{A}_3\mathcal{P}_4\mathcal{A}_3\\
      &\qquad+\mathcal{P}_4(\mathcal{B}_0+\mathcal{C}_0)(\mathcal{P}_1+\mathcal{P}_2+\mathcal{P}_3)+(\mathcal{P}_1+\mathcal{P}_2+\mathcal{P}_3)(\mathcal{B}_0+\mathcal{C}_0)\mathcal{P}_4
       +\mathcal{P}_4(\mathcal{B}_0+\mathcal{C}_0)\mathcal{P}_4\\
      &\qquad+\mathcal{P}_3\big(\overline{\mathcal{B}}_0+\widetilde{\mathcal{C}}_0\big)
       +(\overline{\mathcal{B}}_0+\widehat{\mathcal{C}}_0)\mathcal{P}_3+\mathcal{P}_2\mathcal{B}_0\mathcal{P}_3+\mathcal{P}_3\mathcal{B}_0\mathcal{P}_2+\mathcal{P}_2\mathcal{C}_0\mathcal{P}_3+\mathcal{P}_3\mathcal{C}_0\mathcal{P}_2\\
      &\qquad+(\mathcal{P}_1+\mathcal{P}_2)\widetilde{\mathcal{C}}_0+\widehat{\mathcal{C}}_0(\mathcal{P}_1+\mathcal{P}_2)
       +\mathcal{P}_1\mathcal{C}_0\mathcal{P}_2+\mathcal{P}_2\mathcal{C}_0\mathcal{P}_1+\mathcal{P}_2\mathcal{C}_0\mathcal{P}_2+\overline{\mathcal{C}}_0,\\
      &\mathcal{P}_1(T)=\mathcal{G}_2,\ \mathcal{P}_2(T)=0,\ \mathcal{P}_3(T)=0,\ \mathcal{P}_4(T)=0.
\end{aligned}
\right.
\end{equation}
In this case the solvability of the above system of Riccati equations can be easily obtained, from the standard Riccati equation theory. In fact, noting that the equations for $\mathcal{P}_1(\cdot),\mathcal{P}_2(\cdot),\mathcal{P}_3(\cdot),\\\mathcal{P}_4(\cdot)$ are not coupled. So we can solve firstly $\mathcal{P}_1(\cdot)$, then $\mathcal{P}_2(\cdot)$, thirdly $\mathcal{P}_3(\cdot)$ and finally $\mathcal{P}_4(\cdot)$.

We have the following theorem.

\vspace{1mm}

\noindent{\bf Theorem 3.2}\quad{\it Let $(\mathcal{P}_1(\cdot),\mathcal{P}_2(\cdot),\mathcal{P}_3(\cdot),\mathcal{P}_4(\cdot))$ satisfy (\ref{system of Riccati equation-control independent diffusion}), $\check{\hat{X}}(\cdot)$ be the $\mathcal{G}^1_t\cap\mathcal{G}^2_t$-adapted solution to
\begin{equation}\label{check-hat-X-optimal feedback-control independent diffusion}
\left\{
\begin{aligned}
  d\check{\hat{X}}(t)&=\big[\mathcal{A}_0+\overline{\mathcal{B}}_0+\widetilde{\mathcal{C}}_0+(\mathcal{B}_0+\mathcal{C}_0)(\mathcal{P}_1+\mathcal{P}_2+\mathcal{P}_3+\mathcal{P}_4)\big]\check{\hat{X}}(t)dt\\
                     &\quad+\mathcal{A}_3\check{\hat{X}}(t)dW_3(t),\ t\in[0,T],\\
   \check{\hat{X}}(0)&=X_0,
\end{aligned}
\right.
\end{equation}
$\check{X}(\cdot)$ be the $\mathcal{G}^2_t$-adapted solution to
\begin{equation}\label{check-X-optimal feedback-control independent diffusion}
\left\{
\begin{aligned}
  d\check{X}(t)&=\Big\{\big[\mathcal{A}_0+(\mathcal{B}_0+\mathcal{C}_0)(\mathcal{P}_1+\mathcal{P}_3)\big]\check{X}(t)+\big[\overline{\mathcal{B}}_0+\widetilde{\mathcal{C}}_0
               +(\mathcal{B}_0+\mathcal{C}_0)(\mathcal{P}_1+\mathcal{P}_3)\big]\check{\hat{X}}(t)\Big\}dt\\
               &\quad+\mathcal{A}_2\check{X}(t)dW_2(t)+\mathcal{A}_3\check{X}(t)dW_3(t),\ t\in[0,T],\\
   \check{X}(0)&=X_0,
\end{aligned}
\right.
\end{equation}
$\hat{X}(\cdot)$ be the $\mathcal{G}^1_t$-adapted solution to
\begin{equation}\label{hat-X-optimal feedback-control independent diffusion}
\left\{
\begin{aligned}
  d\hat{X}(t)&=\Big\{\big[\mathcal{A}_0+\overline{\mathcal{B}}_0+\mathcal{B}_0(\mathcal{P}_1+\mathcal{P}_2)\big]\hat{X}(t)
              +\big[\widetilde{\mathcal{C}}_0+\mathcal{B}_0(\mathcal{P}_1+\mathcal{P}_2)+\mathcal{C}_0(\mathcal{P}_1+\mathcal{P}_2\\
             &\qquad+\mathcal{P}_3+\mathcal{P}_4)\big]\check{\hat{X}}(t)\Big\}dt+\mathcal{A}_1\hat{X}(t)dW_1(t)+\mathcal{A}_3\hat{X}(t)dW_3(t),\ t\in[0,T],\\
   \hat{X}(0)&=X_0,
\end{aligned}
\right.
\end{equation}
and $X(\cdot)$ be the $\mathcal{F}_t$-adapted solution to
\begin{equation}\label{close-loop state of the leader-control independent diffusion}
\left\{
\begin{aligned}
  dX(t)&=\Big\{(\mathcal{A}_0+\mathcal{B}_0\mathcal{P}_1)X(t)+(\overline{\mathcal{B}}_0+\mathcal{B}_0\mathcal{P}_2)\hat{X}(t)+\big[\mathcal{B}_0\mathcal{P}_3+\mathcal{C}_0(\mathcal{P}_1+\mathcal{P}_3)\big]\check{X}(t)\\
       &\qquad+\big[\widetilde{\mathcal{C}}_0+\mathcal{B}_0\mathcal{P}_4+\mathcal{C}_0(\mathcal{P}_2+\mathcal{P}_4)\big]\check{\hat{X}}(t)\Big\}dt+\mathcal{A}_1X(t)dW_1(t)\\
       &\quad+\mathcal{A}_2X(t)dW_2(t)+\mathcal{A}_3X(t)dW_3(t),\ t\in[0,T],\\
   X(0)&=X_0,
\end{aligned}
\right.
\end{equation}
and define $(Y(\cdot),Z_1(\cdot),Z_2(\cdot),Z_3(\cdot))$ by (\ref{relation of X and Y}) and (\ref{comparing dW1,dW2,dW3-leader-control independent diffusion}), respectively.
Then (\ref{optimality system-leader-2 dim-without control-control independent diffusion}) holds, and $u_2^*(\cdot)$ given by
\begin{equation}\label{optimal control of the leader-feedback-control independent diffusion-final}
\begin{aligned}
u_2^*(t)&=-N_2^{-1}\left(\begin{array}{cc}C_0P(t)&0\end{array}\right)(\mathcal{P}_1+\mathcal{P}_3)\check{X}(t)\\
        &\quad-N_2^{-1}\Big[\left(\begin{array}{cc}C_0P(t)&0\end{array}\right)(\mathcal{P}_2+\mathcal{P}_4)+\left(\begin{array}{cc}C_0P(t)&0\end{array}\right)\Big]\check{\hat{X}}(t)
\end{aligned}
\end{equation}
is a feedback optimal control of the leader.}

\noindent{\it Proof.}\quad The conclusion is easily obtained from (\ref{relation of X and Y}), (\ref{hat X-control independent diffusion}), (\ref{check-X-control independent diffusion}), (\ref{hat and check-X-control independent diffusion}), (\ref{optimality system-leader-2 dim-without control-control independent diffusion}), (\ref{optimal control of the leader-control independent diffusion-2 dim}). $\Box$

Finally, for the follower, by (\ref{optimal control of the follower-feedback-control independent diffusion}), noting (\ref{new state}) and (\ref{relation of X and Y}), we obtain
\begin{equation}\label{optimal control of the follower-feedback-control independent diffusion-final}
\begin{aligned}
u_1^*(t)=&-N_1^{-1}B_0\big[P(t)\hat{x}^*(t)+\hat{\phi}^*(t)]\\
        =&-N_1^{-1}\Big[\left(\begin{array}{cc}B_0P(t)&0\end{array}\right)\hat{X}(t)+\left(\begin{array}{cc}0&B_0\end{array}\right)\hat{Y}(t)\Big]\\
        =&-N_1^{-1}\Big[\left(\begin{array}{cc}B_0P(t)&0\end{array}\right)+\left(\begin{array}{cc}0&B_0\end{array}\right)\big(\mathcal{P}_1+\mathcal{P}_2\big)\Big]\hat{X}(t)\\
         &-N_1^{-1}\left(\begin{array}{cc}0&B_0\end{array}\right)\big(\mathcal{P}_3+\mathcal{P}_4\big)\check{\hat{X}}(t).
\end{aligned}
\end{equation}

\section{A continuous-time principal-agent problem}

This section is devoted to studying the continuous-time principal-agent problem with overlapping information (Example 1.1 of Section 1), which naturally motivates the research for the problems in previous sections. The financial framework is a generalization of the work by Williams \cite{Wil15}.

In order to apply the results in Section 3, we define $X:=(y,m)^\top$ and then
\begin{equation}\label{state equation--}
\left\{
\begin{aligned}
dX(t)&=\big[\widetilde{r}X(t)+\widetilde{B}e(t)+\alpha_1c(t)+\alpha_2s(t)+\alpha_3d(t)\big]dt+\widetilde{\sigma}_1dW_1(t)\\
     &\quad+\widetilde{\sigma}_2dW_2(t)+\widetilde{\sigma}_3dW_3(t),\ t\in[0,T],\\
 X(0)&=X_0\in\mathbb{R}^2,
\end{aligned}
\right.
\end{equation}
and
\begin{equation}\label{cost principal}
\begin{aligned}
 J_1(e(\cdot),c(\cdot),s(\cdot),d(\cdot))&=\frac{1}{2}\mathbb{E}\left[\int_0^T\big[c^2(t)-e^2(t)+\langle\widetilde{G}_1X(t),X(t)\rangle\big]dt+\langle\widetilde{G}_1X(T),X(T)\rangle\right],\\
 J_2(e(\cdot),c(\cdot),s(\cdot),d(\cdot))&=\frac{1}{2}\mathbb{E}\left[\int_0^T\big[d^2(t)-s^2(t)+\langle\widetilde{G}_2X(t),X(t)\rangle\big]dt+\langle\widetilde{G}_2X(T),X(T)\rangle\right],
\end{aligned}
\end{equation}
where
\begin{equation*}
\left\{
\begin{aligned}
&X_0:=\left(\begin{array}{c}y_0\\m_0\end{array}\right),\ \widetilde{r}:=\left(\begin{array}{cc}r&0\\0&r\end{array}\right),\
 \widetilde{B}:=\left(\begin{array}{c}B\\0\end{array}\right),\ \alpha_1:=\left(\begin{array}{c}0\\-1\end{array}\right),\\
&\alpha_2:=\left(\begin{array}{c}-1\\1\end{array}\right),\ \alpha_3:=\left(\begin{array}{c}-1\\0\end{array}\right),\
 \widetilde{\sigma}_1:=\left(\begin{array}{c}\sigma_1\\\bar{\sigma}_1\end{array}\right),\ \widetilde{\sigma}_2:=\left(\begin{array}{c}\sigma_2\\\bar{\sigma}_2\end{array}\right),\\
&\widetilde{\sigma}_3:=\left(\begin{array}{c}\sigma_3\\\bar{\sigma}_3\end{array}\right),\ \widetilde{G}_1:=\left(\begin{array}{cc}0&0\\0&1\end{array}\right),\
 \widetilde{G}_2:=\left(\begin{array}{cc}1&0\\0&0\end{array}\right).
\end{aligned}
\right.
\end{equation*}

For the follower (agent)'s problem, first the leader (principal) announces his control $s(\cdot),d(\cdot)$. Following the step in Section 3.1, we have
\begin{equation}\label{optimal effort and consumption}
 e^*(t)=\widetilde{B}^\top\big[P(t)\hat{X}(t)+\Phi(t)\big],\ c^*(t)=-\alpha_1^\top\big[P(t)\hat{X}(t)+\Phi(t)\big],
\end{equation}
where $2\times2$-matrix-valued function $P(\cdot)$ satisfies
\begin{equation}\label{Riccati'type equation}
\left\{
\begin{aligned}
 &\dot{P}(t)+P(t)\widetilde{r}+\widetilde{r}^\top P(t)+P(t)\big(\widetilde{B}\widetilde{B}^\top-\alpha_1\alpha_1^\top\big)P(t)^\top+\widetilde{G}_1=0,\ t\in[0,T],\\
 &P(T)=\widetilde{G}_1,
\end{aligned}
\right.
\end{equation}
and $\mathbb{R}^2$-valued, $\mathcal{G}^1_t$-adapted process quadruple $(\hat{X}(\cdot),\Phi(\cdot),\Pi_1(\cdot),\Pi_3(\cdot))$ satisfies FBSDFE
\begin{equation}\label{FBSDFE}
\left\{
\begin{aligned}
    d\hat{X}(t)&=\Big\{\big[\widetilde{r}+\widetilde{B}\widetilde{B}^\top P(t)-\alpha_1\alpha_1^\top P(t)\big]\hat{X}(t)+\big(\widetilde{B}\widetilde{B}^\top-\alpha_1\alpha_1^\top\big)\Phi(t)\\
               &\qquad+\alpha_2\hat{s}(t)+\alpha_3\hat{d}(t)\Big\}dt+\widetilde{\sigma}_1dW_1(t)+\widetilde{\sigma}_3dW_3(t),\\
      -d\Phi(t)&=\Big\{\big[\widetilde{r}+\widetilde{B}\widetilde{B}^\top P(t)-\alpha_1\alpha_1^\top P(t)\big]\Phi(t)+P(t)\alpha_2\hat{s}(t)+P(t)\alpha_3\hat{d}(t)\Big\}dt\\
               &\quad-\Pi_1(t)dW_1(t)-\Pi_3(t)dW_3(t),\ t\in[0,T],\\
     \hat{X}(0)&=X_0,\quad \Phi(T)=0.
\end{aligned}
\right.
\end{equation}

For the leader (principal)'s problem, the state now writes
\begin{equation}\label{state of the principal}
\left\{
\begin{aligned}
    dX(t)&=\Big\{\widetilde{r}X(t)+\big(\widetilde{B}\widetilde{B}^\top-\alpha_1\alpha_1^\top\big)P(t)\hat{X}(t)+\big(\widetilde{B}\widetilde{B}^\top-\alpha_1\alpha_1^\top\big)\Phi(t)+\alpha_2s(t)\\
         &\qquad+\alpha_3d(t)\Big\}dt+\widetilde{\sigma}_1dW_1(t)+\widetilde{\sigma}_2dW_2(t)+\widetilde{\sigma}_3dW_3(t),\\
-d\Phi(t)&=\Big\{\big[\widetilde{r}+\widetilde{B}\widetilde{B}^\top P(t)-\alpha_1\alpha_1^\top P(t)\big]\Phi(t)+P(t)\alpha_2\hat{s}(t)+P(t)\alpha_3\hat{d}(t)\Big\}dt\\
         &\quad-\Pi_1(t)dW_1(t)-\Pi_3(t)dW_3(t),\ t\in[0,T],\\
     X(0)&=X_0,\quad \Phi(T)=0.
\end{aligned}
\right.
\end{equation}
Following the step in Section 3.2, we have
\begin{equation}\label{optimal payment and consumption}
 s^*(t)=-\alpha_2\check{\widetilde{y}}(t)-\alpha_2P(t)\check{\hat{\widetilde{y}}}(t),\ d^*(t)=\alpha_3\check{\widetilde{y}}(t)+\alpha_3P(t)\check{\hat{\widetilde{y}}}(t),
\end{equation}
where $\mathbb{R}^2$-valued, $\mathcal{F}_t$-adapted process quintuple $(\widetilde{p}(\cdot),\widetilde{y}(\cdot),z_1(\cdot),z_2(\cdot),z_3(\cdot))$ satisfies the adjoint equation
\begin{equation}\label{adjoint equation of the principal}
\left\{
\begin{aligned}
 d\widetilde{p}(t)&=\Big\{\big(\widetilde{B}\widetilde{B}^\top-\alpha_1\alpha_1^\top\big)\widetilde{y}(t)+\big[\widetilde{r}^\top+\widetilde{B}\widetilde{B}^\top P(t)+\alpha_1\alpha_1^\top P(t)\big]\widetilde{p}(t)\Big\}dt,\\
-d\widetilde{y}(t)&=\big[\widetilde{r}\widetilde{y}(t)+\big(\widetilde{B}\widetilde{B}^\top-\alpha_1\alpha_1^\top\big)P(t)\hat{\widetilde{y}}(t)+\widetilde{G}_2X^*(t)\big]dt-z_1(t)dW_1(t)\\
                  &\quad-z_2(t)dW_2(t)-z_3(t)dW_3(t),\ t\in[0,T],\\
  \widetilde{p}(0)&=0,\quad \widetilde{y}(T)=\widetilde{G}_2X^*(T).
\end{aligned}
\right.
\end{equation}

Let
\begin{equation}\label{stack}
\begin{aligned}
\mathcal{X}=\left(\begin{array}{c}X^*\\\widetilde{p}\end{array}\right),\ Y=\left(\begin{array}{c}\widetilde{y}\\\Phi^*\end{array}\right),\ Z_1=\left(\begin{array}{c}z_1\\\Pi_1^*\end{array}\right),\
Z_2=\left(\begin{array}{c}z_2\\0\end{array}\right),\ Z_3=\left(\begin{array}{c}z_3\\\Pi_3^*\end{array}\right),
\end{aligned}
\end{equation}
and
\begin{equation*}
\left\{
\begin{aligned}
&\mathcal{A}_0:=\left(\begin{array}{cc}\widetilde{r}&0\\0&\widetilde{r}+\widetilde{B}\widetilde{B}^\top P(t)-\alpha_1\alpha_1^\top P(t)\end{array}\right),\
 \overline{\mathcal{A}}_0:=\left(\begin{array}{cc}\widetilde{B}\widetilde{B}^\top P(t)-\alpha_1\alpha_1^\top P(t)&0\\0&0\end{array}\right),\\
&\mathcal{B}_0:=\left(\begin{array}{cc}0&\widetilde{B}\widetilde{B}^\top-\alpha_1\alpha_1^\top\\\widetilde{B}\widetilde{B}^\top-\alpha_1\alpha_1^\top&0\end{array}\right),\
 \widetilde{\alpha}_2:=\left(\begin{array}{c}\alpha_2\\0\end{array}\right),\ \widetilde{\alpha}_3:=\left(\begin{array}{c}\alpha_3\\0\end{array}\right),\\
&\Sigma_1:=\left(\begin{array}{c}\widetilde{\sigma}_1\\0\end{array}\right),\ \Sigma_2:=\left(\begin{array}{c}\widetilde{\sigma}_2\\0\end{array}\right),\
 \Sigma_3:=\left(\begin{array}{c}\widetilde{\sigma}_3\\0\end{array}\right),\ \mathcal{G}_2:=\left(\begin{array}{cc}\widetilde{G}_2&0\\0&0\end{array}\right),\ \mathcal{X}_0:=\left(\begin{array}{c}X_0\\0\end{array}\right),\\
&\overline{\alpha}_2:=\left(\begin{array}{c}0\\P(t)\alpha_2\end{array}\right),\ \overline{\alpha}_3:=\left(\begin{array}{c}0\\P(t)\alpha_3\end{array}\right),\
 \Lambda_1:=\left(\begin{array}{c}Z_1\\\beta_1\end{array}\right),\ \Lambda_2:=\left(\begin{array}{c}Z_2\\0\end{array}\right),\
 \Lambda_3:=\left(\begin{array}{c}Z_3\\\beta_3\end{array}\right),
\end{aligned}
\right.
\end{equation*}
then we have
\begin{equation}\label{X-X}
\left\{
\begin{aligned}
d\mathcal{X}(t)&=\big[\mathcal{A}_0\mathcal{X}(t)+\overline{\mathcal{A}}_0\hat{\mathcal{X}}(t)+\mathcal{B}_0Y(t)+\widetilde{\alpha}_2s^*(t)+\widetilde{\alpha}_3d^*(t)\big]dt\\
               &\quad+\Sigma_1dW_1(t)+\Sigma_2dW_2(t)+\Sigma_3dW_3(t),\\
         -dY(t)&=\Big\{\big[\mathcal{G}_2\mathcal{X}(t)+\mathcal{A}_0Y(t)+\overline{\mathcal{A}}_0\hat{Y}(t)++\overline{\alpha}_2\hat{s}^*(t)+\overline{\alpha}_3\hat{d}^*(t)\Big\}dt\\
               &\quad-\Pi_1(t)dW_1(t)-\Pi_3(t)dW_3(t),\ t\in[0,T],\\
 \mathcal{X}(0)&=\mathcal{X}_0,\quad Y(T)=\mathcal{G}_2\mathcal{X}(T).
\end{aligned}
\right.
\end{equation}
As Section 3.2, letting $Y(t)=\mathcal{P}_1(t)\mathcal{X}(t)+\mathcal{P}_2(t)\hat{\mathcal{X}}(t)+\mathcal{P}_3(t)\check{\mathcal{X}}(t)+\mathcal{P}_4(t)\check{\hat{\mathcal{X}}}(t)$, then we get
\begin{equation}\label{optimal payment and consumption-2}
\left\{
\begin{aligned}
 s^*(t)&=-\big(\widetilde{\alpha}_2^\top\mathcal{P}_1+\widetilde{\alpha}_2^\top\mathcal{P}_3\big)\check{\mathcal{X}}(t)
        -\big(\widetilde{\alpha}_2^\top\mathcal{P}_2+\widetilde{\alpha}_2^\top\mathcal{P}_4+\overline{\alpha}_2^\top\big)\check{\hat{\mathcal{X}}}(t),\\
 d^*(t)&=\big(\widetilde{\alpha}_3^\top\mathcal{P}_1+\widetilde{\alpha}_3^\top\mathcal{P}_3\big)\check{\mathcal{X}}(t)
        +\big(\widetilde{\alpha}_3^\top\mathcal{P}_2+\widetilde{\alpha}_3^\top\mathcal{P}_4+\overline{\alpha}_3^\top\big)\check{\hat{\mathcal{X}}}(t),
\end{aligned}
\right.
\end{equation}
where $4\times4$-matrix-valued functions $\mathcal{P}_1(\cdot),\mathcal{P}_2(\cdot),\mathcal{P}_3(\cdot),\mathcal{P}_4(\cdot)$ satisfy
\begin{equation}\label{Riccati equation system}
\left\{
\begin{aligned}
      &0=\dot{\mathcal{P}}_1+\mathcal{P}_1\mathcal{A}_0+\mathcal{A}_0^\top\mathcal{P}_1+\mathcal{P}_1\mathcal{B}_0\mathcal{P}_1+\mathcal{G}_2,\ \mathcal{P}_1(T)=\mathcal{G}_2,\\
      &0=\dot{\mathcal{P}}_2+\mathcal{P}_2(\mathcal{A}_0+\overline{\mathcal{A}}_0)+(\mathcal{A}_0+\overline{\mathcal{A}}_0)^\top\mathcal{P}_2
       +\mathcal{P}_1\mathcal{B}_0\mathcal{P}_2+\mathcal{P}_2\mathcal{B}_0\mathcal{P}_1+\mathcal{P}_2\mathcal{B}_0\mathcal{P}_2\\
      &\quad\ +\mathcal{P}_1\overline{\mathcal{A}}_0+\overline{\mathcal{A}}_0^\top\mathcal{P}_1,\ \mathcal{P}_2(T)=0,\\
      &0=\dot{\mathcal{P}}_3+\mathcal{A}_0^\top\mathcal{P}_3+\mathcal{P}_3\mathcal{A}_0
       +\mathcal{P}_1\big(\mathcal{B}_0+\widetilde{\alpha}_3\widetilde{\alpha}_3^\top-\widetilde{\alpha}_2\widetilde{\alpha}_2^\top\big)\mathcal{P}_3
       +\mathcal{P}_3\big(\mathcal{B}_0+\widetilde{\alpha}_3\widetilde{\alpha}_3^\top-\widetilde{\alpha}_2\widetilde{\alpha}_2^\top\big)\mathcal{P}_1\\
      &\quad\ +\mathcal{P}_3\big(\mathcal{B}_0+\widetilde{\alpha}_3\widetilde{\alpha}_3^\top-\widetilde{\alpha}_2\widetilde{\alpha}_2^\top\big)\mathcal{P}_3
       +\mathcal{P}_1\big(\widetilde{\alpha}_3\widetilde{\alpha}_3^\top-\widetilde{\alpha}_2\widetilde{\alpha}_2^\top\big)\mathcal{P}_1,\ \mathcal{P}_3(T)=0,\\
      &0=\dot{\mathcal{P}}_4+\mathcal{P}_4(\mathcal{A}_0+\overline{\mathcal{A}}_0)+(\mathcal{A}_0+\overline{\mathcal{A}}_0)^\top\mathcal{P}_4
       +\mathcal{P}_2\big(\widetilde{\alpha}_2\overline{\alpha}_2^\top-\widetilde{\alpha}_3\overline{\alpha}_3^\top\big)\\
      &\quad\ +\mathcal{P}_4\big(\mathcal{P}_1\mathcal{B}_0-\mathcal{P}_1\widetilde{\alpha}_2\widetilde{\alpha}_2^\top+\mathcal{P}_1\widetilde{\alpha}_3\widetilde{\alpha}_3^\top
       +\mathcal{P}_2\mathcal{B}_0-\mathcal{P}_2\widetilde{\alpha}_2\widetilde{\alpha}_2^\top+\mathcal{P}_2\widetilde{\alpha}_3\widetilde{\alpha}_3^\top
       -\overline{\alpha}_2\widetilde{\alpha}_2^\top+\overline{\alpha}_3\widetilde{\alpha}_3^\top\big)\\
      &\quad\ +\big(\mathcal{B}_0\mathcal{P}_1-\widetilde{\alpha}_2\widetilde{\alpha}_2^\top\mathcal{P}_1+\widetilde{\alpha}_3\widetilde{\alpha}_3^\top\mathcal{P}_1
       +\mathcal{B}_0\mathcal{P}_2-\widetilde{\alpha}_2\widetilde{\alpha}_2^\top\mathcal{P}_2+\widetilde{\alpha}_3\widetilde{\alpha}_3^\top\mathcal{P}_2
       -\widetilde{\alpha}_2\overline{\alpha}_2^\top+\widetilde{\alpha}_3\overline{\alpha}_3^\top\big)\mathcal{P}_4\\
      &\quad\ +\mathcal{P}_4\big(\mathcal{B}_0-\widetilde{\alpha}_2\widetilde{\alpha}_2^\top+\widetilde{\alpha}_3\widetilde{\alpha}_3^\top\big)\mathcal{P}_4
       -\mathcal{P}_1\widetilde{\alpha}_2\widetilde{\alpha}_2^\top\mathcal{P}_2+\mathcal{P}_1\widetilde{\alpha}_3\widetilde{\alpha}_3^\top\mathcal{P}_2
       -\mathcal{P}_2\widetilde{\alpha}_2\widetilde{\alpha}_2^\top\mathcal{P}_1\\
      &\quad\ +\mathcal{P}_2\widetilde{\alpha}_3\widetilde{\alpha}_3^\top\mathcal{P}_1
       +\mathcal{P}_2\big(\mathcal{B}_0-\widetilde{\alpha}_2\widetilde{\alpha}_2^\top-\widetilde{\alpha}_3\widetilde{\alpha}_3^\top\big)\mathcal{P}_3
       -\mathcal{P}_2\widetilde{\alpha}_2\widetilde{\alpha}_2^\top\mathcal{P}_2+\mathcal{P}_2\widetilde{\alpha}_3\widetilde{\alpha}_3^\top\mathcal{P}_2\\
      &\quad\ -\mathcal{P}_1\widetilde{\alpha}_2\overline{\alpha}_2^\top+\mathcal{P}_1\widetilde{\alpha}_3\overline{\alpha}_3^\top
       -\mathcal{P}_2\widetilde{\alpha}_2\overline{\alpha}_2^\top+\mathcal{P}_2\widetilde{\alpha}_3\overline{\alpha}_3^\top
       -\overline{\alpha}_2\widetilde{\alpha}_2^\top\mathcal{P}_1+\overline{\alpha}_3\widetilde{\alpha}_3^\top\mathcal{P}_1
       -\overline{\alpha}_2\widetilde{\alpha}_2^\top\mathcal{P}_2\\
      &\quad\ +\overline{\alpha}_3\widetilde{\alpha}_3^\top\mathcal{P}_2-\overline{\alpha}_2\widetilde{\alpha}_2^\top\mathcal{P}_3+\overline{\alpha}_3\widetilde{\alpha}_3^\top\mathcal{P}_3
       +\overline{\mathcal{A}}_0^\top\mathcal{P}_3-\overline{\alpha}_2\overline{\alpha}_2^\top+\overline{\alpha}_3\overline{\alpha}_3^\top,\ \mathcal{P}_4(T)=0,
\end{aligned}
\right.
\end{equation}
$\mathbb{R}^4$-valued, $\mathcal{G}^2_t$-adapted processes $\check{\mathcal{X}}(\cdot)$ satisfies the SDFE
\begin{equation}\label{check XX}
\left\{
\begin{aligned}
d\check{\mathcal{X}}(t)&=\Big\{\big[\mathcal{A}_0+\big(\mathcal{B}_0-\widetilde{\alpha}_2\widetilde{\alpha}_2^\top
                        +\widetilde{\alpha}_3\widetilde{\alpha}_3^\top\big)\big(\mathcal{P}_1+\mathcal{P}_3\big)\big]\check{\mathcal{X}}(t)\\
                       &\qquad+\big[\overline{\mathcal{A}}_0
                        +\big(\mathcal{B}_0-\widetilde{\alpha}_2\widetilde{\alpha}_2^\top+\widetilde{\alpha}_3\widetilde{\alpha}_3^\top\big)\big(\mathcal{P}_2+\mathcal{P}_4\big)
                        -\widetilde{\alpha}_2\overline{\alpha}_2^\top+\widetilde{\alpha}_3\overline{\alpha}_3^\top\big]\check{\hat{\mathcal{X}}}(t)\Big\}dt\\
                       &\quad+\Sigma_2dW_2(t)+\Sigma_3dW_3(t),\ t\in[0,T],\\
 \check{\mathcal{X}}(0)&=\mathcal{X}_0,
\end{aligned}
\right.
\end{equation}
and $\mathbb{R}^4$-valued, $\mathcal{G}^2_t\cap\mathcal{G}^2_t$-adapted processes $\check{\hat{\mathcal{X}}}(\cdot)$ satisfies
\begin{equation}\label{check hat XX}
\left\{
\begin{aligned}
d\check{\hat{\mathcal{X}}}(t)&=\big[\mathcal{A}_0+\overline{\mathcal{A}}_0+\big(\mathcal{B}_0-\widetilde{\alpha}_2\widetilde{\alpha}_2^\top
                              +\widetilde{\alpha}_3\widetilde{\alpha}_3^\top\big)\big(\mathcal{P}_1+\mathcal{P}_2+\mathcal{P}_3+\mathcal{P}_4\big)\\
                             &\quad-\widetilde{\alpha}_2\overline{\alpha}_2^\top+\widetilde{\alpha}_3\overline{\alpha}_3^\top\big]\check{\hat{\mathcal{X}}}(t)dt+\Sigma_3dW_3(t),\ t\in[0,T],\\
 \check{\hat{\mathcal{X}}}(0)&=\mathcal{X}_0.
\end{aligned}
\right.
\end{equation}

For the follower, by (\ref{optimal effort and consumption}), noting (\ref{stack}), we obtain
\begin{equation}\label{optimal effort and consumption-agent}
\left\{
\begin{aligned}
e^*(t)&=\widetilde{B}^\top\big[P(t)\hat{X}^*(t)+\Phi^*(t)\big]=\Big[\left(\begin{array}{cc}\widetilde{B}^\top P(t)&0\end{array}\right)\hat{\mathcal{X}}(t)
          +\left(\begin{array}{cc}0&\widetilde{B}^\top\end{array}\right)\hat{Y}(t)\Big]\\
     &=\Big[\left(\begin{array}{cc}\widetilde{B}^\top P(t)&0\end{array}\right)+\left(\begin{array}{cc}0&\widetilde{B}^\top\end{array}\right)\big(\mathcal{P}_1+\mathcal{P}_2\big)\Big]\hat{\mathcal{X}}(t)
         +\left(\begin{array}{cc}0&\widetilde{B}^\top\end{array}\right)\big(\mathcal{P}_3+\mathcal{P}_4\big)\check{\hat{\mathcal{X}}}(t),\\
c^*(t)&=-\alpha_1^\top\big[P(t)\hat{X}^*(t)+\Phi^*(t)\big]=\Big[\left(\begin{array}{cc}-\alpha_1^\top P(t)&0\end{array}\right)\hat{\mathcal{X}}(t)
          +\left(\begin{array}{cc}0&-\alpha_1^\top\end{array}\right)\hat{Y}(t)\Big]\\
     &=\Big[\left(\begin{array}{cc}-\alpha_1^\top P(t)&0\end{array}\right)+\left(\begin{array}{cc}0&-\alpha_1^\top\end{array}\right)\big(\mathcal{P}_1+\mathcal{P}_2\big)\Big]\hat{\mathcal{X}}(t)
         +\left(\begin{array}{cc}0&-\alpha_1^\top\end{array}\right)\big(\mathcal{P}_3+\mathcal{P}_4\big)\check{\hat{\mathcal{X}}}(t),
\end{aligned}
\right.
\end{equation}
where $\mathbb{R}^4$-valued, $\mathcal{G}^1_t$-adapted process $\hat{\mathcal{X}}(\cdot)$ satisfies the SDFE
\begin{equation}\label{hat XX}
\left\{
\begin{aligned}
d\hat{\mathcal{X}}(t)&=\Big\{\big[\mathcal{A}_0+\overline{\mathcal{A}}_0+\mathcal{B}_0\big(\mathcal{P}_1+\mathcal{P}_2\big)\big]\hat{\mathcal{X}}(t)
                      +\big[\mathcal{B}_0\big(\mathcal{P}_3+\mathcal{P}_4\big)\\
                     &\qquad-\big(\widetilde{\alpha}_2\widetilde{\alpha}_2^\top-\widetilde{\alpha}_3\widetilde{\alpha}_3^\top\big)\big(\mathcal{P}_1+\mathcal{P}_2+\mathcal{P}_3+\mathcal{P}_4\big)\big]\check{\hat{\mathcal{X}}}(t)
                      -\widetilde{\alpha}_2\overline{\alpha}_2^\top+\widetilde{\alpha}_3\overline{\alpha}_3^\top\Big\}dt\\
                     &\quad+\Sigma_1dW_1(t)+\Sigma_3dW_3(t),\ t\in[0,T],\\
 \hat{\mathcal{X}}(0)&=\mathcal{X}_0.
\end{aligned}
\right.
\end{equation}
And the optimal state equation of the leader is
\begin{equation}\label{XX}
\left\{
\begin{aligned}
d\mathcal{X}(t)&=\Big\{\big(\mathcal{A}_0+\mathcal{B}_0\mathcal{P}_1\big)\mathcal{X}(t)+\big(\overline{\mathcal{A}}_0+\mathcal{B}_0\mathcal{P}_2\big)\hat{\mathcal{X}}(t)
                +\big[\big(\mathcal{B}_0-\widetilde{\alpha}_2\widetilde{\alpha}_2^\top+\widetilde{\alpha}_3\widetilde{\alpha}_3^\top\big)\mathcal{P}_3\\
               &\qquad-\big(\widetilde{\alpha}_2\widetilde{\alpha}_2^\top-\widetilde{\alpha}_3\widetilde{\alpha}_3^\top\big)\mathcal{P}_1\big]\check{\mathcal{X}}(t)
                +\big[\big(\mathcal{B}_0+\widetilde{\alpha}_2\widetilde{\alpha}_2^\top-\widetilde{\alpha}_3\widetilde{\alpha}_3^\top\big)\mathcal{P}_4\\
               &\qquad-\big(\widetilde{\alpha}_2\widetilde{\alpha}_2^\top-\widetilde{\alpha}_3\widetilde{\alpha}_3^\top\big)\mathcal{P}_2-\widetilde{\alpha}_2\overline{\alpha}_2^\top
               +\widetilde{\alpha}_3\overline{\alpha}_3^\top\big]\check{\hat{\mathcal{X}}}(t)\Big\}dt\\
               &\quad+\Sigma_1dW_1(t)+\Sigma_2dW_2(t)+\Sigma_3dW_3(t),\ t\in[0,T],\\
 \mathcal{X}(0)&=\mathcal{X}_0.
\end{aligned}
\right.
\end{equation}

Finally, we rewrite the above Stackelberg equilibrium strategy $(e^*(\cdot),c^*(\cdot),s^*(\cdot),d^*(\cdot))$, with respect to the asset of the principal $y(\cdot)$ and that of the agent $m(\cdot)$. In fact, let
$$\mathcal{P}_k\equiv\big(\mathcal{P}_k^{i,j}\big)_{4\times4},k=1,2,3,4,\mbox{ and }\quad P\equiv\big(P^{i,j}\big)_{2\times2},$$ where $\mathcal{P}^{i,j}$, $P^{i,j}$ denotes the elements of the matrices.

Then by (\ref{optimal payment and consumption-2}) and (\ref{optimal effort and consumption-agent}) we have
\begin{equation*}
\begin{aligned}
 s^*(t)&=\big(-\mathcal{P}_1^{1,1}+\mathcal{P}_1^{2,1}-\mathcal{P}_3^{1,1}+\mathcal{P}_3^{2,1}\big)\check{y}(t)
       +\big(-\mathcal{P}_2^{1,1}+\mathcal{P}_2^{2,1}-\mathcal{P}_4^{1,1}+\mathcal{P}_4^{2,1}\big)\check{\hat{y}}(t)\\
      &\quad+\big(-\mathcal{P}_1^{1,2}+\mathcal{P}_1^{2,2}-\mathcal{P}_3^{1,2}+\mathcal{P}_3^{2,2}\big)\check{m}(t)
       +\big(-\mathcal{P}_2^{1,2}+\mathcal{P}_2^{2,2}-\mathcal{P}_4^{1,2}+\mathcal{P}_4^{2,2}\big)\check{\hat{m}}(t)\\
      &\quad+\left(\begin{array}{c}-\mathcal{P}_1^{1,3}+\mathcal{P}_1^{2,3}-\mathcal{P}_3^{1,3}+\mathcal{P}_3^{2,3}\\
       -\mathcal{P}_1^{1,4}+\mathcal{P}_1^{2,4}-\mathcal{P}_3^{1,4}+\mathcal{P}_3^{2,4}\\
       -\mathcal{P}_2^{1,3}+\mathcal{P}_2^{2,3}-\mathcal{P}_4^{1,3}+\mathcal{P}_4^{2,3}-P^{1,1}+P^{1,2}\\
       -\mathcal{P}_2^{1,4}+\mathcal{P}_2^{2,4}-\mathcal{P}_4^{1,4}+\mathcal{P}_4^{2,4}-P^{2,1}+P^{2,2}\end{array}\right)^\top
       \left(\begin{array}{c}\check{\widetilde{p}}(t)\\\check{\hat{\widetilde{p}}}(t)\end{array}\right),\\
\end{aligned}
\end{equation*}
\begin{equation}\label{Stackelberg equalibrium}
\begin{aligned}
 d^*(t)&=-\big(\mathcal{P}_1^{1,1}+\mathcal{P}_3^{1,1}\big)\check{y}(t)-\big(\mathcal{P}_2^{1,1}+\mathcal{P}_4^{1,1}\big)\check{\hat{y}}(t)
       -\big(\mathcal{P}_1^{1,2}+\mathcal{P}_3^{1,2}\big)\check{m}(t)\\
      &\quad-\big(\mathcal{P}_2^{1,2}+\mathcal{P}_4^{1,2}\big)\check{\hat{m}}(t)-\left(\begin{array}{c}\mathcal{P}_1^{1,3}+\mathcal{P}_3^{1,3}\\
       \mathcal{P}_1^{1,4}+\mathcal{P}_3^{1,4}\\\mathcal{P}_2^{1,3}+\mathcal{P}_4^{1,3}+P^{1,1}\\\mathcal{P}_2^{1,4}+\mathcal{P}_4^{1,4}+P^{2,1}\end{array}\right)^\top
       \left(\begin{array}{c}\check{\widetilde{p}}(t)\\\check{\hat{\widetilde{p}}}(t)\end{array}\right),\\
e^*(t)&=B\big(P^{1,1}+\mathcal{P}_1^{3,1}+\mathcal{P}_2^{3,1}\big)\hat{y}(t)+B\big(\mathcal{P}_3^{3,1}+\mathcal{P}_4^{3,1}\big)\check{\hat{y}}(t)\\
      &\quad+B\big(P^{1,2}+\mathcal{P}_1^{3,2}+\mathcal{P}_2^{3,2}\big)\hat{m}(t)+B\big(\mathcal{P}_3^{3,2}+\mathcal{P}_4^{3,2}\big)\check{\hat{m}}(t)\\
      &\quad+B\left(\begin{array}{cccc}\mathcal{P}_1^{3,3}+\mathcal{P}_2^{3,3}&
       \mathcal{P}_1^{3,4}+\mathcal{P}_2^{3,4}&\mathcal{P}_3^{3,3}+\mathcal{P}_4^{3,3}&
       \mathcal{P}_3^{3,4}+\mathcal{P}_4^{3,4}\end{array}\right)
       \left(\begin{array}{c}\hat{\widetilde{p}}(t)\\\check{\hat{\widetilde{p}}}(t)\end{array}\right),\\
c^*(t)&=\big(P^{2,1}+\mathcal{P}_1^{4,1}+\mathcal{P}_2^{4,1}\big)\hat{y}(t)+\big(\mathcal{P}_3^{4,1}+\mathcal{P}_4^{4,1}\big)\check{\hat{y}}(t)\\
      &\quad+\big(P^{2,2}+\mathcal{P}_1^{4,2}+\mathcal{P}_2^{4,2}\big)\hat{m}(t)+\big(\mathcal{P}_3^{4,2}+\mathcal{P}_4^{4,2}\big)\check{\hat{m}}(t)\\
      &\quad+\left(\begin{array}{cccc}\mathcal{P}_1^{4,3}+\mathcal{P}_2^{4,3}&
       \mathcal{P}_1^{4,4}+\mathcal{P}_2^{4,4}&\mathcal{P}_3^{4,3}+\mathcal{P}_4^{4,3}&
       \mathcal{P}_3^{4,4}+\mathcal{P}_4^{4,4}\end{array}\right)
       \left(\begin{array}{c}\hat{\widetilde{p}}(t)\\\check{\hat{\widetilde{p}}}(t)\end{array}\right),
\end{aligned}
\end{equation}
where $(\check{\hat{y}}(t),\check{\hat{m}}(t),\check{\hat{\widetilde{p}}}(t))\equiv\check{\hat{\mathcal{X}}}(t)$ satisfies (\ref{check hat XX}), $(\hat{y}(t),\hat{m}(t),\hat{\widetilde{p}}(t))\equiv\hat{\mathcal{X}}(t)$ satisfies (\ref{hat XX}) and $(\check{y}(t),\check{m}(t),\check{\widetilde{p}}(t))\equiv\check{\mathcal{X}}(t)$ satisfies (\ref{check XX}).

\section{Concluding remarks}

In this paper, we have discussed the stochastic LQ Stackelberg differential game with overlapping information. This kind of game problem possesses two attractive features worthy of being highlighted. First, the game problem has the asymmetric and overlapping information between the two players, which was not considered in Yong \cite{Yong02}, \O ksendal et al. \cite{OSU13} and Bensoussan et al. \cite{BCS15}.  Stochastic filtering technique is introduced to compute the optimal filtering estimates for the corresponding adjoint processes, which perform as the solution to some nonstandard stochastic filtering equations. Second, the Stackelberg equilibrium is represented in its state estimate feedback form, under some appropriate assumptions on the coefficient matrices in the state equation and the cost functional. Some new system of Riccati equations are first introduced in this paper, to deal with the leader's problem.

Many interesting and important problems remain open. The solvability of the system of Riccati equations (\ref{system of Riccati equation}) is a challenging question. It is worthy to study the numerical approximation
of its solution (\cite{DM97}, \cite{YZ99}). Problems with time delay (\cite{XZ16}), with partial observation (Friedman \cite{Fr72}, Wang et al. \cite{WWX13}, \cite{SWX16}) and of mean-field type (Wang et al. \cite{WZZ14}) which are important and reasonable for applications and more technological demanding in its filtering procedure, are highly desirable for further research. These challenging topics will be considered in our future work.

\end{document}